\documentclass[10pt]{article}
\usepackage[affil-it]{authblk}

\usepackage{amsmath, amsthm, amssymb}
\usepackage{amssymb}
\usepackage{graphicx}
\usepackage{subfigure}
\usepackage[noprefix]{nomencl}

\newcommand\bi{\begin{itemize}}
\newcommand\ei{\end{itemize}}

\def\addsymbol #1: #2#3{$#1$ \> \parbox{5in}{#2 \dotfill \pageref{#3}}\\}

\newtheorem{fact}{Fact}

\newtheorem{defin}{Definition}[section]
          
          \newtheorem{teo}{Theorem}[section]
          
          \newtheorem{con}{Conjecture}
          \newtheorem{cond}{Condition}
          \newtheorem{prop}[teo]{Proposition}
          \newtheorem{lem}{Lemma}[section]
          \newtheorem{rmk}[teo]{Remark}
          \newtheorem{cor}{Corollary}[section]
          
          \newcommand{\bfact}{\begin{fact}}
          \newcommand{\efact}{\end{fact}}
          \newcommand{\Netbk}{{\cal N}^{b,k}}

          \newcommand{\beq}{\begin{equation}}
          \newcommand{\eeq}{\end{equation}}
          \newcommand{\beqn}{\begin{eqnarray}}
          \newcommand{\beqnn}{\begin{eqnarray*}}
          \newcommand{\eeqn}{\end{eqnarray}}
          \newcommand{\eeqnn}{\end{eqnarray*}}
          \newcommand{\bprop}{\begin{prop}}
          \newcommand{\eprop}{\end{prop}}
          \newcommand{\bcor}{\begin{cor}}
          
          \newcommand{\ecor}{\end{cor}}
          \newcommand{\bcon}{\begin{con}}
          \newcommand{\econ}{\end{con}}
          \newcommand{\bcond}{\begin{cond}}
          \newcommand{\Loc}{{\cal L}}
          \newcommand{\econd}{\end{cond}}
          \newcommand{\bteo}{\begin{teo}}
          \newcommand{\eteo}{\end{teo}}
          \newcommand{\brm}{\begin{rmk}}
          \newcommand{\erm}{\end{rmk}}
          \newcommand{\blem}{\begin{lem}}
          \newcommand{\elem}{\end{lem}}

          \newcommand{\ben}{\begin{enumerate}}
          \newcommand{\een}{\end{enumerate}}
          \newcommand{\bei}{\begin{itemize}}
          \newcommand{\eei}{\end{itemize}}
          \newcommand{\bdf}{\begin{defin}}
          \newcommand{\edf}{\end{defin}}

          \renewcommand{\>}{&>&}

          \newcommand{\fr}{\frac}
          \renewcommand{\r}{{\mathbb R}}
          \newcommand{\br}{\bar{\mathbb R}}
          \newcommand{\Z}{{\mathbb Z}}
          \newcommand{\cG}{{\cal G}}

          \newcommand{\R}{{\mathbb R}}
          
          \newcommand{\E}{{\mathbb E}}
          \newcommand{\bs}{{\tilde \b}}
          
          \renewcommand{\P}{{\mathbb P}}

          \newcommand{\N}{{\mathbb N}}

	\newcommand{\M}{{\cal M}}
          \newcommand{\W}{{\cal W}}

          \newcommand{\cE}{{\cal E}}
          
          \newcommand{\h}{{\cal H}}
          \newcommand{\f}{{\cal F}}
          \newcommand{\cV}{\mathbb{Z}_{even}^2}

          \newcommand{\cU}{{\cal U}}
          \newcommand{\cUb}{{\cal U}^{b_\b} }
          \newcommand{\cUbk}{{\cal U}^{b_\b,k_\b} }

          \renewcommand{\b}{\beta}

          \newcommand{\half}{\frac{1}{2}}

          \newcommand{\e}{\epsilon}

          \newcommand{\s}{\sigma}
          \renewcommand{\o}{\Pi}

          \newcommand{\binfty}{\lim_{\beta\uparrow\infty}}

\newcommand\ii{\item}

\newcommand{\btt}{\begin{theorem}}
\newcommand{\ett}{\end{theorem}}
\newcommand{\Net}{\mathcal{N}}

\newcommand{\daw}{\downarrow}
\newcommand{\uaw}{\uparrow}
\newcommand{\raw}{\rightarrow}

\newcommand{\Wl}{\mathcal{W}_l}
\newcommand{\Wr}{\mathcal{W}_r}
\newcommand{\be}{\begin{equation}}
\newcommand{\ee}{\end{equation}}

\newcommand{\Monetwo}{{\cal M}_{(1,2)}^b}
\newcommand{\Mzerotwo}{{\cal M}_{(1,1)}^k}
\newcommand{\Mzero}{{\cal M}}

\newcommand\no{\nonumber}
\newcommand\eps{\epsilon}

\newcommand\bss{{\tilde {\tilde \b}}}

          \newcommand{\cL}{{\cal T}}
          \newcommand\sqr{\vcenter{
          \hrule height.1mm
          \hbox{\vrule width.1mm height2.2mm\kern2.18mm\vrule width.1mm}
          \hrule height.1mm}}        

\title{Brownian Net with Killing}
\author[1]{C.M. Newman}
\author[2]{ K. Ravishankar}
\author[4,3]{E. Schertzer}
\affil[1]{Courant Institute of Mathematical Sciences, New York University, New York, NY 10012, USA.}
\affil[2]{Department of Mathematics, SUNY College at New Paltz, New Paltz, NY 12561, USA.}     
\affil[3]{UPMC Univ. Paris 06,
Laboratoire de Probabilit\'es et Mod\`eles Al\'eatoires, CNRS UMR 7599, Paris, France.}
\affil[4]{Coll\`ege de France,
Center for Interdisciplinary Research in Biology, CNRS UMR 7241, Paris, France.}     

\begin{document}

\maketitle

{\bf Abstract.} Motivated by its relevance for the study of perturbations of
one-dimensional voter models, including stochastic Potts models at
low temperature, we consider diffusively rescaled coalescing
random walks with branching and killing. Our main result is
convergence to a new continuum process, in which the random
space-time paths of the Sun-Swart Brownian net are terminated
at a Poisson cloud of killing points. We also prove existence
of a percolation transition as the killing rate varies. Key issues
for convergence are the relations of the discrete model killing points
and their intensity measure to the continuum counterparts.

\newpage

\tableofcontents

\newpage

\section{Introduction}
\label{intro}

\begin{figure}
\begin{minipage}[b]{0.5\linewidth}
\centering
\includegraphics[scale=.3]{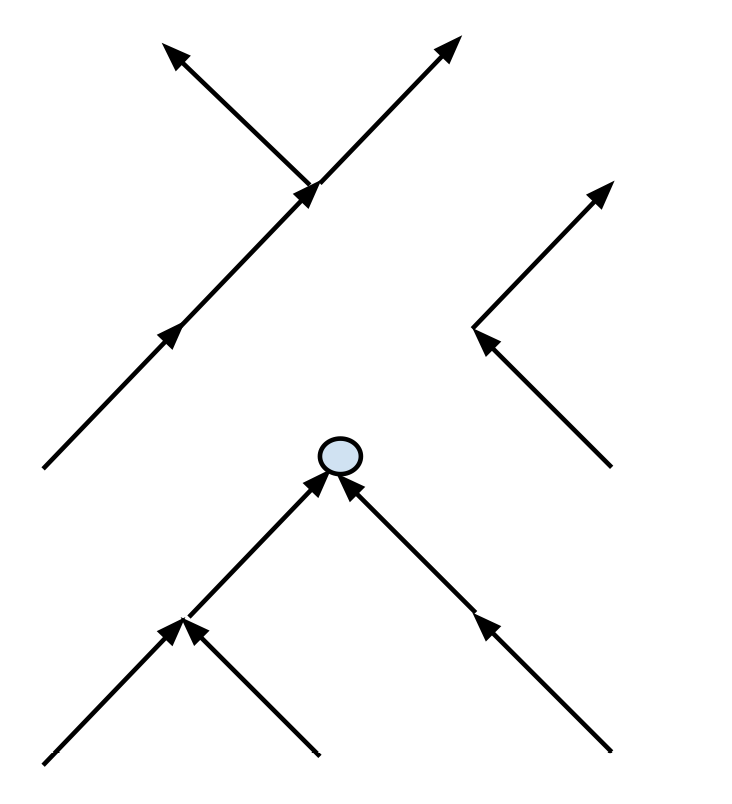}
\caption{BCK (Branching, Coalescing, Killing) configuration.}
\label{sim}
\end{minipage}
\hspace{0.5cm}
\begin{minipage}[b]{0.5\linewidth}
\centering
\includegraphics[scale=.3]{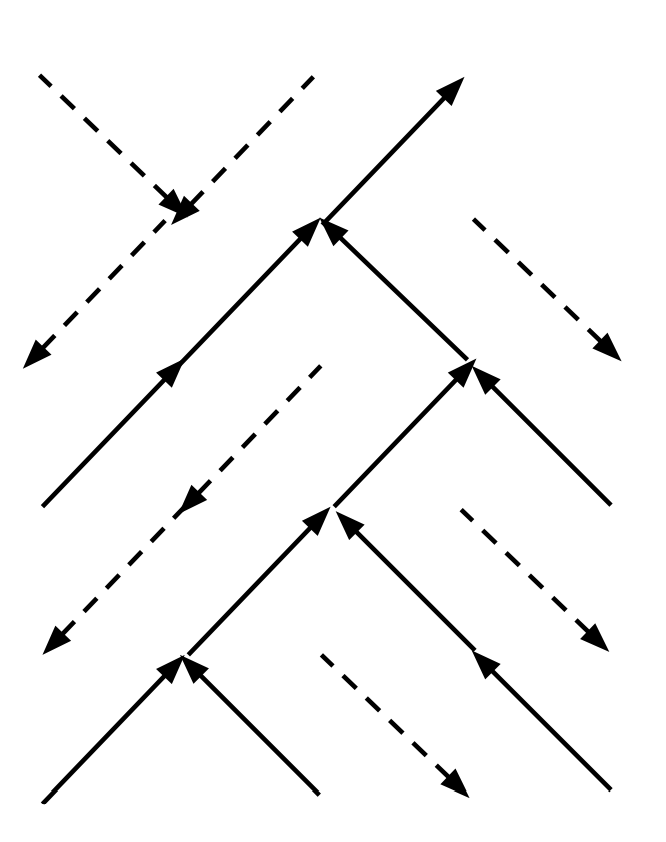}
\caption{Construction of the dual of the discrete web.}
\label{Duality}
\end{minipage}
\end{figure}

{\bf The model.} In the present paper,
we consider the natural scaling limit of a generalization
of the one-dimensional oriented percolation model, as introduced in \cite{MNR13}.
The model 
is parametrized by two non-negative numbers $b,k\geq0$
and can defined as followed.
Let us consider
\be\label{Zeven}
\cV=\{(x,t)\in \Z^2 \ : x+t \ \mbox{is even}\},
\ee
where $x$ is interpreted as a space coordinate
and $t$ as a time coordinate.
Each site $v=(x,t)\in\cV$
has two nearest neighbors
with higher time coordinates: $v_r=(x+1,t+1)$ and $v_l=(x-1,t+1)$.
$v$
is then randomly (and independently for different $v$'s)
connected 
to a subset of its neighbors
$v_r$ and $v_l$ by drawing arrows
according
to the following distribution.
\bi 
\item with probability $b$, draw the
two arrows $(v\raw v_r)$ and $(v\raw v_l)$;
\item do not draw any arrow with 
probability $k$, in which case it is called a killing point;
\item with the remaining probability, draw a
single arrow, its direction being chosen uniformly
at random  (see Figure \ref{sim}).
\ei

If we denote by $\cE^{b,k}$ the resulting random
arrow configuration,
the random directed graph $\cG^{b,k}=(\cV,\cE^{b,k})$
defines a certain type of
one dimensional percolation model
oriented forward in the $t$-direction.
By definition, a path $\pi$ along $\cE^{b,k}$
will denote
a path starting from any site of $\Z_{even}^2$
and following the random arrow configuration
until getting killed (by reaching a killing point) or reaching $\infty$. More precisely, 
a path $\pi$ is the graph of a function
defined on an interval in $\R$
of the form
$[\sigma_\pi,e_\pi)\label{sigfor}$, with $\sigma_\pi,e_\pi\in\Z\cup\{\infty\}$ such that
$e_{\pi}=\infty$ or else
$(\pi(e_\pi),e_\pi)$ is a killing
point; for every integer $t\in[\sigma_\pi,e_\pi)$,
$(\pi(t),t)$ connects to $(\pi(t+1),t+1)$
and $\pi$ is linear between $t$ and $t+1$.

Considering the set of all
the paths along  $\cE^{b,k}$,
one generates
an infinite family $\cU^{b,k}$
that can loosely be described
as a collection of graphs of 
one dimensional coalescing
simple random walks that branch
with probability $b$ and are killed with probability $k$.
A walk at space-time site
$v$ can create two new walks
(starting respectively at $v_l$
and $v_r$)
with probability $b$ and can be killed
with probability $k$; two walks 
move independently when they are apart
but become perfectly correlated (i.e., they coalesce)
upon meeting at a space-time point. 
In the following, $\cU^{b,k}$
will be referred to
as a system of branching-coalescing-killing random
walks (or in short, BCK) with
parameters $(b,k)$.

This model encompasses several classical models from statistical mechanics. 
For $p\in(0,1)$
and $b=p^2$, $k=(1-p)^2$, one recovers the standard 
one dimensional oriented percolation model. When $b=k=0$,
the trajectories in our random graphs are distributed as coalescing random walks.
In the present work, we will investigate the behavior of this model 
when the parameters $b$ and $k$ are non-zero, but small. This generalizes
previous work by Sun and Swart \cite{SS07}, and by the present authors \cite{NRS08},
where the case $k=0$ and small $b$ was investigated.

\bigskip

{\bf Motivations.} Our interest stems from several applications in 
interacting particle systems. 
On the one hand, it is well known 
that the classical voter model  
is dual to coalescing random walks. On
the other hand, Cox, Durrett and Perkins \cite{CDP11}
noticed that several models in ecology --- the spacial Lokta-Volterra model \cite{NP99}, the
evolution of cooperation \cite{OHLN06} --- and in statistical mechanics --- non linear voter models \cite{MoDDGL99}
but also the stochastic Potts model in $d=1$, see \cite{MNR13} ---
can be seen as perturbations (of strength $\eps$) of the voter model in a certain range of 
their parameter space. 
More precisely, they noticed that in a certain range of the transition rates, those models
can be written in the form 
\beqn
\label{rate-general}
\forall j\in\{1,\cdots,q\}, \ \ c_{j}^\eps(\eta,x) \ = \  c_{j}^v(\eta,x)+ c_{j}^{*,\eps}(\eta,x) 
\eeqn
where $c^\eps_{j}(\eta,x)$ is the transition rate of site $x$ to state $j$ given a configuration $\eta$,
$c_{j}^v$
is the transition rate of a standard (possibly non nearest neighbor) voter model and
the remaining term is such that $c_{j}^{*,\eps}(\eta,x) \rightarrow 0$ as $\eps\rightarrow0$.
 
 In general, the interacting particle systems alluded to above 
 are difficult 
 to study either because of a lack of monotonicity (Lotka-Volterra model)
 or because of the intrinsic complexity of
 the model (certain models of evolution of cooperation). Since 
 the voter model is well understood, the idea is to use
 standard results about that model
 to derive some properties  (such as coexistence of species in the Lokta Voterra model)
 of the more complex models described above. This program 
 was carried out successfully in \cite{CDP11}, where
 it is shown that for $d\geq3$,
 the properly rescaled local density of particles converges 
 to the solution of a reaction diffusion equation, and that properties
 of the underlying particle system can be derived from the behavior
 of this PDE. 
 
 The present work is a first step to understand the situation in low dimension, when $d=1$.
 The idea is the following.
Since the voter model
 is dual to coalescing random walks,
 perturbation of the voter model
 should be dual to perturbation 
 of coalescing
 random walks. As a consequence,
 it is natural to consider 
 a system 
 of coalescing
 random walks with an additional perturbative
 branching and killing mechanism. In fact,
 in a future work \cite{NRS13}, it will be shown that very general
 perturbations of voter models in
 dimension 1 converge to a continuum
 object which is constructed via the scaling  
 limit of the BCK: an object that we name the Brownian net with killing and which is
 the central object of this work.
 
 \bigskip
 
 {\bf Scaling of the parameters.} As we shall see below (see Theorem \ref{onewebteo}), the Brownian net
 with killing emerges as the scaling limit of the BCK model when 
 the branching and killing parameters are scaled differently. More precisely,
 if we denote $k_\b$ and $b_\b$
 the killing and branching parameters of the model,
 then 
 if $b_\b= e^{-\b} b$
 and $k_\b =e^{-2\b} k$ for some $b,k\geq0$,
 the (properly rescaled) BCK converges to the Brownian net with killing
 as
 $\b\rightarrow\infty$. 
 
Although scaling branching and killing
 differently might appear quite artificial at first sight,
 we will show in  \cite{NRS13} that such scaling 
 appears naturally in statistical mechanics. In particular, when considering the stochastic Potts model in one dimension,
 we will show that at large inverse temperature $\b$,
 such a model is asymptotically dual
 to the BCK
 with branching and killing parameters
 $$
 b_\b = \frac{q}{2} e^{-\b} \ \ \ k_\b = q e^{-2\b}, 
 $$
 where $q$ is the number of colors in the system. For this reason, we will
 parametrize the BCK with $\b$  and interpret this number as the inverse temperature of
 the model.

\bigskip

{\bf Main results.}
To illustrate our approach, let us first consider the case $b=k=0$. For every $(x,t)\in\Z_{even}^2$,
let us denote by $\pi_{x,t}$
the unique path
starting from the spacial point $x$ at time $t$
and following the random arrow configuration up to $\infty$.
As mentioned 
earlier,
this path
is simply the graph of a simple symmetric random walk. Furthermore,
the paths starting from different locations $(x_i,t_i)$
are coalescing, i.e. they move independently and then coalesce
upon meeting at any space time point. 

Let $n\in\N$ and let $\{(x_i,t_i)\}_{i\leq n}$ be
$n$ points in $\Z_{even}^2$. If one rescales
space and time diffusively by the transformation
\be\label{scal-transf}
S_{\b}(x,t) \ = (x e^{-\b},t e^{-2\b}),
\ee
it is then easy to see that if $\{(x_i^\b,t_i^\b)\in \Z_{even}^2 \}_{\b}$
are such that 
$$
\forall i\leq n, \lim_{\b\uaw\infty} \ \ S_\beta(x_{\b}^i,t_{\b}^i)\rightarrow (x^i,t^i),
$$
then the family of rescaled paths $\{S_\beta\left(\pi_{(x_i^\beta,t_i^\beta)}\right)\}$ converges
to a family of graphs of $n$ coalescing one dimensional Brownian motions
starting from the space-time points $\{(x_i,t_i)\}$ 
--- see T\`oth and Werner \cite{TW98} and Arratia \cite{A81}. 
In \cite{FINR04}, Fontes, Newman, Isopi and Ravishankar strengthened this statement
by proving that the set of all paths in the random configuration (also called the discrete web) converges  
to a continuum object that they called the Brownian web (BW). 

\begin{figure}
\centering
\includegraphics[scale=.3]{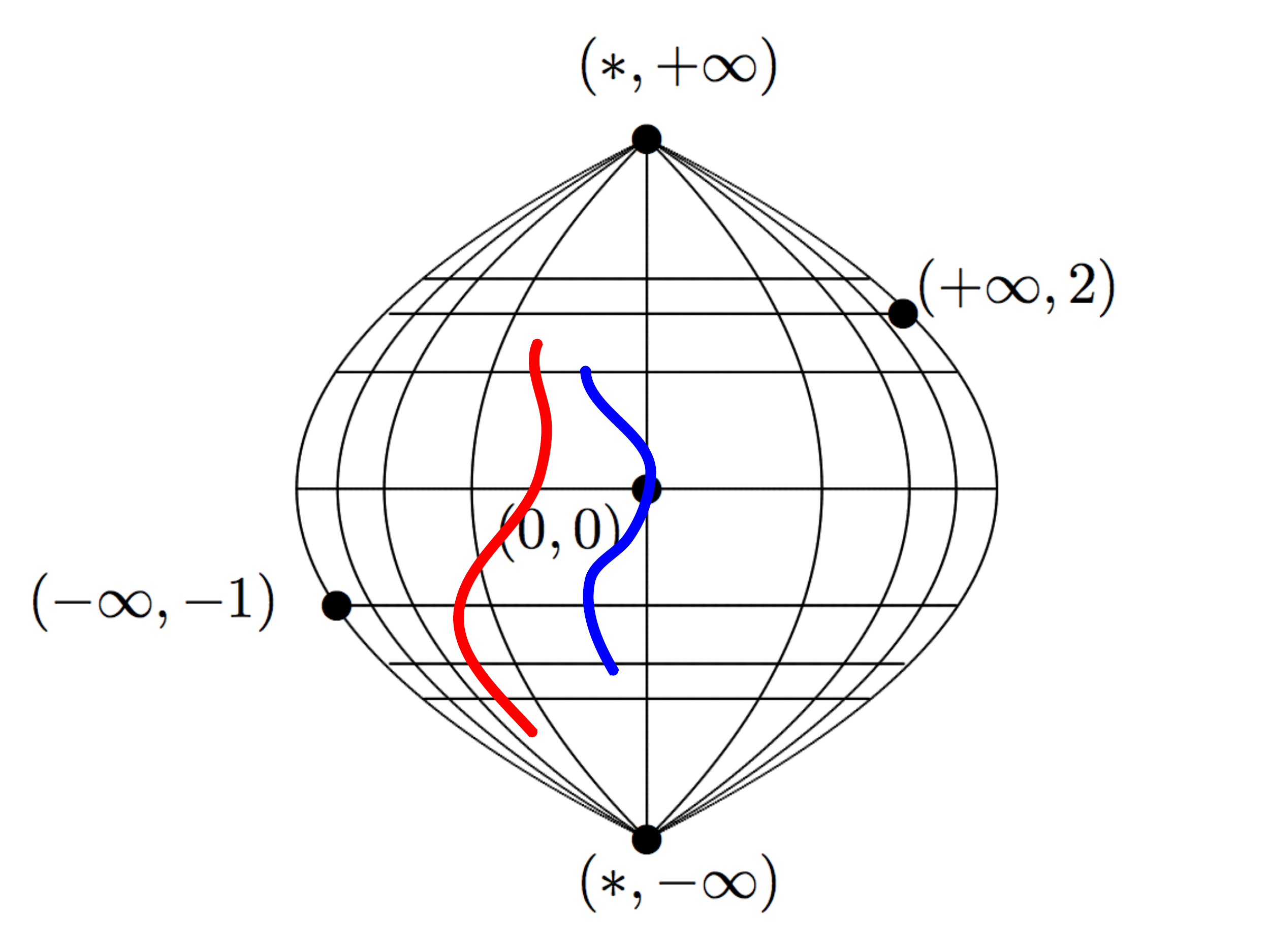}
\caption{Informally, the distance between two paths is obtained by first mapping the two paths into a compact subset of
the plane (blue and red curves) and then taking the distance between the  images.}
\label{distance}
\end{figure}

In Section \ref{the-space} below, we will
precisely define the probability space on which the Brownian web
lives, 
as it was introduced in \cite{FINR04}.
For the sake of presentation, we
start with an informal discussion (see also Fig. \ref{distance}).
The  Brownian web is a random
collection of continuous paths with specified starting points and ending points in space-time.
The paths take values in a metric space $(\br^2,\rho)$
which is a compactification of $\r^2$.  $(\o,d)$ denotes the space
whose elements are paths with specific starting points and ending points. 
The metric
$d$ is then defined as the
maximum of the sup norm of the distance between
two paths, the distance between their respective starting points and 
the distance between their ending
points. (In particular, when no killing occurs,
as in the discrete and Brownian webs, the ending
point of each path is $\infty$). 
The Brownian web takes
values in the metric space $(\h,d_\h)$,
whose elements are compact collection of paths in
$(\o,d)$ with $d_\h$ being the induced Hausdorff metric. Thus the
Brownian web is an
 $(\h,\f_\h)$-valued random variable,
 where $\f_\h$ is the Borel $\s$-field  associated to the metric
 $d_\h$. 

In ~\cite{FINR04},
the authors define the 
the Brownian web as the scaling limit 
of the discrete web.
More precisely, it is shown that there exists  
a random variable $\cal W$ on $({\cal H},d_{\cal H})$ --- the Brownian web ---
such that 
$$
S_\b( W ) \rightarrow {\cal W} \ \ \mbox{in law,}
$$
where $W$ denotes the discrete web. Building on their approach,
we will show an analogous result for the BCK when
$b,k$
are nonzero:

\bteo
\label{onewebteo}
Let $b,k\geq0$ and
$\lim_{\b\uaw\infty} k_\beta e^{2\b} = k$,
$\lim_{\b\uaw\infty} b_\beta e^{\b} = b$.
There exists a random variable $\Netbk$ taking values in $(\h,\f_\h)$
so that 
$$
S_{\b}( {\cal U}^{b_\b,k_\b} )\ \raw \  \Netbk
\ \ \ \ \text{in law.}
$$
\eteo

In the following, 
we will say that 
${\cal N}^{b,k}$
percolates (more accurately we should say ``percolates from the origin'')
if there exists a 
path in ${\cal N}^{b,k}$ starting from the origin
and reaching to $\infty$. The following
result is a continuum analog of 
a result proved for the discrete BCK in \cite{MNR13}.
\bteo
\label{perco}
\begin{enumerate}
\item $\Net^{b,k}$ is identical in law with $S_{\ln(b)}({\cal N}^{1,k/b^2})$.
\item There exists $k_c\in(0,\infty)$ such that 
\begin{itemize}
\item for every $k<k_c$, $\P(\Net^{1,k} \ \mbox{ percolates})>0$.
\item for every $k>k_c$, $\P(\Net^{1,k} \ \mbox{ percolates})=0$.
\end{itemize}
\end{enumerate}

\eteo

\brm Note that the convergence result, Theorem \ref{onewebteo}
immediately
implies the
scaling property
$$
S_{\ln(b)}(\Net^{b,k})  \ =  \  \Net^{1,\frac{k}{b^2}} \ \ \mbox{in law}.
$$ 
Indeed, 
 for
$e^{\b} b_\b \rightarrow b$ and $ e^{2\b} k_\b \rightarrow k$,
we have
$$
S_{\b + \ln(b)}( {\cal U}^{b_\b,k_\b} ) \rightarrow \Net^{1,\frac{k}{b^2}},
$$
whereas 
$$
S_{\b  + \ln(b)}( {\cal U}^{b_\b,k_\b} )  = 
S_{\ln(b)}\circ S_ {\b}( {\cal U}^{b_\b,k_\b} ) \rightarrow S_{\ln(b)}(\Net^{b,k}) .
$$
\erm

\bigskip

{\bf Outline.} The rest of the paper is organized as follows.
In  Section \ref{BNK}, we give two alternative constructions of
the Brownian net with killing.  
In Section \ref{Properties}, we prove the percolation property
of the killed Brownian net stated in Theorem \ref{perco}. Finally,
in Section \ref{Invariance::Principle},
we show that the BCK
converges (under proper rescaling)
to the Brownian net with killing, as stated in Theorem \ref{onewebteo}.
A key issue in proving
that convergence, as discussed at the beginning of Section 4 and then treated in Sections 4.1 and 4.2, is the convergence of the Poissonian
killing points in such a way that if a continuum path is killed somewhere,
then the discrete approximating paths are killed as well.

\section{Construction of the Brownian Net with Killing}
\label{BNK}

In this section, we give
two alternative constructions of the object $\Net^{b,k}$.
We first briefly outline 
the ideas behind those two constructions.

\bigskip

{\bf One web, two markings.} In a first approach, 
in the spirit of \cite{NRS08}, 
we start by describing the 
scaling limit in the special case
$b=k=0$. In this particular setting,
where the branching and killing parameters are turned off,
the discrete BCK becomes an infinite family
of coalescing random walks, also referred to 
as a discrete web --- see \cite{FINR04}.

As we shall 
see in Section \ref{ooo},
the Brownian web can be
used to construct
the natural scaling limit for the BCK, 
the Brownian net with killing $\Net^{b,k}$
in the case 
where $b,k\neq0$. 
At the discrete level, 
the idea is simple
and relies on the idea that the BCK can be constructed 
by starting
with a discrete web and then
turn on the killing and branching parameters.
Effectively, one starts with a discrete web,
and then independently at each site
\bi
\item removes an arrow with probability $k$.
\item or adds an extra arrow with probability $b$.
\ei
In Section \ref{special_points}, 
we will explain how one can
generalize this construction to the continuum level.
More precisly, we will identify 
some geometrical random configurations in the Brownian
Web, playing
a role analogous to the arrows of the discrete web (the BCK with no branching or killing).
Finally,
we will briefly explain how those ``continuum arrows" can 
be added or removed by two Poissonian markings in order to generate the Brownian net with killing $\Net^{b,k}$ from the Brownian web ${\cal W}$.
\bigskip

{\bf One net, one marking.} Alternatively, our second construction (presented in Section \ref{netconstruction}) will only make use 
of a single set of marks, while using a richer underlying structure: the standard Brownian
net (with no killing), an object introduced using different methods
by Sun and Swart \cite{SS07}  and analysed by Newman, Ravishankar and Schertzer \cite{NRS08}.   
In Section \ref{standard-net}, we first 
present the three distinct direct constructions of the standard Brownian 
net, due to Sun and Swart. 
In Section \ref{onenet}, we will show how one can 
directly construct the killed 
Brownian net from 
the standard Brownian net by adding one extra set of marks.
This direct construction
will be based on some intermediate results
provided in Sections \ref{special-force} and \ref{sec:bc}.
We note that those will also be relevant for other parts of this paper. 
In Section
\ref{special-force}, we give some results about the special points of the Brownian net.
In Section \ref{sec:bc}, we introduce 
the branching-coalescing point set
of the standard Brownian net. 

\bigskip

Finally, in Section
\ref{equivalence:construction},
the equivalence between our two alternative constructions of
the killed Brownian net is established.

\subsection{The Space $({\cal H}, d_{\cal H})$}
\label{the-space}

As in \cite{FINR04}, we will define the Brownian net with killing 
as a random compact set of paths. In this section, we
briefly 
outline the construction of the space of compact sets.
For more details, the interested reader may refer 
to \cite{FINR04}.

First define $(\bar\R^2,\rho)$ to be 
the compactification of $\R^2$
with
$$
\rho((x_1,t_1),(x_2,t_2)) \ = \ | \phi(x_1,t_1) - \phi(x_2,t_2)| \vee |\psi(t_1)-\psi(t_2)|,
$$
where
$$
\phi(x,t) = \frac{\tanh(x)}{1+|t|}  \ \ \ \mbox{and}  \ \ \  \psi(t) = \tanh(t).
$$
In particular, we note that the mapping
$(x,t)\rightarrow(\phi(x,t),\psi(t))$
maps $\bar\R^2$ onto a compact subset of $\R^2$.

Next, let $C[t_0,t_1]$
denote the set of continuous functions from 
$[t_0,t_1]$ to $[-\infty,+\infty]$. From there, we define the set 
of continuous paths in $\R^2$ (with a prescribed
starting and ending point) as
$$
\Pi := \cup_{t_0\leq t_1} C[t_0,t_1]\times \{t_0\}\times\{t_1\}.
$$
Finally, 
we equip this set of paths with a metric
$d$, defined as the
maximum of the sup norm of the distance between
two paths, the distance between their respective starting points and 
the distance between their ending
points. (In particular, when no killing occurs,
as in the forward Brownian web, the ending
point of each path is $\infty$).
More precisely,
if for any path $\pi$,
we denote by $\sigma_\pi$
the starting time of $\pi$
and by $e_\pi$
its ending time,
we have
\beqnn\label{diistance}
d(\pi_1,\pi_2) \ = \
 | \psi(\sigma_{\pi_1}) - \psi(\sigma_{\pi_2})  | 
\vee
| \psi(e_{\pi_1}) - \psi(e_{\pi_2})  | 
\vee 
\max_{t} \  \rho((\bar\pi_1(t),t),(\bar\pi_2(t),t)),
\eeqnn
where $\bar \pi$ is the extension of $\pi$
into a path from $-\infty$ to $+\infty$
by setting
$$
\bar \pi(t) = \left\{ \begin{array}{c} 
\pi(\sigma_\pi) \ \ \mbox{for $t<\sigma_\pi$},  \\
\pi(e_\pi) \ \ \mbox{for $t>e_\pi$.}  \end{array} \right.
$$
Finally, let ${\cal H}$ denote the set of compact subsets of $({\cal H},d_{\cal H})$
where $d_{\cal H}$ is 
the
Hausdorff
metric
\beqnn\label{dh}
d_{\cal H}(K_1,K_2) = \max_{\pi_1\in K_1} \min_{\pi_2\in K_2} d_{\cal H}(\pi_1,\pi_2) \vee  
\max_{\pi_2\in K_2} \min_{\pi_1 \in K_1} d_{\cal H}(\pi_1,\pi_2).   
\eeqnn
In \cite{FINR04}, it is proved that $({\cal H},d_{\cal H})$ is Polish. In the following, 
we will construct the Brownian net with killing as a random element of this space.

\subsection{One Web, Two Markings}
\label{ooo}

\subsubsection {The Brownian Web}
\label{Forward_Web}
As mentioned in the introduction, 
the Brownian web is the scaling limit of the discrete web
under diffusive space-time scaling
and is defined as a random element
of $({\cal H},d_{\cal H})$.
The next theorem, taken from~\cite{FINR04}, gives some of the key
properties of the BW.

          \bteo
          \label{teo:char}
          There is an \( ({\cal H},{\cal F}_{{\cal H}}) \)-valued random variable
          \(
          {\W} \)
          whose distribution is uniquely determined by the following three
          properties.
          \begin{itemize}
                    \item[(o)]  from any deterministic point \( (x,t) \) in
          $\r^{2}$,
                    there is almost surely a unique path \( {B}_{(x,t)} \)
          starting
                    from \( (x,t) \).

                    \item[(i)]  for any deterministic, dense countable subset
        \(
          {\cal
                    D} \) of \( \r^{2} \), almost surely, \( {\W} \) is the
          closure in
                    \( (\Pi, d) \) of \( \{ {B}_{(x,t)}: (x,t)\in
                    {\cal D} \}. \)
  \item[(ii)]  for any deterministic  $n$ and \((x_{1}, t_{1}),
        \ldots,
                    (x_{n}, t_{n}) \), the joint distribution of \(
                    {B}_{(x_{1},t_{1})}, \ldots, {B}_{(x_{n},t_{n})} \) is that
                    of coalescing Brownian motions from those starting points (with unit diffusion
        constant).

          \end{itemize}
          \eteo

Note that (i) provides a practical construction of the Brownian web. For $\cal D$ as defined above,  construct coalescing Brownian motion paths starting from $\cal D$. This defines a {\it skeleton} for the Brownian web that is denoted by $\W({\cal D})$. $\W$ is simply defined as the closure of this precompact set of paths.

\subsubsection{The Backward (Dual) Brownian Web}
\label{backwardBW}

The discrete web is defined on $\Z_{even}^2$ 
and oriented forward in time. 
There is a natural dual system of paths
defined on the complementary lattice 
$\Z_{odd}^2:=\{(x,t) \ : \ x+t \ \ \mbox{is odd}\}$ and oriented in the reverse direction.
Indeed, given a realization 
of $W$, at any vertex $(x,t)\in\Z_{odd}^2$,
we can construct a (backward) arrow 
configuration by
rotating
the (forward) arrow
at $(x,t-1)\in\cV$ through $180^\circ$ and then
translating its starting point to $(x,t)$
(see Fig. \ref{Duality}).
This
defines a system of
coalescing random walks starting from $\Z_{odd}^2$, running backward in
time without crossing the forward discrete web paths.
Furthermore, 
the resulting system of paths
is easily seen to be identically
distributed with the original model,
after a $180^\circ$ rotation followed
by a unit
translation in the $x$ or $t$ direction. 

This family will be referred to as the backward discrete web, 
and
the backward (dual)
BW $\hat{\mathcal{W}}$ may be defined analogously as a
functional of the (forward) BW $\mathcal{W}$.
More precisely for a countable dense deterministic set of space-time points,
the backward BW path from each of these is the (almost surely) unique
continuous curve (going backwards in time) from that point that
does not cross (but may touch) any of the (forward) BW paths;
$\hat{\mathcal{W}}$
is then the closure of that collection of paths. The
first part of the next proposition
states that the ``double BW'', i.e., the pair
$(\mathcal{W},\hat{\mathcal{W}})$, is the diffusive
scaling limit of the corresponding discrete pair
$(W,\hat{W})$ (after diffusive scaling and 
as the scale parameter $\b \to \infty$).
Convergence in the sense of weak convergence of probability
measures on $(\h,\f_\h) \times (\hat \h,\hat \f_\h)$ was proved in~\cite{FINR04};
convergence of finite dimensional distributions and
the second part of the proposition were already contained
in~\cite{TW98}.

\bprop
\label{forwardandbackward1}
\begin{enumerate}
\item Up to a reflection with respect to the $x$-axis,
$\W$ and $\hat\W$ are identically distributed.
\item Invariance principle : $\left(S_\eps(W),S_\eps(\hat{W})\right) \to
(\mathcal{W},\hat{\mathcal{W}})$ as $\eps\to0$.
\item For any (deterministic)
pair of points $(x,t)$ and $(\hat x, \hat t)$ there is
almost surely a unique forward path $B$ starting from $(x,t)$ and a
unique backward path $\hat B$ starting from $(\hat x,\hat t)$.
\end{enumerate}

\eprop

The next proposition, cited from~\cite{STW00}, which gives the
joint distribution of a single forward and single backward BW path,
has an extension to the joint distribution of finitely many
forward and backward paths. We remark that that extension
can be used to give a characterization (or construction)
of the double Brownian web $(\W,\hat \W)$ analogous to the one for
the (forward) BW from
Theorem~\ref{teo:char} --- see~\cite{STW00, FINR05} for more details.

\bprop
\label{forwardandbackward2}
\begin{enumerate}
\item \underline{Distribution of $(B,\hat B)$}: Let $(B_{ind},\hat B_{ind})$
be a pair of independent forward and backward Brownian motions starting
at $(x,t)$ and $(\hat x,\hat t)$ and let $(R_{\hat B_{ind}}(B_{ind}),\hat B_{ind})$ be
the pair obtained after reflecting (in the Skorohod sense) $B_{ind}$ on
$\hat B_{ind}$, i.e., $R_{\hat B_{ind}}(B_{ind})$ 
is the following function of $u\in[t,\hat t]$:
 \begin{equation}
R_{\hat B_{ind}}(B_{ind})=\left\{ \begin{array}{ll}
B_{ind}(u)-0 \wedge \min_{t \leq v \leq u } ( B_{ind}(v)-\hat B_{ind}(v))\
\ \ \ \textrm{on $\{B_{ind}(t) \geq \hat B_{ind}(t) \}$}, \\
B_{ind}(u)-0 \vee \max_{t \leq v \leq u } ( B_{ind}(v)-\hat B_{ind}(v))\
\ \ \ \textrm{on $\{B_{ind}(t) < \hat B_{ind}(t) \}$}. \end{array}
\right.
\end{equation}
Then
\begin{equation}
(R_{\hat B_{ind}}(B_{ind}),\hat B_{ind})=(B,\hat B) \ \ \textrm{in law},
\end{equation}
where $B$ is the path in $\W$ starting at $(x,t)$ and $\hat B$ is 
the path in $\hat W$ starting at $(\hat x,\hat t)$.
\item Similarly,
\begin{equation}
(B_{ind},R_{B_{ind}}(\hat B_{ind}))=(B,\hat B) \ \ \textrm{in law}.
\end{equation}
\end{enumerate}

\eprop

\subsubsection{Special Points of the Brownian Web}
\label{special_points}

While there is only a single path starting from (and no path passing through) any deterministic point in
$\R^2$ in both the forward and backward webs, there exist random
points $z \in \R^2$ with one or more than one path
passing through or starting from~$z$. As we shall see later,
those ``special points" of the Brownian web will
play a key role in our construction of the Brownian net with killing.
We remark that here when we say that a path starts from $(x_0,t_0)$,
it does not preclude the possibility that it is the $t\geq t_0$
continuation of a path that passes through $(x_0,t_0)$.

        We  start by describing the ``types'' of points $(x,t)\in\r^2$,
        whether deterministic or not. 
	We say that two paths $B,B'\in\W$ are equivalent paths entering 
	$z=(x,t)$, denoted by 
	\begin{equation}
	\label{str-equ}
	B=_{in}^z B'
	\end{equation} 
	iff $B=B'$ on $[t-\e,t]$
	for some $\e>0$. The relation $=^z_{in}$ is a.s. an equivalence relation on the set
        of paths in $\W$ entering the point $z$ 
	and we define $m_{in}(z)$ as the number of those equivalence 
	classes. ($m_{in}(z)=0$ if there are no paths entering $z$.)
        $m_{out}(z)$ is defined as the number of distinct paths starting from $z$.
	For $\hat \W$, $\hat m_{in}(z)$ and $\hat m_{out}(z)$ are defined similarly.

          \bdf
          The type of $z$ is the pair $(m_{in}(z),
          m_{out}(z))$.
          \edf

	\begin{figure}[ht!]
\begin{minipage}[b]{0.5\linewidth}
\centering
\includegraphics[scale=.3]{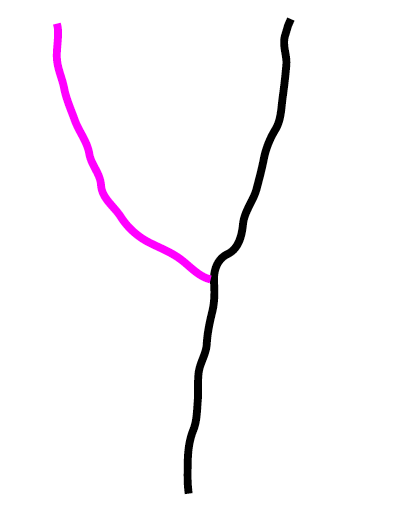}
\caption{A schematic diagram of a right $(m_{\mbox{{\scriptsize\emph{in}}}},m_{\mbox{{\scriptsize
          \emph{out}}}})=(1,2)$ point.
          In this example the incoming forward path connects to the rightmost
          outgoing path
          (with a corresponding dual connectivity for the backward paths), the left
          outgoing path is a newly born path.}
\label{fig12}
\end{minipage}
\hspace{0.5cm}
\begin{minipage}[b]{0.5\linewidth}
\centering
\includegraphics[scale=.3]{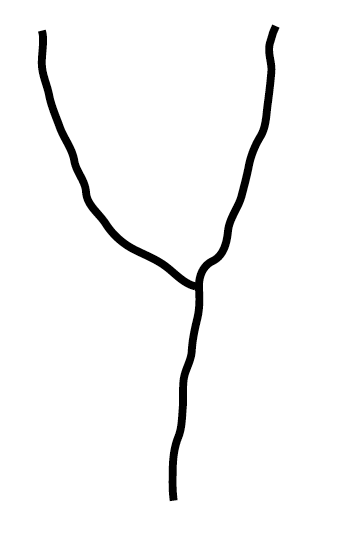}
\caption{A branching point is created by attaching
the left branch to the incoming path. 
}
\label{br}
\end{minipage}
\end{figure}

          The following results from \cite{TW98} (see also ~\cite{FINR05}) specify what types
          of points are possible in the Brownian web.

 \bteo
          \label{teo:types}
          For the Brownian web, almost surely, every $(x,t)$ is one of
          the
          following types, all of which occur: $(0,1)$, $(0,2)$, $(0,3)$,
          $(1,1)$,
          $(1,2)$, $(2,1)$.
\eteo

\bprop
\label{special-web}
	Let $t$ be a deterministic time. Almost surely for every
	realization of the Brownian web, any point of the form $(x,t)$
	(with $x$ deterministic or random) is either of type
	$(0,1)$, $(1,1)$ or $(0,2)$.
\eprop

\bprop
          \label{prop:dual}
          For the Brownian web, almost surely for {\em every}
          $z$ in $\r^2$, $\hat m_{in}
          (z)=m_{out}(z)-1$
          and
          $\hat m_{out}(z)=
          m_{in}(z)+1$.
          \eprop


Finally,
we close this section by proving
an easy lemma about the special points
of the Brownian web.
This result will 
be useful
in the rest
of the paper
but can be skipped at first
reading.

\blem
\label{special-web}
Let $w$ be a path in $\W$ starting 
from a point $z=(x,t)$.
If
$z$ is not a $(0,3)$ point, then
there exists a sequence
of paths $w_n\in\W$
starting strictly below
time $t$ (i.e. with $\sigma_{w_n}<t$)
and coalescing with
$w$ at $\tau_n>t$,
with $\tau_n \rightarrow t$. 
\elem

\begin{proof}
We first claim that 
if $z$ is not of type $(0,3)$,
for any outgoing path,
there exists a sequence $w_n$
starting below time $t$ such that
$w_n \rightarrow w$ in the topology $(\Pi,d)$.
By Theorem \ref{teo:types},
we need to show that
our result holds for points
of the type $(*,1)$ and $(*,2)$ ---
where $*$ can be 0,1 or 2 when 
there is a single outgoing path, and $0$
or $1$ for two outgoing paths.

Let us first show that
our claim holds when $z$ is of type $(*,1)$. Let us consider
any sequence
$w_n$ in $\W$
with starting point
$z_n = (x_n,\sigma_{w_n}) \rightarrow z$ and $\sigma_{w_n}<t$.
Since $\W$ is a compact set 
of paths,
from any subsequence
extracted from $w_n$,
one can extract a further
subsequence 
converging to a 
path starting from $z$.
Since $w$ is the unique 
such path, we conclude that
the sequence $w_n$
must converge to $w$. Let us now show
that our claim holds
for points of the type $(*,2)$. W.l.o.g.
we aim to show that
the left outgoing path
can be approximated 
by a sequence 
$w_n$ with $\sigma_{w_n}<t$.
Let $\hat w$ be the left-most backward 
path in $\hat \W$
starting from the point $z$
and let us consider
a sequence $w_n$
starting in the left region 
delimited
by $\hat w$.
Since dual paths in $\hat \W$
do not cross paths in the 
forward web $\W$, 
$w_n$
hits the time line $\R\times\{t\}$ either
strictly to the left of $z$ --- i.e.
at a point $(x_n,t)$
such that $x_n<x$ --- or at $z$.
In order to see that the latter does not occur,
we distinguish between two cases
\begin{enumerate}
\item[(i)] $z$ is of type $(0,2)$. Since there is no incoming Brownian web path at $z$, 
our claim is obvious.
\item[(ii)] $z$ is of type $(1,2)$. Recall that we chose $\hat w$ to be the left-most
backward path starting from $z$. For $(1,2)$ points, the incoming path is necessarily
to the right of this backward path (squeezed between the left-most and 
the right-most backward paths), implying that $w_n$
can not coalesce with $w$ before time $t$ (the time coordinate of $z$). Hence
if $w_n$ hits the point $z$, this would imply the existence of two non-equivalent
paths entering $z$, which would contradict the definition
of a $(1,*)$ point.
\end{enumerate}
Since points of type $(*,2)$ are necessarily of
one of the previous two types (again by Theorem \ref{teo:types}),
it follows that 
we must have $x_n<x$. 
This implies that
$w_n$ coalesces
with the left outgoing 
path before meeting the 
right outgoing path
starting from $z$. By
reasoning as
for the $(*,1)$ case,
this easily implies
that $w_n$ must converge to
the left outgoing path. 

It remains to show that
if
$w_n$
converges to $w$ with $\sigma_{w_n}<t$,
then $w_n$
must coalesce
with $w$
at a time $\tau_n$,
with $\tau_n\rightarrow t$. 
In order to see 
that, 
let us consider some arbitrary
small
$\eps>0$. In \cite{FINR04},
it is proved that
the set of paths 
with starting time 
$\sigma_w \leq t$,
hits the line
$\R\times\{t+\eps\}$ at
only locally
finitely many points.
Hence,
$w_n(t+\eps)$
must be stationary
after some rank,
implying
that $w_n(t+\eps)=w(t+\eps)$
and that $w_n$ must have
coalesced
with $w$
before time $t+\eps$.
Since $\eps$
is arbitrarily
small,
this ends the proof
of our lemma.

\end{proof}

\subsubsection{Construction of the Brownian Net with Killing}
\label{oneweb}

In our construction of the Brownian net with killing,
$(1,2)$ points and $(1,1)$
points will play a special role. 
First, it is important to realize that points of type $(1,2)$
can be characterized
in two ways, both of which will play a crucial role
in our construction of the Brownian net with killing. 1)  By Proposition
\ref{prop:dual}, 
$z\in \r^2$ is of
type $(1,2)$ precisely if both a forward and a backward path pass
through~$z$. 2) A single incident path continues along exactly
one
of the two outward paths --- with the choice determined intrinsically.
It is either left-handed or
right-handed according to
whether the continuing
path is to the left or to the right of the other
outgoing
path. For a left (resp., right) 
$(1,2)$ point $z$, the right (resp, left) outgoing path
will be referred to as the {\it newly born path} starting from $z$.
See Fig. \ref{fig12} for a schematic
diagram of the ``right-handed'' case.

In particular, one can add a ``continuum'' 
arrow at a given $(1,2)$
point $z$ by not only connecting an incoming
path to
the continuing outgoing path,
but also to the newly born path starting from $z$ --- see Fig. \ref{br}. 
In the discrete 
picture, this amounts to adding an arrow
which induces a ``macroscopic'' effect in the web.
This particular construction will be discussed later in this section.

On the other hand, 
a point of type $(1,1)$ can be 
seen as a single arrow passed through 
by a path, and
whose suppression 
will induce a macroscopic effect 
on the Brownian web, i.e., the suppression 
of this arrow will disconnect this path into two large components.

\bigskip

As already mentioned at the beginning of this section,
the discrete BCK can be constructed from the discrete web
by randomly adding and removing arrows.
Now that we have identified the random structures playing
the role of arrows at the continuum level,
we need to implement some random selection
of those structures. This will be done by first
defining two uniform measures on
those sets of "arrows": one measure for the 
$(1,2)$ points --- the local time measure, related to
the intersection local time between
the forward and backward paths of the 
Brownian web; see Proposition
\ref{ltm} below ---
and one for $(1,1)$ points --- the time length measure; see Proposition \ref{tlm}. 

\bprop[Time Length Measure]
\label{tlm}
Given a realization of the Brownian web, 
for every Borel set $E\subset\R^2$, define
$$
\cL_{web}(E) \ := \ \int_{-\infty}^\infty \ |\{x \ : \ (x,t)\in E\cap (1,1)\}| dt,
$$
where $(1,1)$ refers to the set of $(1,1)$ points of the Brownian web and $|\cdot|$
denotes cardinality. $\cL_{web}$
defines a $\sigma$-finite measure on $\R^2$
so that for every $\pi\in\W$
\be\label{e-cl}
\forall  \ x_1<x_2,t_1<t_2, \ \ \cL_{web}(\{(x,t)\in[x_1,x_2]\times[t_1,t_2] \ :  \pi(t) =x  \} ) 
\ = \ \int_{t_1 \vee \sigma_\pi}^{t_2 \vee \sigma_\pi} 1_{\pi(t)\in[x_1,x_2]}  \ \ dt.
\ee
\eprop

In other words,
if one considers the set defined as the intersection of
a rectangle with a Brownian web path,
the time length measure of this set
is given by the time spent by this path in the rectangle.

\begin{proof}
It is straightforward to check that $\cL_{web}$ defines a measure.
To prove the $\sigma$-finite property,
let 
us consider a dense countable deterministic set ${\cal D}:=\{ z_1,\cdots, z_n,\cdots\}$
in $\R^2$, and define ${\cal D}_n=\{ z_1,\cdots,z_n\}$. 
In \cite{FINR04} it was shown that  
for every path $\pi\in{\cal W}$ and every time $t$
such that $\sigma_\pi<t$,
we must have $(\pi(t),t)\in\W({\cal D})$, where $\W({\cal D})$
refers  temporarily to the trace of the skeleton ${\cal W}({\cal D})$.
It follows that $(1,1)\subset\W({\cal D})$ and that
\beqnn
\cL_{web}(E) & = &   \int_{-\infty}^\infty \ |\{x \ : \ (x,t)\in E\cap (1,1) \cap \W({\cal D})\}| dt \\
 	  & = &    \lim_{n\uaw\infty}\int_{-\infty}^\infty \ |\{x \ : \ (x,t)\in E\cap (1,1) \cap \W({\cal D}_n)\}| dt, 
\eeqnn 
where the last identity follows from the monotone convergence theorem. 
Finally, since $\W({\cal D}_n)$ meets any horizontal line $\R\times\{t\}$ at $n$ points at the most
(recall that there is a.s. a unique path starting from any deterministic point ),
this proves the $\sigma-$finite property for our measure.

It remains to show (\ref{e-cl}).
For every realization of the Brownian web $\W$,
let us consider the set of times
$$
T = \{t \ : \ \exists \pi\in\W  \ \ \mbox{with $\sigma_\pi<t$ and $(\pi(t),t)$ is not of type $(1,1)$} \ \}.
$$
In \cite{FINR04}, \cite{TW98}, it was shown that 
for every deterministic
time $t$, any point $(x,t)$, with $x\in\R$, is either of the type $(0,1)$, $(0,2)$ or $(1,1)$. Hence,
at a determistic time $t$, any point $(x,t)$ such that there exists some
$\pi\in\W$ with $\sigma_\pi<t$ so that 
$x=\pi(t)$ must be a $(1,1)$ point, which amounts to saying that
$t$ does not belong to the set $T$ almost surely. 
By a direct application
of the Fubini Theorem, the set of times which belongs
to $T$ must have zero Lebesgue measure, implying that
\beqnn \cL_{web}(\{(x,t)\in[x_1,x_2]\times[t_1,t_2] \ :  \pi(t) =x  \} ) 
& = & \int_{t_1 \vee \sigma_\pi}^{t_2 \vee \sigma_\pi} 1_{\pi(t)\in[x_1,x_2]}1_{(\pi(t),t)\in (1,1)}  \ \ dt  \\ 
& = & \int_{t_1 \vee \sigma_\pi}^{t_2 \vee \sigma_\pi} 1_{\pi(t)\in[x_1,x_2]}  \ \ dt,  \\ 
\eeqnn
where the first equality follows directly from the definition of $\cL_{web}$.

\end{proof}

The following Proposition is cited from \cite{SSS13}.

\bprop[Local Time Measure]
\label{ltm}
Almost surely for every realization of the Brownian web,
there exists a unique $\sigma$-finite measure $\Loc$, concentrated on the set of points of type $(1,2)$
in $\W$, such
that for each $\pi\in\W$ and $\hat\pi\in\hat\W$:
\beqnn
\Loc(\{(x,t)\in\R^2 \ : \ \sigma_\pi < t < \sigma_{\hat\pi}, \ \pi(t) = x = \hat\pi(t) \} )\\
  =  \ \lim_{\epsilon \downarrow 0} \frac{1}{\epsilon} |\{t: \sigma_\pi < t < \sigma_{\hat\pi}, \  \ |\pi(t)-\hat \pi(t)| \leq \epsilon\}|, \
\eeqnn
where $\sigma_{\hat \pi}$ is the starting time of $\hat \pi$
and the limit in the RHS exists and is finite. 
\eprop

Next, conditioned on a realization of the Brownian Web $\W$,
and with the two random measures $\cL_{web}$ and $\Loc$ 
on hand,
we select a random set of marks
by defining two independent Poisson Point Processes
with respective intensity measure $k\cdot \cL_{web}$ and $b\cdot \Loc$.
In the following those two sets of points
will be respectively denoted by
 $\Mzerotwo$ and $\Monetwo$. 
It is easy to show that the local time measure
and time length measure of any 
bounded open set $O$ is infinite, 
implying that the number of marks
in such a set is infinite (but countable) 
almost surely.

\bigskip

We are now ready to construct the Brownian net with killing 
in which arrows 
are added to the Brownian web at the $\Monetwo$ points 
and arrows are suppressed at the points of $\Mzerotwo$.
More precisely, let $\{z_m\}$ be an enumeration of $\Monetwo$. 
We define a partial Brownian net
$\Netbk_{n}$ by
introducing branching  at the points of the partial
marking $\{z_m\}_{m=1}^n$ and killing at $\Mzerotwo$.
As explained in Section \ref{special_points}, if the $(1,2)$ point in Fig.~\ref{fig12}
is marked, then the set $\Netbk_{n}$ will include not only paths that
connect to the right outgoing path (as in the original web)
but also ones that connect to the left outgoing path. Furthermore,
such a path $\pi$ will be killed
when it encounters a point in $\Mzerotwo$, i.e.,
the ending time of the path will be
$$
e_\pi := \inf\{t>\sigma_\pi \ : \ \exists (x,t)\in\Monetwo  \ \ \mbox{with} \ \ \pi(x)=t \}.
$$
Finally, we define the Brownian net with killing $\Netbk_{web}$
(where the web subscript indicates
that the object has been constructed from the Brownian web)
as the closure of
$\bigcup_{n=1}^\infty \Netbk_{n }\bigcup_{z\in\R^2}\{z\}$,
where by a slight abuse of notation, $\{z\}$  
refers to the trivial path
starting and ending at $z$.

\subsection{One Net, One Marking.}
\label{netconstruction}

\subsubsection{The Standard Brownian Net}
\label{standard-net}
When $k=0$,  Sun and Swart give three alternative constructions 
of the Brownian net which are called
respectively the {\em hopping}, {\em wedge, and mesh characterizations}.
Those constructions are all based on the 
construction of two coupled drifted 
Brownian webs $(\W_l,\W_r)$. Following \cite{SS07}, we call $(l_1, \ldots, l_m; r_1,
\ldots, r_n)$ a collection of left-right coalescing Brownian
motions if $(l_1, \ldots, l_m)$ is distributed as coalescing
Brownian motions each with drift $-b$, $(r_1, \ldots, r_n)$ is
distributed as coalescing Brownian motions each with drift $+b$, paths
in $(l_1, \ldots, l_m; r_1, \ldots, r_n)$ evolve independently when
they are apart, and the interaction between $l_i$ and $r_j$ when they
meet is a form of sticky reflection. More precisely, for any $L \in
\{l_1, \ldots, l_m\}$ and $R\in\{r_1, \ldots, r_n\}$, the joint law of
$(L,R)$ at times $t>\sigma_L\vee \sigma_R$ is characterized as the unique
weak solution of

\begin{eqnarray}
d L(t)\ & = &  \ d B_l - b dt, \nonumber \\[0.1in]
d R(t)\ & = &\ d B_r +b d t, \nonumber \\[0.1in]
d \langle B_l,B_r\rangle(t)\ & = & \ 1_{L(t)=R(t)} \ dt, \nonumber \\[0.1in]
\forall t\ \geq \ T_{R,L}:= \inf\{t>\sigma_L\vee \sigma_R : L(t)=R(t)\}
, & & R(t) \ \geq \ L(t), \label{lrsde}
\end{eqnarray}
where $B^{l}, B^{r}$ are two standard
Brownian motions. We then have the following characterization of the left-right
Brownian web from \cite{SS07}.

\bteo
\label{T:lrwebchar}{\rm (Characterization of the Left-Right Brownian Web).} There exists an $(\mathcal{H}^2, \mathcal{F}_{\mathcal{H}}^2)$-valued random variable $(\W_l,
\W_r)$, called the standard left-right Brownian web (with parameter $b>0$), whose distribution
is uniquely determined by the following two properties:
\begin{itemize}
\item[{\rm (a)}] $\W_l$, resp.\ $\W_r$, is distributed as the standard
Brownian web, except tilted with drift $-b$, resp.~$+b$.

\item[{\rm (b)}] For any finite deterministic set $z_1, \ldots, z_m,
z'_1, \ldots, z'_n\in\R^2$, the subset of paths in $\W_l$ starting from
$z_1, \ldots, z_m$, and the subset of paths in $\W_r$ starting from
$z'_1, \ldots, z'_n$, are jointly distributed as a collection of
left-right coalescing Brownian motions.
\end{itemize}

\eteo

Similar to the Brownian web, the left-right Brownian web $(\W_l, \W_r)$
admits a natural dual $(\hat\W_l, \hat\W_r)$ which is equidistributed
with $(\W_l, \W_r)$ after a rotation
by $180^{\circ}$. In particular, 
$(\W_l,\hat\W_l)$ and $(\W_r, \hat\W_r)$ are
pairs of tilted double Brownian webs.

\bigskip

{\bf Hopping:} The basic idea of the hopping construction
of the standard Brownian net $\Net^b$
consists in concatenating 
paths of the right and left webs. 
More precisely, given two paths $\pi_1, \pi_2\in\Pi$, let $\sigma_{\pi_1}$ and $\sigma_{\pi_2}$ be the starting times of those paths.
For $t >
\sigma_{\pi_1} \vee \sigma_{\pi_2}$ (note the strict inequality),
 $t$ is called
an intersection time of $\pi_1$ and $\pi_2$ if $\pi_1(t)=\pi_2(t)$. By
hopping from $\pi_1$ to $\pi_2$, we mean the construction of a new
path by concatenating together the piece of $\pi_1$ before and the
piece of $\pi_2$ after an intersection time. Given the left-right
Brownian web $(\W_l, \W_r)$, let $H(\W_l\cup\W_r)$ denote the set of paths
constructed by hopping a finite number of times between paths in
$\Wl\bigcup\Wr$.
$\Net^b$
is then constructed as the closure
of $H(W_l\bigcup W_r)$.

\medskip

As pointed out earlier, Sun and Swart gave two other 
constructions of the net based on the pair $(\W_l,\W_r)$. 
In order to describe those two
constructions,
we first 
introduce the notion
of meshes and wedges.

\noindent
{\bf Wedges:} Let $(\hat\W_l, \hat\W_r)$ be the dual left-right Brownian
web almost surely determined by $(\W_l, \W_r)$. For a path $\hat\pi \in
\hat\Pi$, let $\hat\sigma_{\hat \pi}$ denote its (backward) starting
time. Any pair $\hat l\in\hat\W_l$, $\hat r\in \hat\W_r$ with 
$\hat r(\sigma_{\hat l}\wedge \sigma_{\hat r}) < \hat
l(\sigma_{\hat l}\wedge \sigma_{\hat r})$ defines an open set
\be\label{wedge}
W(\hat r, \hat l) = \{(x,u)\in\R^2 : T<u< \sigma_{\hat l}\wedge
\sigma_{\hat r},\ \hat r(u) < x < \hat l(u) \},
\ee
where $T := \sup\{t< \sigma_{\hat l}\wedge \sigma_{\hat r}: \hat
r(t) =\hat l(t)\}$ is the first (backward) hitting time of $\hat r$
and $\hat l$, which might be $-\infty$. Such an open set is called a
{\em wedge} of $(\hat\Wl, \hat\Wr)$.\medskip

\noindent
{\bf Meshes:}  By definition, a {\em mesh} of $(\Wl,\Wr)$ is an open set of
the form
\be\label{mesh}
M=M(r,l)=\{(x,t)\in\R^2:\sigma_l<t<T_{l,r},\ r(t)<x<l(t)\},
\ee
where $l\in\Wl$, $r\in\Wr$ are paths such that $\sigma_l=\sigma_r$,
$l(\sigma_l)=r(\sigma_r)$ and $r(s)<l(s)$ on
$(\sigma_l,\sigma_l+\epsilon)$ for some $\e>0$.  We call
$(l(\sigma_l),\sigma_l)$ the bottom point, $\sigma_l$ the bottom time,
$(l(T_{l,r}),T_{l,r})$ the top point,
 $
 T_{l,r} = \inf\{t>\sigma_l : l(t)=r(t)\}
 $
 the top
time, $r$ the left boundary, and $l$ the right boundary of $M$.\medskip

Given an open set $A\subset \R^2$ and a path $\pi\in\Pi$, we say $\pi$
{\em enters} $A$ if there exist $\sigma_\pi <s<t$ such that $\pi(s)
\notin A$ and $\pi(t)\in A$. We say $\pi$ {\em enters $A$ from the
outside} if there exists $\sigma_\pi <s<t$ such that $\pi(s) \notin
\bar A$, the closure of $ A$, and $\pi(t)\in A$. We now recall the following
characterizations of the Brownian net from \cite{SS07}.

\bteo\label{teo-sun-swart}
Let $\Net^b$ be the standard Brownian net with branching parameter $b$
and
let $(\Wl,\Wr)$
be the left and right webs in $\Net^b$.
\ben
\item $\Net^b$ is identical with the set of paths in $\Pi$ which do
not enter from the outside any wedge of $(\hat\Wl, \hat\Wr)$.
\item $\Net^b$ is identical with the set of paths in $\Pi$ which do
not enter from the outside any mesh of $(\Wl, \Wr)$.
\een
\eteo

\subsubsection{Special Points of the Standard Brownian Net}
\label{special-force}

In this section ,
we will assume no killing in the Brownian net; i.e.,
we will work under the hypothesis 
$k=0$. Under this hypothesis,  
we
recall the classification of the special points
of the Brownian net, as described in \cite{SSS08a}. 
As in the Brownian web,
the classification of special points will be based on 
the local geometry of the Brownian net.
Of special interest
to us will be the points with a deterministic
time coordinate.

\bigskip

Recall the notion of (strong)
equivalence of paths
in the Brownian web, as defined in Section \ref{special_points} --- see (\ref{str-equ}). 
In order to classify
the special points of the Brownian Net,
we will need
to consider the following (weaker) definition
of equivalence between paths entering and leaving a point.

\bdf[Equivalent Ingoing and Outgoing Paths]
\label{eq:in-out}
Two paths $\pi,\pi'\in(\Pi,d)$
are said to be (weakly) equivalent paths entering a point $z\in\R^2$, or in short 
$\pi \sim_{z}^{in} \pi' $, iff there exists a sequence $\{t_n\}$
converging to $t$ such that
$t_n<t$ and
$
\pi(t_n) = \pi'(t_n)
$
for every $n$.
Equivalent paths exiting a point $z$, 
denoted by
$\pi \sim_{z}^{out} \pi' $,
are defined 
analogously by finding a sequence $t_n>t$
with $\{t_n\}$ converging to $t$ and $\pi(t_n)=\pi'(t_n)$. 
\edf

Despite the notation,
these are not in general equivalence relations on the spaces of all
paths entering (resp, leaving) a point.
However, in \cite{SSS08a},
it is shown that
that
if $(\W_l,\W_r)$ is a left-right Brownian web, then a.s. for all
$z\in\R^2$, the relations $\sim_z^{in}$
and
$\sim_z^{out}$
actually
define equivalence relations on the set of paths in $\W_l \cup \W_r$
entering (resp., leaving) z, and the equivalence classes of paths in $\W_l \cup \W_r$
entering (resp., leaving)
$z$ are naturally ordered from left to right. Moreover,
the authors 
gave a complete classification of points
$z\in\R^2$ according to the structure of the equivalence classes in $\W_l \cup \W_r$ entering (resp., leaving)
$z$, in the spirit of the classification of special points of the Brownian web in Theorem \ref{teo:types}.

In general, such an equivalence class may be of three types. If it contains only
paths in $\W_l$ then we say it is of type $l$, 
if it contains only paths in $\W_r$ then we say it is of
type $r$, and if it contains both paths in $W_l$ and $W_r$
then we say it is of type $p$, standing for
pair. To denote the type of a point $z\in\R^2$ in a Brownian net $\Net^b$, 
we first list the incoming
equivalence classes in $W_l\cup W_r$
from left to right and then, separated by a comma, the outgoing
equivalence classes.

In \cite{SSS08a},
the authors showed that there are 20 types of special points in the Brownian net.
However, the next proposition states
that at deterministic times, 
there are only three types of points, 
namely the types (o, p), (p, p)
and (o, pp), where an o means that there are no incoming paths in $\W_l\cup\W_r$ at $z$.
We cite the following result from 
Theorems 1.7 and 1.12 in \cite{SSS08a}.

\bprop[Geometry of the Net at Deterministic Times]
\label{special-net}
Let $\Net^b$ be a Brownian net and let $(\W_l,\W_r)$ be the left-right Brownian web and the dual
left-right Brownian web associated with $\Net^b$.
Let $t$ be a deterministic time.
\ben
\ii For every deterministic point $x$, the point $(x,t)$ is of type $(0,p)$.
\ii Each point $(x,t)$ (with $x$ deterministic or random)
is either of type
(o, p), (p, p) or (o, pp), 
and all of these types occur.
\ii Every path $\pi\in\Net^b$ starting from the line $\R\times\{t\}$,
is squeezed between an equivalent pair of right-most and left-most paths,
i.e.,
there exists $l\in \W_l(z)$ and $r\in \W_r(z)$ so that
$l\sim^z_{out} r$ such that 
$l \leq \pi \leq r$ on $[t,\infty)$.
\ii  Any point $(x,t)$ entered by a path $\pi\in\Net^b$ with $\sigma_\pi<t$
is of type $(p,p)$. Moreover, $\pi$
is squeezed between an equivalent pair of right-most and left-most paths, i.e.,
there exist
 $l\in\W_l, r \in \W_r$ with $r\sim_{in}^z l$
and $\eps > 0$ and such that 
$l \leq \pi \leq r$ on $[t - \eps,t]$.
\een
\eprop

\subsubsection{Branching-Coalescing Point Set}
\label{sec:bc}
We first recall the definition of the
branching coalescing (BC) point set as
as introduced in \cite{SS07}.
To ease the notation, we hide the dependence on $b$.

\bdf{\rm (The Branching-Coalescing Point Set).}\label{df-age}
Let $(x,t)\in\R^2$. If $\sigma_\pi$ denotes the starting time of  a path $\pi$, 
define $\xi^{S}(\cdot)$ (resp., $\xi^{(x,S)}(\cdot)\}$),
the Branching-Coalescing 
point set starting
at time $S$ (resp., at $(x,S)$), as
\beqn
\xi^{S}(T)                    & =  &  
\{y  \ : \   \exists \pi\in\Net^{b}, \  \ \sigma_\pi \ = \ S, \ \pi(T)=y \}, \ \  \text{for} \ \  \   S \leq \  T, \no\\
\xi^{(x,S)}(T)               & =   & \{y  \ : \ \exists \pi\in\Net^{b}, \sigma_\pi = S, \pi(S)=x , \ \pi(T)=y \ \}, \ \ \ \text{for } \ \ S  \ \leq  \ T, \no\\
\xi^{(x,S)}(T), \ \xi^{S}(T) & =  \O & \text{for} \ \  T < S  \label{xis}.
\eeqn
\edf

One of the most striking properties
of the branching-coalescing point set starting at time 
$S$ 
is that the point set
is almost surely locally finite
at any deterministic time 
$T>S$ (even 
if there are
infinitely many paths  
starting from any given point $(x,S)$).
Furthermore, according to
the next proposition cited from \cite{SS07}
for $k=0$, the expected value of the
density of the branching-coalescing
point set
can be computed explicitly.

\bprop
\label{prop-xi-finite}
Let $T$ be a deterministic time. For almost
every realization of the Brownian net $\Net^b$, for every $S<T$,
the set $\xi^S(T)$
is locally finite. Moreover, when $S$ is also deterministic, we have
\beqn
\label{canon}
\E\left( |\xi^S(T)\cap[x_1,x_2]| \right) \ = \ \left(\fr{e^{-b^2(T-S)}}{\sqrt{\pi (T-S)}} \ + \ 2b \phi(b\sqrt{2(T-S)})\right) \cdot (x_2-x_1),
\eeqn
where $\phi$ is the cumulative distribution function of a standard normal r.v.
\eprop

\begin{proof}
For deterministic $S$,
(\ref{canon}) 
is Proposition 1.12 in \cite{SS07} together
with the scaling property 
\beqn
\Net^b = S_{1/b}(\Net^1)  \ \ \ \text{in law} \ \label{scaling-prpoert}
\eeqn
(see Theorem 1.1 of \cite{SS07}). 
For random $S<T$, take a rational number $s$
with $S<s<T$. The result follows by 
noting that $\xi^S(T)\subset \xi^s(T)$
and that $\xi^s(T)$
is a.s. locally finite, by what we have just
stated for the deterministic case.
\end{proof}

\subsubsection{Second Construction of the Brownian Net with Killing}
\label{onenet}

Having the standard Brownian net at hand,
we now need to turn on the killing
mechanism. To do that, we introduce
the time length measure of the Brownian net.

\bprop[Time Length Measure of the Brownian Net]
\label{tlm-net}
For almost every realization of the standard Brownian net,
for every Borel set $E\subset\R^2$, we define
\be
\cL_{net}(E) = \int_\R |\{x : (x,t)\in E\cap(p,p)\}| dt \label{ex-cl},
\ee
where $(p,p)$ refers to the set of $(p,p)$ points for the net $\Net^{b}$.
$\cL_{net}$
defines a $\sigma$-finite measure such that for every $\pi\in\Net^b$ 
\be\label{cL-def}
\forall \ x_1<x_2, t_1<t_2, \ \ \cL_{net}(\{(x,t)\in[x_1,x_2]\times[t_1,t_2]\ :  \ t>\sigma_\pi \ \ \mbox{and} \  \pi(t) =x  \} ) 
\ = \ \int_{t_1\vee \sigma _{\pi}}^{t_2 \vee \sigma_{\pi}} \ 1_{\pi(t)\in[x_1,x_2]} \ dt.
\ee
\eprop
\begin{proof}
This proposition generalizes Proposition \ref{tlm}
and the proof goes along the same lines. As in that proof, 
we consider a countable deterministic dense set ${\cal D}$
and let ${\cal D}_n$ be a nested sequence of finite sets converging to ${\cal D}$.
By definition,
for any point $(x,t)\in(p,p)$ (where $(p,p)$ denotes the set of $(p,p)$
points), there must exist $\pi_l\in\Wl$,
such that 
$\pi_l(t)=x$, with $\sigma_{\pi_l}<t$.
As in the proof of Proposition \ref{tlm}, 
it follows
that
for every path $\pi\in{\cal W}_l$ and any time $t$
such that $\sigma_\pi<t$,
we must have $(\pi(t),t)\in \W_l({\cal D})$ --- where, by a slight abuse of notation, $\W_l({\cal D})$ refers here to the trace 
of the left web paths starting from ${\cal D}$. It follows that
$(p,p)\subset \W_l({\cal D})$ and that
\beqnn
\cL_{net}(E) & = &   \int_{-\infty}^\infty \ |\{x \ : \ (x,t)\in E\cap (p,p) \cap \W_l({\cal D})\}| dt \\ 
 	  & = &    \lim_{n\uaw\infty}\int_{-\infty}^\infty \ |\{x \ : \ (x,t)\in E\cap (p,p) \cap \W_l({\cal D}_n)\}| dt.
\eeqnn 
This proves the $\sigma$-finite property for our measure.
We leave the reader to convince herself that
(\ref{cL-def})
can be proved along the same lines as 
was (\ref{e-cl}), using the structure of the special points of the Brownian net
at deterministic times --- see Proposition \ref{special-net}.
 
\end{proof}

Given a realization of the standard Brownian net $\Net^b$,
we define the set of the killing marks $\Mzero^k$ as
a Poisson point process with intensity measure
$k \cL_{net}$. Finally, we define
the killed Brownian net as the union
of 
(1) all the paths
$\pi\in\Net^b$ killed
at 
\be\label{str-e}
e_{\pi}:=\inf\{t>\sigma_{\pi} : (\pi(t),t)\in\M^k\}
\ee
and (2) for every $z\in\R^2$, the trivial path whose starting point
and ending point coincide with $z$. In the following,
this construction
will be denoted as $\Net^{b,k}_{net}$. Finally, we close this 
section by listing some of the properties
of the marking set $\M^k$.

\bprop\label{mks}
Let $S$ be a deterministic time.
\ben
\ii Let
$\M^k(S)$ be
the set of killing  
points 
which are touched
by the set of paths in the standard 
Brownian net $\Net^b$
with starting time $S$.  The set of points
$\M^k(S)$
is a.s. locally finite.
\ii Let $e$ be a killing time (i.e., the time coordinate 
of some mark in ${\cal M}^k$). For every
rational $\delta t$, the point set $\xi^{e-\delta t}(e)$
is a.s. locally finite.
\een
\eprop
\begin{proof}

We start by proving the first statement. Recall that the set 
of killing points $\M^k$ 
is defined as a PPP with intensity measure 
$k \cL_{net}$. Hence, 
given a realization
of the Brownian net $\Net^b$,
the expected number of
those points 
in a box of the form $[-L,L]\times[S,T]$
is given by
$$
I = k \ \int_{S}^T  |\xi^S( u ) \cap [-L,L]|  \ d u.
$$
Now by (\ref{canon}), the value of this number
averaged over realizations of the Brownian net is finite.
It follows that $\M^k(S)$ is a.s. locally
finite. 

The second statement can be established by the same reasoning, noting that 
the property holds a.s.
for any deterministic $e$ and by definition
of the Brownian net time length measure. 
\end{proof}

\subsection{Equivalence Between the Constructions}
\label{equivalence:construction}
In this section, we show the equivalence between 
the two alternative constructions, $\Net^{b,k}_{web}$ and 
$\Net^{b,k}_{net}$, of the Brownian net with killing. 
As a first step,
we start by showing the compactness of ${\mathcal N}^{b,k}_{net}$.

\bprop
\label{prop:compact}
$\Net_{net}^{b,k}$ is a compact set in $(\Pi,d)$.
\eprop 
\begin{proof}
By compactness of the standard Brownian net, 
the set $\Net_{net}^{b,k}$ is pre-compact
--- by a combination of Arzela-Ascoli
and the fact that the modulus 
of continuity
of a path  is larger than the one of its killed version.
Hence,
we only need to prove that $\Net_{net}^{b,k}$
is closed, which amounts
to proving that 
for every sequence 
$\{\tilde \pi_n\}\in\Net^{b,k}_{net}$
converging to a path $\tilde \pi$, the limiting path also belongs to 
$\Net^{b,k}_{net}$.
The path $\tilde \pi_n$ (resp., $\tilde \pi$)
is characterized by a path $\pi_n\in\Net^b$ (resp., $\pi\in\Net^b$)
and an ending time $e_n$ (resp., $e$), 
so that $\pi_n$ converges to $\pi$
and $e_n$ converges to $e$. In the following, we restrict ourself 
to the case $\sigma_\pi<e$, the other case being obvious
since trivial paths have been included in $\Net^{b,k}_{net}$.
By definition
of the killed Brownian net $\Net_{net}^{b,k}$, the point $(\pi_n(e_n),e_n)$ 
belongs to the set of marks $\M^k$. It remains to show that 
(1) $(\pi(e),e)$ is  a mark and (2) that
this is the first mark encountered by $\pi$ (strictly after time $\sigma_{\pi}$).

First, since $\sigma_\pi<e$, we can always find a rational
time $S$
such that $\sigma_{\pi_n} < S < e_n$ for large enough $n$. For such $n$,
$(\pi_n(e_n),e_n) $ belongs to $\M^k(S)$ --- the set 
of killing points attained by paths starting at $S$ --- with $\M^{k}(S)$ being
locally finite by Proposition \ref{mks}(1).
Hence, after some $n$,
the sequence $(\pi_n(e_n),e_n)$ is fixed
and coincides with a certain $z=(x,e)\in\M^k$, 
which must also be hit by the 
path $\pi$. This proves that $(\pi(e),e)$ is  a mark.

Second, let us assume that $\pi$
passes though a killing point $z=(x',e')$
with $e'$, with $e\geq e'>\sigma_{\pi}$. In order to 
prove that $e=e'$, we need to show that 
$\pi_n$ must also pass through this point. 
Let us take 
$\delta t>0$ rational so that 
$\sigma_\pi<e'-\delta t < e'$. Proposition \ref{mks}(2)
implies that after some $n$,
$(\pi_n(e'),e')$ is stationary in $n$, which implies that   
$\pi_n$ is killed at $z$, as claimed earlier.

\end{proof}

We are now ready to show the equivalence between $\Net^{b,k}_{web}$ and $\Net^{b,k}_{net}$.
We will make 
use of the fact that the result has already been established in the case $k=0$ \cite{NRS08}. 
In the rest of this section,
we will then assume 
that the
Brownian web $\W$ and the standard Brownian net $\Net^b$ are then coupled, 
with the coupling being specified by
the marking construction
described in Section \ref{oneweb}, in the special case $k=0$.

\blem
Under the coupling between $\W$ and $\Net^b$ described above, ~$\cL_{net}=\cL_{web}$.
\elem
\begin{proof}
Recall that the time length measure of the Brownian net
has been defined as
$$
\cL_{net}(E) \ = \ \int_\R |\{x  \ : \  (x,t)\in E \ \text{and} \ \  (x,t)\in (p,p)\}| dt \ ,
$$
whereas for the Brownian web we have
$$
\cL_{web}(E) \ = \ \int_\R |\{x  \ : \  (x,t)\in E \ \text{and} \ \  (x,t)\in (1,1)\}| dt.
$$
Thus, by the Fubini's theorem, it is sufficient to establish that at deterministic times, the set of $(1,1)$ points coincides
with the
set of $(p,p)$
points. In the following $t$
will refer to a deterministic time. 
First, if $z=(x,t)$ is of type $(1,1)$, there exists a path $\pi\in\Net^b$ (in fact a web path)
with $\sigma_{\pi}<t$ so that $\pi$ passes through $z$. By Proposition \ref{special-net}(4),
such a point must be of type $(p,p)$. Conversely, suppose that $z$ is of type $(p,p)$.
By definition there must exist 
$(\pi_l,\pi_r)\in(\Wl,\Wr)$
entering the point $z$ with $\sigma_{\pi_l},\sigma_{\pi_r}<t$.
W.l.o.g., one can assume that those paths 
originate from a dense countable deterministic set ${\cal D}$ and
by the hopping construction of the
Brownian net given in Section \ref{standard-net}, $(\pi_l,\pi_r)$ forms a pair of sticky Brownian motions. 
For such a pair, it is well known that there is no open interval on which the two paths coincide. 
As a consequence,
we can find some $s<t$ such that 
$\pi_l(s)<\pi_r(s)$. Next, we consider a web path $w$  starting from
$\frac{\pi_l(s)+\pi_r(s)}{2}$. Since this path belongs to the Brownian net (recall
that the Brownian net is constructed from the Brownian web through our marking
construction),
our path is squeezed between the 
rightmost and the leftmost path $\pi_r$ and $\pi_l$ and
thus must enter the point $z$.
It follows that $z$ must be a point of type $(a,b)$
in the Brownian web, with $a,b\geq1$.
Finally,
Proposition \ref{special-net}(3)-(4) imply that for a point
of type $(p,p)$, we can only have $a,b\leq1$.
As a consequence, $z$
must be a point of type $(1,1)$. 
\end{proof}

The previous lemma implies that
the two sets
$\M^{k}$ and $\M^{k}_{(1,1)}$ are identical in law and can be coupled, where 
$\M^{k}$ and $\M^{k}_{(1,1)}$ are
respectively the marking points of the Brownian net
and the marked $(1,1)$
points of the Brownian web,
as defined in Sections \ref{onenet} and \ref{oneweb} respectively. 
In the following,
we will assume that $\Net^{b,k}_{net}$
is constructed from the pair $(\Net^{b},\M_{(1,1)}^k)$
and under this assumption,
we will show that $\Net_{net}^{b,k}$
and $\Net_{web}^{b,k}$
coincide (a priori, these might have been distinct because
the orders of killing and taking finite concatenations of web paths
are done differently). First, it is easy to see that 
$\Net_{web}^{b,k} \subset \Net_{net}^{b,k}$. Indeed, by definition,
every path
$\pi\in \Net^{b,k}_{web}$
can be approximated
by a sequence of paths $\pi_n$,
where $\pi_n$ 
is a concatenation 
of web paths killed at some point in $\M_{(1,1)}^k$.
Since $\pi_n$ belongs to $\Net_{net}^{b,k}$, the compactness
of
$\Net_{net}^{b,k}$ (see Proposition \ref{prop:compact}) implies that 
$\pi$ must also belong to $\Netbk_{net}$. 
It now remains to prove that 
$\Net^{b,k}_{net}\subset\Net^{b,k}_{web}$.
First, any trivial path belongs to $\Net^{b,k}_{web}$. 
On the other hand,
for
any non-trivial path $\tilde \pi\in \Net^{b,k}_{net}$,
we can find a pair
$(\pi,e)$ with $\pi\in\Net^b$,
so that the path $\tilde \pi$
is obtained by killing the path $\pi$
at time $e>\sigma_{\pi}$. On the other hand,
by the marking construction of the standard Brownian 
net,
there exists a sequence $\pi_n$
converging to $\pi$, with
$\pi_n$ being constructed by patching together some
web segments at marked $(1,2)$ points. 
By reasoning as in the third paragraph
of the proof of Proposition \ref{prop:compact},
one can easily show that for $n$ large enough,
the $\pi_n$'s
must all be killed at the same time $e$. Hence, the resulting killed path $\tilde \pi_n$,
must
belong to $\Net_{web}^{b,k}$. Since
$\Net_{web}^{b,k}$ is closed (by definition),
the limiting path $\tilde \pi$ must be
also be in $\Net_{web}^{b,k}$, implying
that $\Net_{net}^{b,k}\subset\Net_{web}^{b,k}$. 

\section{Percolation of the Brownian Net with Killing}
\label{Properties}

In this section, 
we prove Theorem \ref{perco}. Recall the definition of 
the BC point set (see Definition \ref{sec:bc}) for the standard Brownian net $\Net^b$. 
Analogously,
we define
$\xi^{k,S}$ (resp., $\xi^{k,(x,S)}$) -- 
the branching coalescing point set with killing starting from $\R\times\{S\}$ (resp., $(x,S)$) -- 
using the Brownian net with killing $\Net^{b,k}$.

\blem
Let $B=\{(x_i,0)\}_{1 \leq i \leq n}$, and let $\xi^{k,B}(t)$ be the BC  point set with killing (with $b=1$) starting from the set $B$ evaluated at time $t$. Finally, let $\tilde \xi^{k,B}(t)$ be the point set obtained by taking the union of $n$
independent BC point sets with killing, the $i^{th}$ point set starting from $\{(x_i,0)\}$.  Then, for every $t\geq0$, the set $\tilde \xi^{k,B}(t)$ stochastically dominates $\xi^{k,B}(t)$.
\elem

\begin{proof}
First, it is not hard to show the discrete analog of this lemma. 
Let us denote by ${\cal U}^{b,k}(B)$
the discrete net starting from $B$.
Let $B:= \{z_1, \cdots z_n\}$ be a subset  $\Z^2_{even}$ consisting of $n$ distinct  points. 
The discrete analog is proved by constructing a natural coupling between ${\cal U}^{b,k}(B)$ -- an object that we will refer to as the joint net -- 
and the $n$ independent discrete nets $\{{\cal U}^{b,k}_i(z_i)\}$, so that 
\begin{equation}\label{coup-3}
\mbox{Trace}({\cal U}^{b,k}(A))\subseteq \cup_i \ \mbox{Trace}({\cal U}^{b,k}_i(z_i)).
\end{equation}

Let ${\cal E}^{b,k}_i, 1 \leq i \leq n$, be the arrow configuration 
underlying ${\cal U}^{b,k}_i(z_i)$, i.e.,  ${\cal E}^{b,k}_i$ is the graphical representation 
of the $i^{th}$ discrete net ${\cal U}^{b,k}_i(z_i)$. We then construct the joint discrete net  according to the following simple rule:
if a site is occupied by one or more independent discrete net particles -- i.e., is a point in the trace of ${\cal U}_i^{b,k}(z_i)$ -- then the joint discrete net particle at that site follows ${\cal E}^{b,k}_M$ where $M$ is the minimum of the set of labels of independent net particles at that site.  
It is easy to see that this coupling has the two desired properties, i.e., that (1) ${\cal U}^{b,k}(A)$ is distributed as a joint discrete net, and (2) (\ref{coup-3}) is satisfied. 

Next, let $x_i \in \R , 1 \leq i \leq n$ and choose a sequence $x_i^\beta$ of even integers, 
so that $S_\beta(x_i^\beta,0)$ converges to $(x_i,0)$ as $\beta\uaw\infty$
and let $B^\beta = \{(x_1^\beta,0),(x_2^\beta,0), \cdots, (x_n^\beta,0)\}$. Finally,
let $b_\beta,k_\beta$
such that $b_\beta e^{\beta}\rightarrow 1$, and $k_\beta e^{2\beta}\rightarrow k$.
The invariance principle for the killed Brownian net  (see Theorem \ref{onewebteo}) implies that 
the coupling  $S_\beta\left({\cal U}^{b_\b,k_\b}(B^\beta), {\cal U}_{1}^{b_\b,k_\b}((x^\beta_1,0)),\cdots,  {\cal U}^{b_\b,k_\b}_n((x^\beta_n,0))\right)$ -- as defined above -- is tight and converges in distribution along subsequences to some measure $\tilde P$. The marginal ${\cal U}^{b_\b,k_\b}(B^\beta)$ converges in distribution (after rescaling) to a ``joint Brownian net" starting from $\left((x_1,0),(x_2,0), \cdots (x_n,0)\right)$, namely ${\cal N}^{b,k}(B)$, while the marginal $\left({\cal U}_1^{b_\b,k_\b}((x^\beta_1,0)),\cdots,  {\cal U}_n^{b_\b,k_\b}((x^\beta_n,0))\right)$ converges in distribution (after rescaling) to $n$ independent Brownian nets $\{\Net^{b,k}_i((x_i,0))\}_i$ starting respectively from  $(x_1,0),(x_2,0), \cdots (x_n,0)$. By the Skorohod representation theorem,
one can assume w.l.o.g. that this convergence is a.s., and on this probability space, we must have 
$$\mbox{Trace}({\cal N}^{b,k}(B)) \subseteq \cup_{i=1}^n \mbox{Trace}({\cal N}_i^{b,k}((x_i,0))).$$
Our lemma then immediately follows.

\end{proof}

{\bf Proof  of Theorem \ref{perco}}. First, an easy coupling argument shows that $p_0^{k}:= P({\cal N}^{1,k} \hbox{ percolates })$ is monotone non-increasing in $k$. (This is obvious at the discrete level, and the
property can be easily extended using the invariance principle of Theorem \ref{onewebteo}.)  
Given the monotonicity of $p_0^k$ in $k$,  to show the existence of a critical killing rate $0 < k_c < \infty$ (for percolation),  it is sufficient to show that $p_0^k = 0$ for some $k < \infty$ and $p_0^k > 0$ for some $k >0$. 

We first  show that $p_0^k =0$ for large enough $k$ using a branching process bound.
In the following, for any  $Y\subset\R$,  we define $(Y,t):=\{(x,t) : \ x\in Y\}$. 
Given a realization of the net $\Net^{b,k}$, let us define inductively the sequence $\{ X_i \}_{i\geq0}$ as follows :
$$
X_0 = 0, \ \ \  X_{n} = \xi^{k,(X_{n-1},n-1)}(n),  \ \ \ \mbox{for $n\geq 1$}, 
$$
where $\xi^{k,B}(t)$ denotes the BC point set with killing at time $t$ starting from the set $B\subset\R^2$.
Note that $X_n$ coincides with $\xi^{k,(0,0)}(n)$ and that we need to show that 
for $k$ large enough $|X_n|=0$ after some $n$.
In order to show that, let us now consider a new process $\{\tilde X_i\}$
defined inductively by taking 
$\tilde X_{0}=0$ and such that 
$$
\tilde X_{n} = \tilde \xi^{k,(\tilde X_{n-1},n-1)}(n),
$$
where  $\tilde \xi^{k,(\tilde X_{n-1},n-1)}(n)$ 
is generated by considering $|\tilde X_{n-1}|$ independent BC point sets with killing starting from the distinct points of 
$(\tilde X_{n-1}, n-1)$ -- those BC point sets with killing being generated by $|\tilde X_{n-1}|$ Brownian nets, independent from
each other, from $\tilde X_{n-1}$ and also from all the nets used to generate the $\tilde X_i$'s, for $i\leq n-1$.

Given $X_{n-1}$ (resp., $\tilde X_{n-1}$), the Markov property
and translation invariance (in time) implies that  the set $\xi^{k,(X_{n-1},n-1)}(n)$ (resp,  $\tilde \xi^{k,(\tilde X_{n-1},n-1)}(n)$)
is identical in distribution with 
$$
\xi^{k,B}(1)  \ \mbox{with $B=(X_{n-1},0)$}  \ (\mbox{resp}.,  \tilde \xi^{k,B}(1)  \ \mbox{with $B=(\tilde X_{n-1},0)$}),
$$
where $\tilde \xi^{k,B}$ is defined as in the previous lemma.
By the previous lemma, this implies that
the sequence $\{|X_i|\}$ is stochastically dominated by 
$\{|\tilde X_i|\}$. Furthermore, 
translation invariance (in space)
implies that
the sequence $\{|\tilde X_i|\}$ defines a Galton Watson branching process 
with offspring distribution 
$|\xi^{k,(0,0)}(1)|$. 
We now show that $\E(|\xi^{k,(0,0)}(1)|)<1$
for $k$ large enough, hence showing a.s. extinction
of $|X_n|$.

First, 
the variable
$
|\xi^{(0,0)}(1)|
$
(corresponding to a BC set with no killing)
stochastically dominates  the random variable $|\xi^{k,(0,0)}(1)|$. Since 
$
\E(|\xi^{(0,0)}(1)|) < \infty
$
-- this can be shown along the same lines as Proposition 1.12 in \cite{SS07} --
the sequence $\{|\xi^{(k,(0,0))}(1)|\}_k$ is uniformly integrable. 
Next, the left and right web paths starting at $(0,0)$ in the net $\Net^b$ (with no killing) interact as sticky Brownian motions. In the interval $[0,1]$, it is well known that the set of contact times
of those two paths has strictly positive 
Lebesgue measure, or equivalently
positive time length measure.
Since any killing along that set would imply that  $|\xi^{k,(0,0)}(1)|=0$, we obtain that 
$\lim_{k\uaw\infty} |\xi^{k,(0,0)}(1)|=0$ in probability.
Since  $\{|\xi^{k,(0,0)}(1)|\}_k$ is uniformly integrable, we get that 
$$
\lim_{k\uaw\infty} \E( |\xi^{k,(0,0)}(1)| ) =0, 
$$ 
hence showing the desired result.

Now we prove that for small enough $k, p_0^k >0$. We proceed as in the discrete case by dynamic renormalization using the analogous result for an oriented site percolation model in $\Z_{even}^2$. 
Here, we only present a sketch of the argument, since the arguments are almost identical in the discrete and continuous case. Further details can be found in \cite{MNR13}. 
We define boxes $B(m,n)$ which are $6m$ wide and $n$ high. The first row of boxes are arranged with their bottom edge along the $\{t=0\}$ axis and $(0,0)$ is at the midpoint of the bottom edge of the box containing the origin. The second row of  boxes are arranged by shifting them by $3m$ so that the midpoints of the bases of these boxes lie above the left or right edge of the box below.  For every point $v$ in the interval $I = [-2m,2m]$ of the bottom edge of the box containing the origin, let $A_v$ be the event that there exists at least one path from $v$ into the set $[-3m,-m]\times \{n \} $ and at least one path into the set $[m,3m] \times \{n\}$,
with both paths staying in the box. Following the approach in the discrete case it is sufficient to prove that $\min_{v \in I} P(A_v) >p_c$, where $p_c$ is the critical value for independent oriented site percolation in $\mathbb{Z}^2$. We now outline an argument for obtaining this result. Consider ${\cal N}^{1,0}$ and set  $n = \frac{m}{4}$. Let $S(-2m)$ be the event $\{ r_{-2m} (n) \in [m,3m] \cap ( \sup_{0 \leq s \leq n} r_{-2m}(s) \leq 3m) \cap ( \inf_{0 \leq s \leq n} r_{-2m} \geq -3m)\}$
 where $r_{-2m}(.)$ is the right web path starting from $-2m$. 
 Similarly define $S(2m )$ to be the event 
 $\{ l_{2m} (n) \in [-3m,-m] \cap( \sup_{0 \leq s \leq n} l_{2m}(s) \geq -3m) \cap (\inf_{0 \leq s \leq n} l_{2m}(s) \leq 3m)\}$ where $l_{2m}(.)$ is the left web path starting from $2m$. It easily follows from the central limit theorem and reflection principle that $P(S(-2m) \cap S(2m)) $ can be made arbitrarily close to one by choosing $n$ large enough. Note that if $S(2m) \cap S(-2m)$ occurs then $l_{2m}$ and $r_{-2m}$ remain in the box $B(n,m)$ and cross at some point $(x,t)$ where $-3m < x < 3m$ and $0 < t < n$ and for any  $v \in I$, $l_v$ will coalesce with $l_{2m}$ at some time $0 < s \leq t$. Therefore if we consider  paths $\pi_v^1$ obtained by following $l_v$ up to time $t$ and then following $r_{-2m}$ and  $\pi_v^2$ obtained by following $l_v$ up to time $t$ and then following $l_{2m}$, we obtain two paths in ${\cal N}^{1,0}$.
This shows that  given $\epsilon >0$, for large enough $n$, $\inf_{v \in I} P(A_v) > 1 -\epsilon/2$. Since the time measure ${\cal T}_{net}(\{ (x,t) \in [-3m,3m] \times [0,n]: \pi(t) =x  \hbox{ for some } \pi \in \cup_{v \in I} A_v\} )$ is finite, it follows that for small enough $k$, $E(e^{-k {\cal T}_{net}(\{ (x,t) \in [-3m,3m] \times [0,n]: \pi(t) =x  \hbox{ for some } \pi \in \cup_{v \in I} A_v\})} ) \geq 1- \epsilon/2$. Thus by choosing $n$ large enough and $k$ small enough we can ensure that $\min_{v \in I} P(A_v) > 1 -\epsilon$. Choosing $\epsilon < 1- p_c$ we obtain our result.

\section{Invariance Principle}
\label{Invariance::Principle}

In this section, we prove Theorem \ref{onewebteo}. Before going into the details of the proof, 
we outline the main steps leading
to our result. First of all, we will show that the set of discrete killing points
converge (in a sense to be made more precise later) to the set 
of continuum killing points. This will be achieved in Section \ref{conv:killing:POINTS}. 
Next, to prove Theorem \ref{onewebteo}, we will combine  this result with a theorem 
in \cite{SS07} stating that
the Brownian net with no killing is the scaling limit of the discrete net (again
with the killing mechanism turned off) , i.e.,
that 
$$
S_{\b}( {\cal U}^{b_\b,0} ) \ \rightarrow \ \Net^b  \ \ \ \mbox{in law},
$$
whenever $e^{\b} b_\b\rightarrow b$ and where $\Net^b$
is a Brownian net
with branching parameter $b$. 

We note that the two previous results combined do not
imply Theorem \ref{onewebteo} by themselves. The main difficulty can be loosely explained as follows.
Let us imagine that a path of the discrete net passes at a microscopic distance
from a killing point, but without passing through it. At the continuum level,
the limiting path passes though the corresponding killing point
and, by our definition of the Brownian net with killing, it would be killed at that point. 
Thus,
one needs to show that a continuum path passing through a killing point $z$
actually
corresponds to a discrete path being killed at a discrete point approximating
the killing point. This will require some technical work, which will be mostly done in Section
\ref{STL}. The proof is divided into five consecutive steps.

\subsection{Step 1 : Convergence of the Set of Killing Points}
\label{conv:killing:POINTS}

In this section, we show convergence of
the set of discrete killing points to its continuum
analog. First of all, we need to define a proper
type of convergence. Recall that the set of killing points is defined
as a Poisson Point Process (PPP) with intensity measure
taken as the time
length measure of the Brownian net --- as defined
in Proposition \ref{tlm-net}. As previously mentioned,
this measure is $\sigma$-finite but the 
measure of any open set $O\subset\R^2$ is
infinite, implying that
the set of killing points is everywhere dense a.s.. 
In the following, we will need to
apply some cut-off procedure to the standard 
Brownian
net and its discrete analog, in such a way that 
the time length measure of the resulting ``reduced" Brownian net
will be finite for any bounded open set $O$. Finally, convergence will be proved under this cut-off.
Our procedure is based on $age(x,t)$,
the age of a space-time point $(x,t)$, defined as
\beq
\label{age}
age(x,t)=\sup\{\delta\geq 0\ :  \ x\in \xi^{t-\delta}(t) \}.
\eeq
The age of a path is defined as the age of
its starting point.
In the following, the $\delta$-shortening of the net $\Net^b$
will refer to the set of
paths obtained from $\Net^b$
after removal of all  paths
whose age (before removal) is less than $\delta$.

Given a realization of the killed Brownian net (with killing rate $k$),
define $K$ to be the point process on $\R^2\times (\R^{+}\cup\{\infty\}\setminus\{0\})$
\be\label{def-K}
K \ :=  \cup_{\eps>0} \{(x,t,\eps) \ : \ (x,t)\in{\cal M}^k  \ \ \text{and} \ \ age(x,t) = \eps  \}, 
\ee
where ${\cal M}^k$ denotes the set of
killing points in the killed 
Brownian net $\Net^{b,k}$.
Analogously, define on $\R^2\times (\R^{+}\cup\{\infty\}\setminus\{0\}$), the point process
$$
K_{\b} := \cup_{\eps>0}\{(x,t,\eps) \ : \ (x,t)\in {\cal M}^k_{\b} \ \ \text{and} \ \ age_{\b}(x,t) = \eps \cdot e^{2\b} \},
$$
where ${\cal M}^k_{\b}$ denotes
the set of killing points for the BCK, killed at a rate $k_{\b}$. (Note that
any discrete point $(x,t,\eps)$ with a third coordinate $\eps>0$, corresponds
to a killing point $(x,t)$ whose  ``macroscopic" age
is equal to $\eps$.) In this section, we will prove the following convergence statement. 

\bprop
\label{p:kill}
If
$k_\b \ e^{2\b} \ \rightarrow k$, then
\beqn
S_\b\left( \cUb , K_\b \right) \ \Longrightarrow \   \left(\Net^b , K\right)  \ \ 
\text{in law},
\eeqn  
where
the convergence of the second coordinate means that the rescaled random measure $\mu_\b ~ := \sum_{(x,t,\eps)\in K_\b} \delta_{(S_\b(x,t),\eps)}$
converges to $\mu = \sum_{(x,t,\eps)\in K} \delta_{(x,t,\eps)}$ on the space of Radon measures on 
 $\R^2\times(\R^{+}\cup\{\infty\}\setminus\{0\})$.
\eprop

\bigskip

Let us now consider
a rectangle of the type $\square = [x_1,x_2]\times[t_1,t_2]$ in $\R^2$. By definition
of the time length measure $\cL$ (see Proposition \ref{tlm}),
the measure of the $\delta$-shortening 
of the net in the rectangle $\square$, is given by
\beqn
\label{continuous-tl}
\cL^\delta(\square) & : = &\int_{t\in[t_1,t_2]} \ \ |\{(x,t) \ : \ \mbox{$x_1\leq x \leq x_2$ and $age(x,t)\geq \delta$} \}| dt \nonumber \\
 & = &\int_{t\in[t_1,t_2]} |\xi^{t-\delta}(t)\cap[x_1,x_2]| \ dt,
\eeqn
and by definition, the set $K$ is distributed as a Poisson Point Process (PPP)
on $\R^2\times (\R^{+}\cup\{\infty\}\setminus\{0\})$
with intensity measure: 
\beqnn
\forall \eps>0, \  \ 
I\left( \square \times [\eps,\infty )\right) \ = k \ \cL^\eps(\square).
\eeqnn
On the other hand,
let us define $\cL_{\b}^\eps$ to be the discrete analog of $\cL^\eps$, i.e., the 
random measure counting 
the number
of points of (macroscopic) age at least $\eps$ in the rectangle $\square$:
\beqnn
\cL_{\b}^\eps(\square)
& := &  \  k_\b \ |\{(x,t)\in\Z_{even}^2  \ :  \ \left(  x,t \right) \in \square \ , \  age_{\b}(x,t)\geq\eps\}|, \label{discrete-tl}
\eeqnn
where $age_\b(x,t)$ is defined analogously to $age(x,t)$, but with respect to the discrete net $\cUb$. 
The point process
$K_\b$ defines a Bernoulli Point Process (BPP) with intensity measure
$k_\b \ \cL^{\eps_\b}_\b(\square)$, with  $\eps_\b:=e^{2\b} \eps$ (where the intensity measure of a BPP is defined
as the average number of points in the set under consideration).
In the rest of this section, we show that this discrete intensity measure converges 
to the continuous one after proper rescaling.

\bprop\label{def-tl}

For every $\square=[x_1,x_2]\times[t_1,t_2]$,
let $|\square|$ be the area of $\square$.
\be
\label{exp:l}
\forall \eps>0, \ \ \E\left(\cL^\eps\left(\square\right) \right)=  \ |\square| \ \left(\fr{e^{-b^2\eps}}{\sqrt{\pi \eps}} \ + \ 2b \phi(b\sqrt{2\eps})\right),
\ee
where $\phi$ is the cumulative distribution function of
a standard normal and $b$
is the branching parameter
of the Brownian net under consideration. 
\eprop

\begin{proof}
We have
\begin{eqnarray}
\E\left(\cL^\eps\left(\square\right) \right) 
& = & \int_{t\in[t_1,t_2]}  \ \E\left( \ |\xi^{t-\eps}(t)\cap[x_1,x_2]| \  \right) dt \\
& = &  \ |\square| \ \left(\fr{e^{-b^2\eps}}{\sqrt{\pi \eps}} \ + \ 2b \phi(b\sqrt{2\eps})\right), \label{prop-1.12}
\end{eqnarray}
where the last equality follows from (\ref{canon}). 
\end{proof}

Next, define $\xi^t_{\b}(\cdot)$ -- called the 
discrete Branching Coalescing (BC) point set starting at time $0$ --
to be 
the discrete analog of $\xi^t(\cdot)$ as defined in 
(\ref{xis}). More precisely,
for $s\leq t\in\N$,
$\xi^s_{\b}(t)$
will refer to 
the discrete BC point set 
valued at time $t$,
generated by all the paths starting on the time line $\R\times \{s\}$, i.e.,
\beq
\label{xis-d}
\xi^s_{\b}(t)  \ = \ \{x \ : \exists \  p_\b \in\cUb \ \ \text{with starting time $\sigma_{p_\b} = s$ s.t. 
$p_\b(t) \ = \ x$} \}. 
\eeq

\blem \cite{SSS13}
\label{c:walk}
Let $\{b_\b\}$ be a sequence of
real positive numbers with $b_{\b} e^{\b} \rightarrow b$ as $\b\rightarrow\infty$. 
Let us consider a family
of one dimensional random walks starting at the origin and characterized by
the following transition probabilities:
\beqnn
P^1_\b  \  =  \  \frac{1}{4}(1+b_\b)^2 \ , \
P^{-1}_\b  \  =  \frac{1}{4}(1-b_\b)^2 \ , \ 
P^0_\b  \  =  \ \frac{1}{2}(1-b_\b^2), 
\eeqnn
where $P^{\Delta}_\b$ is the probability for a $\Delta$-increment. 
If
we denote by
$v_\b$ the first $(-1)$-hitting time 
of the random walk, and if $\{t_\b\}$
is such that
$t_\b e^{-2\b} \rightarrow t$, then
\beqn
\binfty \ e^\beta \ \P(v_\b \ > \ t_\b ) \ = \ \frac{\exp(-b^2 t)}{\sqrt{\pi t}} + 2b\phi(b\sqrt{2t}).
\eeqn
\elem
\begin{proof}
This is part of Lemma 6.13 in \cite{SSS13}.
\end{proof}

\bcor
\label{c:CLT}
Let $(x_\b,t_\b)\in\Z_{even}^2$ be such that $S_\b(x_\b,t_\b)\rightarrow(x,t)$
as $\b\rightarrow\infty$, and $\{\eps_\b\}_\b$ be
such that $\eps_\b e^{-2\b}\rightarrow\eps$. Then
\beqn
\lim_{\b\uparrow\infty} \ e^\b \ \P\left( x_\b \ \in \ \xi_\b^{t_\b-\eps_\b}(t_\b)\right) \ = \  \fr{e^{-b^2\eps}}{\sqrt{\pi \eps}} \ + \ 2 b \phi(b\sqrt{2\eps}),
\eeqn
whenever $b_\b e^{\b} \rightarrow b$.
\ecor

\begin{proof}
By the wedge characterization of the net (at the discrete level) \footnote{Recall that there is no forward path entering a wedge from the outside.},
\beqn
\{x_\b \in\xi^{t_\b-\eps_\b}_{\b}(t_\b)\} 
& = & 
\{ \mbox{$\hat l_\beta$ and $\hat r_\beta$
do not meet on $[t_\b-\eps_\b,t_\b]$} \}, \no
\eeqn
where $\hat l_\b$ and $\hat r_\b$
are respectively the leftmost and rightmost paths (of
the dual discrete net at inverse temperature $\beta$)
starting from $(x_\b+1,t_\b)$ and $(x_\b-1,t_\b)$ respectively,
forming a (discrete) wedge around the point $(x_\b,t_\b)$. 
In particular, if we define
\beqnn
X_\b(u)  & := & \left( \frac{ \hat l_\b-\hat r_\b}{2}(t_\b - u) \right) -1,
\eeqnn
then $X_\b$ is identical in  distribution
to the random walk defined 
in Lemma \ref{c:walk}, which implies 
that
\begin{eqnarray*}
 \lim_{\beta\uaw\infty} \ \ e^{\beta} \ \P(\{x_\b \in\xi^{t_\b-\eps_\b}_{\b}(t_\b)\} )    \ 
& = &
\lim_{\beta\uaw\infty}  e^{\beta}  \ \ \P( \
\min_{u\in\{0,\cdots,\eps_\b\}}  \ \ \frac{\hat l_\b - \hat r_\b}{2}(t_\b-u) \  > \  0 \ )  \\
& = &
\lim_{\beta\uaw\infty}  e^{\beta} \ \P( \
\min_{u\in\{0,\cdots,\eps_\b\}}  \ X_{\b}(u)  \  > \  -1 \ )  \\
& = &
\lim_{\beta\uaw\infty} e^{\b} \P( v_\b > \eps_\b ) \\
& = &
 \fr{e^{-b^2\eps}}{\sqrt{\pi \eps}} \ + \ 2 b \phi(b\sqrt{2\eps}),
\end{eqnarray*}
where $v_\b$ is defined as in the previous lemma.
\end{proof}

\blem
\label{age>eps}
Let $\eps>0$ be deterministic.
Almost surely, the set of times $t$ such that
there exists a point $(x,t)$ with age 
exactly equal to $\eps$ has $0$-Lebesgue measure.

Furthermore, for every deterministic $t$,
for almost every realization
of the Brownian net,
$t$ does not belong to this set.

\elem

\begin{proof}

Let us define ${\mathbb T}$ to be the set of times $t$ such that 
the two following conditions are fulfilled.
\begin{enumerate}
\item At time $t$, any path $\pi$ starting from a point $z$
with time coordinate $t$
is squeezed between an 
equivalent outgoing pair of leftmost and rightmost paths $(l,r)$
(i.e., $\sigma_l=\sigma_r=t$, and there exists 
a  sequence $t_n\downarrow t$ such
that $l(t_n)=\pi(t_n)=r(t_n)$). 
\item There are no $(0,3)$ points for the left web $\W_l$ on the line $\R\times\{t\}$.
\end{enumerate}
By Propositions \ref{special-web} and \ref{special-net},
a deterministic time belongs to those two sets.
Thus, Fubini's Theorem implies the two sets defined 
by those two conditions have full Lebesgue measure,
implying that ${\mathbb T}$ has also full Lebesgue measure almost surely.
Furthermore,
for every deterministic $t$,
$t$ belongs to ${\mathbb T}$ a.s.

Next, let $t\in{\mathbb T}+\eps$. For such a time, we 
now prove that 
there is no point with time coordinate $t$ and with age exactly equal 
to  
$\eps$. Since ${\mathbb T}+\eps$ 
has full Lebesgue measure, this will end the proof of our lemma. In order to prove our claim, we take a point
$\bar z=(\bar x,t)$ with age greater or equal to $\eps$,
and show that this point must have an age strictly greater than $\eps$.
By definition of the age of a point, there must exist a point $z=(x,t-\eps)$
and a path $\pi$ starting from $z$ and passing through $\bar z$.
In order to prove our result, we only need to
construct a path $\pi'$
starting strictly below $t-\eps$ and also passing through $\bar z$. 

By definition of the set  ${\mathbb T}$, and
since $t\in{\mathbb T}+\eps$, 
there exists 
a pair $(l,r)\in(\W_l,\W_r)$ with $l\sim_z^{out}r$
(i.e., there exists a sequence $t_n\downarrow t-\eps$ such
that $l(t_n)=\pi(t_n)=r(t_n)$) 
such that $\pi$ is squeezed between 
$r$ and $l$ at time $t-\eps$. 
Now,
by Lemma \ref{special-web},
and since $z$ is not a $(0,3)$
point for the left web $\W_l$,
the path $l$
can be approximated from below
by a sequence of paths $\{l_m\}$ --- i.e., 
$\sigma_{l_m}<t-\eps$
and $l_m$ coalesces with $l$
at $\tau_m$, with $\tau_m \rightarrow t-\eps$.
Combining this and what has been said in the previous
paragraph,
one can find $m$ and $n$
large enough such that
$$
\tau_m \leq t_n \leq t.
$$
In this
situation,
$l_m$ meets the path $\pi$
at $t_n$ and
by the hopping construction of the standard Brownian
net (see Section \ref{standard-net}),
the path $\pi'$ obtained by hopping
from $l_m$ to $\pi$ at $t_n$
also belongs $\Net^b$.
By construction, the path $\pi'$ starts strictly
below $t-\eps$ and must hit
the point $z$.   This ends the proof
of the lemma.

\end{proof}

We are now ready to show convergence to
the time length measure of the reduced Brownian
net.

\bprop\label{CV_TL_1}
Let $\{\eps_\b\}_\b$ be a sequence such that 
$\eps_\b e^{-2\b} \ \rightarrow \eps$, and let 
$\{\square_\b\}_\b$ be a sequence
of (deterministic) space-time rectangles 
such that $S_\b(\square_\b)$ converges
to a rectangle $\square$.
Then
\beq\label{e:cv_tl_1}
\lim_{\beta\uaw\infty}  \ \left(S_\b(\cUb)  \ , \  e^{-2\b} \cdot \cL_{\b}^{\eps_\b}(\square_\b)\right) \ =  \ 
\left(\Net^b \ , \  \cL^\eps(\square)\right)  \ \ \ \mbox{in law}.
\eeq
\eprop

\begin{proof}
By the Skorohod
Representation Theorem,
we can assume a coupling
under which the 
rescaled discrete net
converges to the Brownian
net
a.s..
Under this coupling, 
it is sufficient to prove that 
\be
\label{eq-inu-0}
\liminf \ e^{-2\b} \cdot \cL_{\b}^{\eps_\b} (\square_\b) \ \geq \ \cL^\eps (\square)
\ee
and then that 
\be
\label{eq-inu-1} 
\lim \E\left( e^{-2 \b} \cdot \cL_{\b}^{\eps_\b} (\square_\b)\right) \ = \ \E\left(\cL^\eps (\square)\right). 
\ee

We start with (\ref{eq-inu-0}). 
Let  us consider the deterministic box $\square=[x_1,x_2]\times[t_1,t_2]$.
We first show that a.s. for every realization of the Brownian net,
there exists a set of full Lebesgue measure  ${\mathbb T}$
such that 
for every $t\in{\mathbb T}$ the following holds : 
for every sequence $\{t_\b\}_\b$ such that 
$t_\b e^{-2\b} \rightarrow t $,
\beq
\label{claim-item}
\liminf \  \left| \ \xi^{t_\b-\eps_\b}_\b(t_\b) \ \cap \ [x_1 e^\b \ , \  x_2 e^\b] \ \right| \ \geq \ 
\left| \ \xi^{t-\eps}(t) \ \cap \ [x_1,x_2] \ \right|. 
\eeq
In order to prove this,
let us first consider the set of times ${\mathbb T}'$
consisting of the $t$'s such that there is no point
with age exactly equal to $\eps$ and with time coordinate $t$. In
particular, Lemma \ref{age>eps} guarantees that 
this set has full Lebesgue measure.
In the following, we will consider the intersection
of this set
with the set of times such that 
$\xi^{t-\eps}(t)$ contains neither
$\{x_1\}$ nor $\{x_2\}$. It is easy to see that
this holds a.s. at any deterministic time, implying that 
the resulting set ${\mathbb T}$ (the
intersection 
of the two sets we just discussed) has also full Lebesgue measure.

We now show that for every $t\in{\mathbb T}$,
if $x\in\xi^{t-\eps}(t) \cap [x_1,x_2]$,
there exists a
sequence
$x_\b e^{-\b} \rightarrow x$
such that
$x_\b\in\xi_\b^{t_\b-\eps_\b}(t_\b) \cap [x_1 e^\b,x_2 e^\b]$
for $\b$ large enough. 
First, by definition, we can find $\pi$,
with $\sigma_\pi\leq t-\eps$
such that
$\pi(t)=x$ and a
path $p_\b\in\cUb$
approximating the path $\pi$. We denote by $x_\b$
the position of that path at time $t_\b$, i.e., $x_\b = p_\b(t_\b)$.
By definition
of the set ${\mathbb T}$, we can assume w.l.o.g.
that 
$\sigma_\pi < t-\eps$,
implying that,
for $\b$ large enough,
we must have $\sigma_{p_\b} \leq t_\b-\eps_\b$ or equivalently
that $x_\b\in\xi_\b^{t_\b-\eps_\b}(t_\b)$.
Since $S_{\b}(x_\b,t_\b) \rightarrow (x,t)$ and 
since  $x\neq x_1,x_2$ (by definition of the set $\mathbb T$), 
we must have
$x_\b\in\xi^{t_\b-\eps_\b}(t_\b) \cap [x_1 e^\b,x_2 e^\b]$,
as previously claimed. This ends the proof of (\ref{claim-item}). 

Let $t_\b=[t e^{2\beta}]$, where $[x]$ denotes the integer part
of $x$ and let us define the piecewise constant function
$$
f_\b(t) := \left| \ \xi^{t_\b-\eps_\b}_\b(t_\b) \ \cap \ [x_1 e^\b \ , \  x_2 e^\b] \ \right|. 
$$
The previous result implies
that that for almost every realization of the Brownian net,
the function  $f_\b$ converges pointwise to 
$$
f(t) := \left| \ \xi^{t-\eps}(t) \ \cap \ [x_1,x_2] \ \right|.
$$
As a consequence,
Fatou's Lemma implies that 
\beqn
\liminf  \ e^{-2\b} \cdot  \cL_{\b}^{\eps_\b} (\square_\b) & = & 
\liminf \ e^{-2\b} \cdot  \sum_{t_\b \in \Z  \ : \ t_\b e^{-2\b}\in[t_1,t_2]}
\left| \ \xi_{\b}^{t_\b-\eps_\b} (t_\b) \cap [x_1 e^\b,x_2 e^\b] \ \right| \no
\\
& \geq & 
 \int_{t_1}^{t_2} \ 
\left| \ \xi^{t-\eps} (t) \cap [x_1,x_2] \ \right| \ dt \no \\
& = &
\cL^\eps (\square).   \no
\eeqn

We now turn to the proof of
(\ref{eq-inu-1}). By translation invariance
we get that
\beqn
\lim \E(e^{-2\b} \cdot \cL_{\b}^{\eps_\b}) (\square_\b)  & = &
\lim \sum_{\substack{(x_\b,t_\b)\in\Z_{even}^2 \\ S_{\beta}(x_\b,t_\b)\in\square}}  \P(\{x_\b \in\xi^{t_\b-\eps_\b}_{\b}(t)\}) \ e^{-2\beta} \no \\
& =&
\lim e^{-3\beta} \left| \{z\in\Z_{even}^2 \ : \  S_{\b}(z)\in\square\} \right|\ \cdot  \  \lim \P( \ \{ 0 \in\xi^{-\eps_\b}_{\b}(0)\}  \  ) e^{\b}. \no
\eeqn
Finally, Corollary \ref{c:CLT} implies that
\beqn
\lim  \ \E( \ e^{-2\b} \ \cdot \ \cL_{\b}^{\eps_\b}) (\square_\b)  & = & 
\lim e^{-3\beta} \left|\{ z\in\Z_{even}^2 \ : \  S_{\b}(z)\in\square \}\right|\ \cdot  \  \left( \fr{e^{-\eps}}{\sqrt{\pi \eps}} \ + \ 2 \phi(\sqrt{2\eps}) \right).  \no \\
& = &
|\square| \ \left(\fr{e^{-b^2 \eps}}{\sqrt{\pi \eps}} \ + \ 2 b \phi(b \sqrt{2\eps})\right)=\E( \cL^\eps (\square) )  \no
\eeqn
by Proposition \ref{def-tl}.
\end{proof}

To close this section 
with the promised convergence
result for the set of killing points,
we  use the following easy lemma cited
from \cite{SSS13}  --- see Lemma 6.17 therein.

\blem{\rm (Weak Convergence of Coupled Random Variables).}
\label{JAN}
Let $E$ be a Polish space, let $\{F_i\}_{i\in I}$ be a finite or countable collection of Polish spaces
and for each $i \in I$, let $f_i : E \raw F_i$ be a measurable function. Let $X,X_{\b}, Y_{\b,i}$ be 
random
variables such that $X,X_\b$ take values in E and $Y_{\b,i}$ takes values in $F_i$.
Then
$$
\forall i \ \  \{(X_\b, Y_{\b,i})\}_\b \Rightarrow (X, f_i(X)) \ \ \mbox{in law}
$$
implies 
$$
\{(X_\b, \{Y_{\b,i}\})_{i\in I})\}_\b
\Rightarrow
(X, \{f_i(X)\}_{i\in I})  \ \ \text{in law}.
$$
\elem

\bigskip

{\bf Proof of Proposition \ref{p:kill}.}
By the Skorohod Representation Theorem and Lemma \ref{JAN},
we can construct a coupling such that 
(\ref{e:cv_tl_1}) holds jointly for every rational $\eps>0$ and
every rational rectangle $\square$ a.s..
By definition of the marking ${\cal M}^k$,
we already noted that
the process $K$
defines a PPP on $\R^2\times(\R^{+}\cup\{\infty\}\backslash\{0\})$, 
whose intensity measure is given by:
\beqn\label{I(B)}
\forall \eps>0, \  \ 
I\left( \square \times [\eps,\infty )\right) \ = k \ \cL^\eps(\square),
\eeqn
where $\square\subset\R^2$ is a 
rectangle and $k$
is the killing rate in the Brownian net. 
Analogously, the rescaled
point process $S_\b(K_{\b})$ is a BPP on $\R^2\times(\R^{+}\cup\{\infty\}\backslash\{0\})$
with intensity measure:
\beqnn
\forall \eps > 0,  \  \
I_{\b}\left( \square \times [\eps,\infty) \right)  & := &  
k_{\b}  \ |\{ z \in \Z_{even}^2 : age_{\b}(z) \geq  \eps \ e^{2\b}  \ , \   S_\b(z) \in \square \}|  \\
& = & k_{\b} \ e^{2\b} \  \cdot \ e^{-2\b} \ \cL^{\eps_\b}_{\b}\left( \ S_\b^{-1}(\square) \ \right), \   \ 
\text{and} \ \ \eps_\b \ =  \ \eps \cdot e^{2\b}. 
\eeqnn
Proposition \ref{CV_TL_1} implies that for every rational $\eps>0$ and for every rational rectangle $\square$, we 
have
$$
I_{\b}\left( \square \times [\eps,\infty) \right)  \ \ \rightarrow  \ \ I\left( \square \times [\eps,\infty )\right) \ \ \mbox{a.s.}
$$
It is then easy to establish (by standard arguments) that a.s. for any realization of the net, the same result holds
for any relatively compact Borel set $B\subset\R^2\times(\R^+\cup\{\infty\}\setminus\{0\})$ (random or not) such that $I(\partial B)=0$. 
This easily implies that, conditional on a realization of the Brownian net, for any such set $B$,
$$
|S_\b(K_\b)\cap B| \rightarrow |K\cap B| \ \ \mbox{in law.}   
$$
This concludes the proof of our proposition (see e.g., Proposition 1.22 in \cite{K91}).

\subsection{Step 2 : Relating Discrete and Continuum Killing Events}
\label{STL}

In the following, $L$ will refer to a fixed (but arbitrarily)
large positive real number.
We define 
\be\label{thetaL}
\theta^{L} = \inf\{t > 0 \ : \ \exists x\in\xi^{0}(t)\cap[-L,L] \ \mbox{s.t.} \ (x,t)\in\M^k \}.
\ee
and its discrete analog 
$$
\theta_\beta^{L} = \inf\{t > 0 \ : \ \exists x\in\xi_\beta^{0}(t)\cap[-L_\b,L_\b] \ \mbox{s.t.} \ (x,t) \ \mbox{is a killing point.} \}
$$
where $L_{\b} = [L e^\b]$.
To ease the notation, we will sometimes drop the superscript $L$,
hence writing
$\theta^{L}\equiv\theta$ and $\theta^{L}_\b\equiv\theta_\b$. In the following, we consider an arbitrary 
sequence $\{\beta_n\}_{n\in\N}$ with $\lim_{n\uaw\infty} \beta_n =\infty$.
By a slight abuse of notations, 
we will drop the subscript $n$ and for any sequence $\{X_{{\beta_n}}\}_n$, 
we will write  $\{X_{\b_n}\}=\{X_\b\}$ and $\lim_{\b\uaw\infty} X_\beta \equiv \lim_{n\uaw\infty} X_{{\beta_n}}$.

\bigskip

By the Skorohod Representation
Theorem and Proposition
\ref{p:kill},
we can construct a coupling 
between $\{(\cUb,K_\b)\}$ and $(\Net^b,K)$
such that
\be\label{coucou}
S_\b( \cUb,K_\b ) \ \rightarrow \  (\Net^b,K) \ \ \text{a.s.}
\ee
From now on, we will work under this coupling and we will prove the following result:

\bprop\label{kill-kill}
Let us consider 
the set of paths $\pi\in\Net^b$ (net with no killing) such that 
$\sigma_\pi=0$, $age(\pi)>0$ and $\max_{t\in[0,\theta^L]} |\pi(t)|<L/2$.
Under the coupling (\ref{coucou}), there is a subsequence -- call it $\{\tilde \beta_n\}$ -- such that 
the following properties hold a.s..

\begin{enumerate}
\item $\{e^{-2\bs}\theta_\bs\}$ converges to $\theta$ a.s..
\item For every $\pi$ in the set,
there exists a sequence $\pi_{\bs_n}\in{\cal U}^{\bs_n}$ 
such that (i)
$
\lim_{\bs\uaw\infty}d(\pi,S_\bs(\pi_\bs)) = 0
$
uniformly in $\pi$, and (ii) such that 
for $\tilde \beta$ large enough (again uniformly in $\b$)
$$
\pi   \ \mbox{is killed at $\theta^L$} \  \Longleftrightarrow  \ \pi_\bs \  \mbox{is killed at $\theta_\bs^{L}$.} \
\footnote{recall that  a path $\pi\in\Net^b$ is said to be killed at $z$ iff $z$ is the first killing mark
encountered by $\pi$.}$$ 
\end{enumerate}
\eprop

\bigskip

This will be shown by proving 
that at time $\theta_\b$,
the branching coalescing point set 
is sparse -- see Lemma \ref{l:sparse} below -- i.e.,
that there are no points at a microscopic distance from
one another.
In particular, this will imply that
that if $\pi$
passes through a killing point
at $\theta$, then 
the analog must occur at
time $\theta_\b$
for the path $\pi_\b$.

\blem\label{thetab}
$\theta>0$ a.s. and there exists a.s. a unique $x$ such that $(x,\theta)\in{\cal M}^k$. Furthermore, writing $L_\b=\left[e^\b L\right]$
$$
S_\b(\sum_{\begin{array}{c} x'\in \xi_{\b}(\theta_\b)\cap[-L_{\b},L_\b] \\ (x',\theta_\b)\ \mbox{is discrete killing point} \end{array} } \delta_{x'},\ \theta_\b) \ \rightarrow \ (\delta_{x} ,\theta) \ \ \ \mbox{in law.}
$$
 \elem

\begin{proof}
The fact that $\theta>0$ follows directly from Proposition \ref{mks}(1).
First, the a.s. uniqueness of $(x,\theta)$ follows directly from the properties of $K$. Next, note that 
\beqnn
\theta \ = \ \inf\{t>0 \ : \ (x,t,\eps)\in K, \ \ \eps \geq t , \  |x|\leq L\}, \\
\theta_\b \ = \ \inf\{t>0 \ : \ (x,t,\eps)\in K_\b, \ \ \eps \geq t, \   |x|\leq L e^{\b}\}.
\eeqnn
For every $\gamma>0$, let us introduce 
\beqnn
\theta(\gamma) \ = \ \inf\{t>\gamma \ : \ (x,t,\eps)\in K,\  \eps \geq t, \   |x|\leq L\}, \\
\theta_\b(\gamma) \ = \ \inf\{t>\gamma e^{2\b} \ : \ (x,t,\eps)\in K_\b, \ \eps \geq t,  \  |x|\leq L e^{\b}\}.
\eeqnn
We first claim that our lemma holds replacing 
$\theta$ (resp., $\theta_\b$) by $\theta(\gamma)$ (resp., $\theta_\b(\gamma)$).
First, by the choice of our coupling, 
the (rescaled) point set $K_\b$ converges to $K$. 
Secondly, the set 
$$
B:=\{(x,t,\eps)\in K : \ t>\gamma, \   \eps \geq t, \   |x|\leq L\}
$$
is a Borel set of $\R^2\times(\R^+\cup\{\infty\}\setminus\{0\})$.  In order to prove the claim made earlier, it
suffices to prove that 
$I(\partial B)=0$ a.s., where $I$ is the measure on $\R^2\times(\R^{+}\cup\{\infty\}\setminus\{0\})$ defined  in (\ref{I(B)}). 
First, the set $\{(x,t,\eps) : \ |x|=L, \ \eps>\gamma\}$ 
has no mass -- there cannot be any mark with spatial coordinate exactly equal to $L$ or $-L$.
Furthermore,
Lemma \ref{age>eps} (see second paragraph therein), guarantees that for any deterministic time $t$,
the age of a point of the form $(x,t)\in(\xi^0(t),t)$ must be strictly greater than $t$. By Fubini's Theorem,
this directly implies that
$I(\{(x,t,t) \ : \ t\geq \gamma\})=0$. 
This proves the claim made earlier, i.e., that our result holds if we replace 
$\theta$ (resp., $\theta_\b$) by $\theta(\gamma)$ (resp., $\theta_\b(\gamma)$).

In order to conclude, we need to prove that $\theta$ and $\theta(\gamma)$ (and their discrete analogs)
coincide with high probability for $\gamma$ small enough, i.e., we will show that 
$$
 \lim_{\gamma\daw0} \P(\theta<\gamma)=0 , \ \  \mbox{and} \  \  \lim_{\gamma\daw0}  \limsup_{\beta\uaw\infty}\P(e^{-2\beta}\theta_{\b} <\gamma)=0. 
$$
The first identity is a direct consequence of the fact that $\theta>0$ a.s.
We now show
the second identity. 
Let us define $\cL_{\b,L}(\gamma)$ to be the time length
measure accumulated by the set of discrete paths starting from
the origin
in the macroscopic box $[-Le^\b,Le^\b]\times[0,\gamma e^{2\b}]$,
i.e.,
$$
\cL_{\b,L}(\gamma) \ = \ \sum_{t\in \N:  t\leq \gamma e^{2\b}}  \ |\{x\in\xi^0_{\b}(t) : |x|\leq Le^\b \}|.
$$
For every $\delta>0$,
\beqnn
\P( e^{-2\beta} \theta^\b <\gamma )
 & \leq & \P(\cL_{\b,L}(\gamma) \geq e^{2\b}  \delta) +
 \P( e^{-2\b} \theta_\b \leq \gamma,  \  \cL_{\b,L}(\gamma) \leq e^{2\b}  \delta ) \\
 & \leq & 
  \P(\cL_{\b,L}(\gamma) \geq e^{2\b} \delta) +
1-(1-k_\b)^{\delta e^{2\b}}  
\eeqnn
implying that 
\beqn
\limsup_{\b\uaw\infty} \P( e^{-2\beta} \theta^\b <\gamma )
 & \leq &
\limsup_{\b\uaw\infty} \P(\cL_{\b,L}(\gamma) \geq e^{2\b}  \delta) +
1- \lim_{\b\uaw\infty} (1-k_\b)^{\delta e^{2\b}}
 \no\\
 & \leq & 
\limsup_{\b\uaw\infty}  \P(\cL_{\b,L}(\gamma) \geq e^{2\b}  \delta)  + 1-\exp(-k \delta)\no\\
 &\leq&
\limsup_{\b\uaw\infty}  \P(\cL_{\b,L}(\gamma) \geq e^{2\b}  \delta)  + k \delta \label{yes1}.
\eeqn
Finally, Chebytchev's inequality
implies that
\beqn
\limsup_{\b\uaw\infty}   \P(\cL_{\b,L}(\gamma) \geq e^{2\b}  \delta) & \leq &
\limsup_{\b\uaw\infty}  \frac{e^{-2\b}}{\delta}  \E(\cL_{\b,L}(\gamma) ) \no \\ 
& = &\limsup_{\b\uaw\infty}  \frac{e^{-2\b}}{\delta} \sum_{t\in \N: t\leq \gamma e^{2\b}}  \E(  |\{ x\in \xi^0_{\b}(t) \  :  |x|\leq Le^\b \}| )  \no \\
& = & \limsup_{\b\uaw\infty}  \frac{e^{-2\b}}{\delta} \sum_{t\in \N: t\leq \gamma e^{2\b}} \sum_{x : |x|\leq L e^\b} 
\P(  \{ x \in\xi^0_{\b}(t) \} ). \label{yes2}
\eeqn
Using translation invariance,
for every $(x,t)\in\Z_{even}^2$,
\beq
\P\left( x \in \ \xi^{0}_{\b}(t) \right) \ = \ \P\left( x_0(t) \in \ \xi^{0}_{\b}(t)  \right), 
\eeq
where $x_0(t)=0$ if $t$ is even, and is equal to $1$ otherwise.
By Corollary \ref{c:CLT}, for any sequence $\{t_\b\}$ 
such that $e^{-2\beta} t_\b\rightarrow t$,
we have
$$
\limsup_{\b\uaw\infty}  e^{\b} \P\left( x_0(t) \in \ \xi^{0}_{\b}(t)  \right)  = \  \fr{e^{-b^2t }}{\sqrt{\pi t}} \ + \ 2 b \phi(b\sqrt{2 t}). 
$$
Using the approximation
of an integral by a Riemann sum, (\ref{yes2}) implies that 
\beqnn
\limsup_{\b\uaw\infty}   \P(\cL_{\b,L}(\gamma) \geq e^{2\b}  \delta)  
& \leq & \limsup_{\b\uaw\infty}  \frac{e^{-2\b}}{\delta} \sum_{t\in \Z: t\leq \gamma e^{2\b}}   L e^{\b} \P(  \{ x_0(t) \in\xi^0_{\b}(t) \} ) 
\\
& = & 
\frac{L}{\delta} \int_{0}^\gamma [\fr{e^{-b^2 t}}{\sqrt{\pi t}} \ + \ 2 b \phi(b\sqrt{2t}) ]dt.
\eeqnn
As a consequence, from (\ref{yes1}),
we get that 
$$
\limsup_{\b\uaw\infty} \P( e^{-2\beta} \theta^\b <\gamma ) \leq \frac{L}{\delta} \int_{0}^\gamma 
[\fr{e^{-b^2 t}}{\sqrt{\pi t}} \ + \ 2 b \phi(b\sqrt{2 t})] d t \ + \ k\delta
$$
for arbitrary small $\delta$ and $\gamma$. Letting $\gamma\rightarrow0$ and then $\delta\rightarrow0$
yields our result.

\end{proof}

In the following lemma, we compare 
the discrete BC point set issued 
from (1) the line $\Z_{even}\times\{0\}$ (with $\Z_{even}= \{x \in\Z  \ : \ \mbox{x is even} \}$),  and (2)  a Bernoulli
field whose intensity is ``finite at the macroscopic level''. 
Clearly, the first one must stochastically dominate the other. The next result shows 
that when the intensity of the Bernoulli field is high enough,
the two point sets at a given positive (macroscopic) time are equal with high probability.

\blem
\label{sparse:first}
Let $\delta>0$ and let $\{a_n\}$ be a positive sequence
diverging to $\infty$. For every $n$, we consider 
a sequence $a_\b^n\in[0,1]$ such that 
$a_\b^n e^\b \rightarrow a_n$, and we denote 
by $\{A_\b^n\}$ the
Bernoulli random field on $\Z_{even}$,
with
intensity $a_\b^n$.
Let us denote by $\xi^{A_\b^n}_\b$
the branching-coalescing point set
starting on the set of
points $A_\b^n$ at time $0$, and
\beqnn
G_\b^n & : =  & \{ \xi_{\b}^0(\delta e^{2\b}) \cap [-L e^{\b}, L e^\b]  \neq \xi_{\b}^{A_\b^n}(\delta e^{2\b})  \cap [-L e^{\b}, L e^\b] \}.
\eeqnn
Then
\beqnn
\lim_{n\uaw\infty} \limsup_{\b\uaw\infty} \P\left( G_\b^n \right) & = & 0.
\eeqnn
\elem

\begin{proof}
We will assume w.l.o.g. that 
$e^{2\b}\delta$ is an even integer.
Using translation invariance of the discrete net, we have
\beqnn
\P(G_\b^n) & \leq & L e^{\b}  \P(\hat r_\b < \hat l_\b \ \mbox{on $[0,e^{2\b}\delta]$}, \ \ \mbox{and} 
\ \ A_\b^n \cap [\hat r_\b(0),\hat l_\b(0)] =\emptyset ) 
\eeqnn
where $\hat r_\b,\hat l_\b$ are the right-most and left-most dual paths 
respectively from $-1$ and $+1$ at the (even) time $e^{2\b}\delta$ ---
so that $(\pm 1,e^{2\b}\delta)$ is in $\Z_{odd}^2$.
Thus, for every $\gamma>0$, we have
\beqnn
\P(G_\b^n)  & \leq &
 L e^{\b} \cdot \ [ \  \P(\hat r_\b < \hat l_\b \ \mbox{on $[0,\delta e^{2\b}]$}, \ \ \mbox{and} 
\ \ |\hat r_\b(0)-\hat l_\b(0)| <\gamma e^\b)  \\ 
& + & \ \P(\hat r_\b < \hat l_\b \ \mbox{on $[0,e^{2\b} \delta]$}, \ \ \mbox{and} 
\ \ |\hat r_\b(0)-\hat l_\b(0)| \geq e^\b \gamma, \  \ A_\b^n \cap [\hat r_\b(0),\hat l_\b(0)] =\emptyset )  \   ]
\\
& \leq &
 L e^{\b} \cdot \ [ \  \P(\hat r_\b < \hat l_\b \ \mbox{on $[0,e^{2\b} \delta]$}, \ \ \mbox{and} 
\ \ |\hat r_\b(0)-\hat l_\b(0)| <e^\b \gamma)  
 \ +   \ \P(\hat r_\b < \hat l_\b \ \mbox{on $[0,e^{2\b}\delta]$}) \cdot (1 -  a_\b^n)^{e^\b\gamma}   \   ].
\eeqnn
By Lemma \ref{c:walk},
\beqnn
\lim_{\b\uaw\infty}  e^{\b} \P(\hat r_\b < \hat l_\b \ \mbox{on $[0,e^{2\b}\delta]$}) \cdot (1 -  a_\b^n)^{e^\b\gamma} 
\leq  \left[ \fr{e^{-b^2\delta}}{\sqrt{\pi \delta}} \ + \ 2 b \phi(b\sqrt{2\delta}) \right] \cdot  \exp(-a_n \gamma ). 
\eeqnn
Thus, letting $n\uaw\infty$ on both sides 
of the inequality yields that 
$$
\lim_{n\uaw\infty} \limsup_{\b\uaw\infty}  e^{\b}   \ \P(\hat r_\b < \hat l_\b \ \mbox{on $[0,e^{2\b}\delta]$}) \cdot (1 -  a_\b^n)^{e^\b\gamma} = 0.  
$$
On the other hand,
\beqnn
y(\gamma) \ 
& := & \ \limsup_{n\uaw\infty} \limsup_{\b\uaw\infty} e^{\b}   \cdot \ \  \P\left(\hat r_\b < \hat l_\b \ \mbox{on $[0,e^{2\b}\delta]$}, \ \ \mbox{and} 
\ \ |\hat r_\b(0)-\hat l_\b(0)| <e^\b \gamma \right) \\
& = &
 \left[ \fr{e^{-b^2\delta}}{\sqrt{\pi \delta}} \ + \ 2 b \phi(b\sqrt{2\delta}) \right]  \cdot \ \limsup_{n\uaw\infty} \limsup_{\b\uaw\infty}   \cdot \ [ \  \P( |\hat r_\b(0)-\hat l_\b(0)| <e^{\b}\gamma \ | \ \hat r_\b < \hat l_\b \ \mbox{on $[0,e^{2\b}\delta]$} )]  
\eeqnn
and one needs to show that  the limit
on the RHS of the second equality goes to $0$ as $\gamma$ goes to $0$. 
In other words, one needs to show that, knowing that the two walks separate 
separate for a macroscopic interval of time  $[0,\delta]$,
they must be at a non-negligible macroscopic distance at time $\delta$.  
This can be established by standard arguments.

\end{proof}

Using this result, we next show that at low temperature (i.e., for large $\b$),
a discrete killing point encountered by some (discrete path) starting at time $0$,
is very likely to be isolated. In order to quantify isolation,
we
introduce 
a measure of the sparseness of the branching-coalescing
point set $\xi^{0}_{\b}(t)$ in the interval $[-L,L]$ :
$$
D_\b(L,t) \ = e^{-\b} \  \min_{ \substack{x \neq y \\ x,y \in \xi_{\b}^{0}(t)\cap[-L ,L ] }}  | x - y |. 
$$ 
The next lemma states that for large $\b$,
the branching-coalescing point set starting at 
time $0$ is sparse at time $\theta_\b$.

\begin{figure}
\centering
\includegraphics[scale=.5]{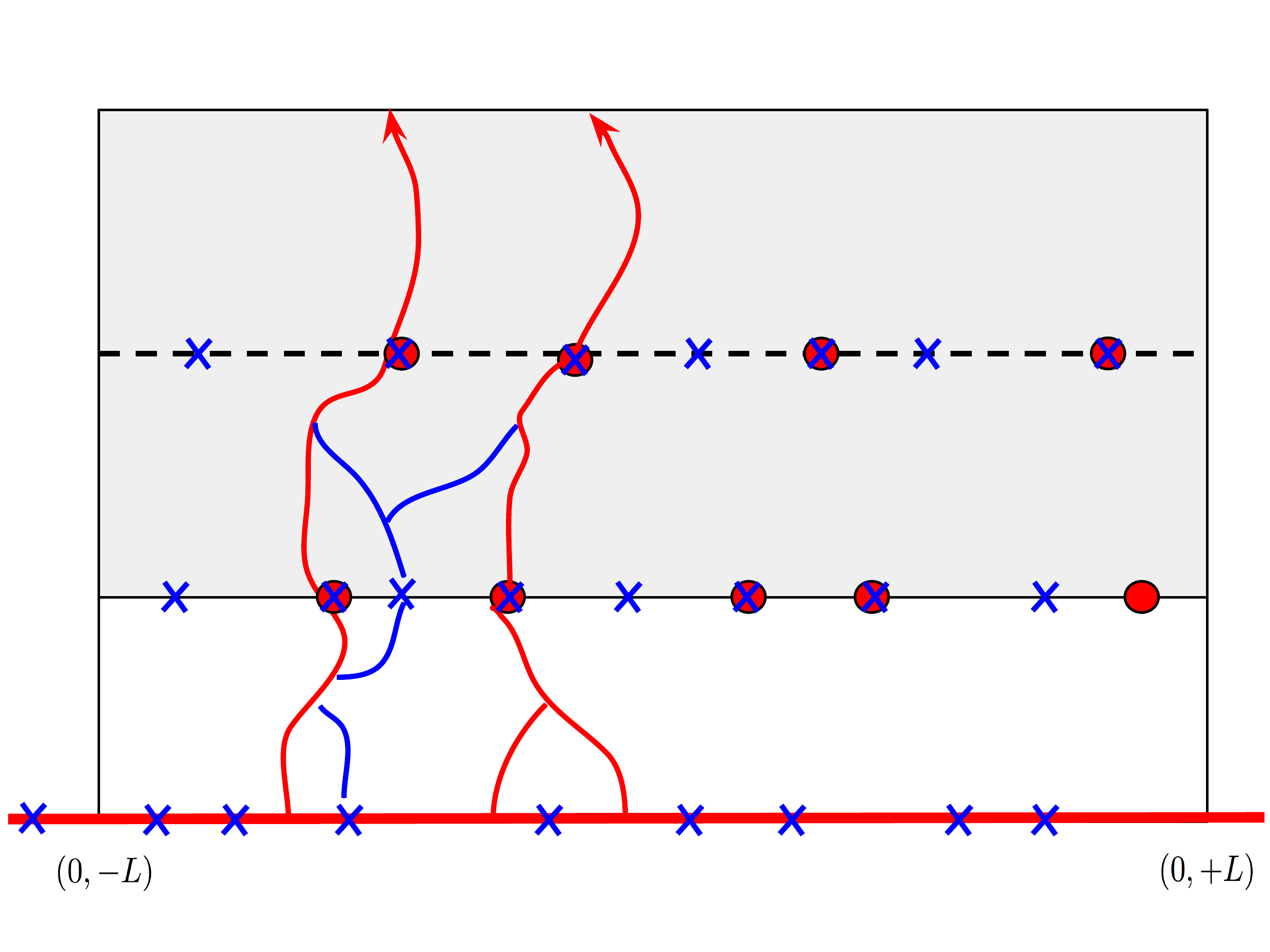}
\caption{Schematic diagram for two coupled discrete BC point sets with different branching rates. The high branching rate BC point set  (crosses) starts from a Bernoulli field at time $0$. The low branching rate BC point set starts from the whole line (red line). 
With high probability, from $T^{-}$ to $T^+$, the blue BC point set dominates the red one on the interval $[-L,L]$. The blue crosses 
along any horizontal line form a Bernoulli field.}
\label{dessin2}
\end{figure}

\blem
\label{l:sparse}
\be\label{lmlm}
\limsup_{\eps\daw0}  \ \limsup_{\b\uaw\infty}  \ \ \P\left( \ D_\b(L e^\b,\theta_\b) \ < \ \eps \right) \ = \ 0.  
\ee

\elem

\begin{proof}

Let $T^{-},T^{+}>0$ with $T^{-}<T^{+}$ and define $\square_{L,T^\pm}=[-L,L]\times[T^{-},T^{+}]$. 
Let us also define the event
$B_{\eps,\b,L,T^\pm}$ to be the event
that a path of the discrete net (with no killing) starting from time $0$ passes through
a killing point in 
the macroscopic box $\square_{L,T^\pm}$ and such that, at that time,
the particle is at a macroscopic
distance less than $\eps$
from another particle starting from the origin, i.e.,
$$
B_{\eps,\b,L,T^\pm} 
\ = \ \{\exists \ x\neq x',t\in\N : S_{\b}(x,t) \in \square_{L,T^\pm}, (x,t)\in{\cal M}^{k_\b}_{\b}, \ x,x' \in \ \xi^{0}_{\b}(t), 
\ |x-x'| < e^{\b} \eps \ \},
$$
where ${\cal M}^{k_\b}_{\b}$ denotes the set of discrete killing points.
We have
\beqn
\P( D_\b(L e^\b,\theta_\b) \ < \ \eps )  & \leq &  \P( B_{\eps,\b,L,T^\pm})
+ \P\left( e^{-2\b} \theta^{L_\b}_{\b} \geq T^{+} \right) 
+ \P\left( e^{-2\b} \theta^{L_\b}_{\b} \leq T^{-} \right) .
\label{dfg}
\eeqn
In the remainder of the proof,
we will first show that for every $L,T^-,T^+>0$,
\beqn
\label{e:b}
\lim_{\eps\daw0} \ \  \  \ \limsup_{\b\uparrow\infty} \ \P( B_{\eps,\b,L,T^\pm} )  \ \ = \ \ 0.
\eeqn
Then, the result will follow by noting that when we let $T^{+}$ and $T^-$ go
to infinity and $0$ respectively, the 
two remaining terms in the RHS of (\ref{dfg}) both go to $0$. Indeed,
since $e^{-2\b} \theta^{L_\b}_{\b}$ converges to $\theta^L$ 
(by Lemma \ref{thetab}), and
$\theta\in(0,\infty)$, we
get,
$$
\lim_{T^+\uaw\infty}\limsup_{\b\uaw\infty} \ \P\left( e^{-2\b} \theta^{L_\b}_{\b} \geq T^{+} \right) \ \ =
 \ \ \lim_{T^-\daw0}\limsup_{\b\uaw\infty} \P\left( e^{-2\b} \theta^{L_\b}_{\b} \leq T^{-} \right)  \ =  \ 0. 
$$
We now turn to the proof of (\ref{e:b}), which will be shown
using a coupling method. Before going into
technical details, we start by giving a brief outline of the
argument (see also Fig. \ref{dessin2}).
As we shall see below, estimating  $B_{\eps,\b,L,T^\pm}$
requires 
obtaining an upper-bound for the probability of finding two elements of the branching-coalescing BC
point set (starting on the line $\{t=0\}$) in a given interval of macroscopic size $\eps$,
i.e.,  we will need to estimate the
correlation function
\be\label{corr-ff}
\P\left( x,x' \in \ \xi^{0}_{\b}(t) , \ \ \mbox{ for a pair $(x,x')$ s.t. $|x-x'|<e^{\b} \eps$}\right),
\ee
for any $t\in[T^-,T^+]$.
In general, the consecutive points 
of the BC are correlated in a nontrivial way. However, in 
\cite{SS07},
it is shown that in the particular 
case where
one starts the BC point set at time $t=0$ 
from a Bernoulli point set $A_\b$ on $\Z_{even}$
with intensity
$$
 \frac{4 b_\b}{(1+b_\b)^2},  
$$
then 
$(\xi^{A_\b}_\b(n), n\in\N)$ -- the BC point set starting from $\{(0,x) :  x\in A_\b\}$ -- 
defines a stationary process modulo parity 
\footnote{Recall the BC set takes values in $\Z_{even}$ for $t$ even and $\Z_{odd}$ for $t$ odd. We can make the process $(\xi^{A_\b}_\b(n), n\in\N)$  stationary by translating the BC at odd times by one unit of space to the right}. Unfortunately, we see from (\ref{corr-ff}), that
we need to consider the 
probability of finding the correlation function at time $t$ 
for the BC set starting from the entire line $\{(x,t)\in\Z_{even}^2 \ : \ t=0\}$. 
However, Lemma \ref{sparse:first} will guarantee
us that 
if we start with 
a Bernoulli field $\bar A_\b$
with high intensity
\be\label{aaa}
\bar a_\b = \frac{4 \bar b_\b}{(1+\bar b_\b)^2} >>  \frac{4 b_\b}{(1+b_\b)^2},  
\ee
then the BC point processes starting respectively from the line $\{t=0\}$ and from the Bernoulli set $\bar A_\b$ locally coincide with a probability
close to $1$. On the one hand, this Bernoulli set is the invariant measure
for the the BC point process with intensity $\bar b_\b$. On the other hand,
this new process
stochastically dominates the original BC (with branching probability $b_\b$) since $\bar b_\b$ must be much larger than $b_\b$
for (\ref{aaa}) to be satisfied.
We can then find an upper bound for (\ref{corr-ff}) using
this new BC point process (with a higher branching rate) and assuming independence between sites.

\medskip

Let us now 
give a more precise argument. 
Let $\bar b>b$, and let us consider 
a sequence of branching rates $\{\bar b_\b\}$,  such that 
$e^\b \bar b_\b\rightarrow \bar b$, as $\beta\uaw\infty$.
Finally we denote by $\bar A_\b$ the Bernoulli random field on $\Z_{even}$
with intensity 
$$
\bar a_\b \ = \ \frac{4 \bar b_\b}{(1+\bar b_\b)^2}. 
$$
In particular, we note that $e^\b \bar a_\b$ converges to $4\bar b$. Furthermore,
$\bar a_\b$ is chosen in such a way that if we denote by $\bar \xi^{\bar A_\b}_\b$
the branching-coalescing point set
with branching probability $\bar b_\b$ starting on the set of
points $\bar A_\b$ at time $0$, then the process 
$(\bar \xi^{\bar A_\b}_\b(n), n\in\N)$ defines a stationary process.

\medskip

For $\bar b_\b > b_\b$, the process $\bar \xi_\b$
stochastically dominates $\xi_\b$, since there exists an obvious coupling
under which  $\xi_\b(n)\subset\bar \xi_\b(n)$, for every $n\in\N$.
As a consequence,
\beqn
\P(B_{\eps,\b,L,T^\pm}) & \leq  &  \P\left(\xi_{\b}^0(T^{-} e^{2\b}) \cap [-L e^{\b}, L e^\b]  \neq \xi_{\b}^{\bar A_\b}( T^{-} e^{2\b})  \cap [-L e^{\b}, L e^\b]   \right) \no  \\ 
&  & + \P( C_{\eps,\b,L,T^{\pm}}  ) \no \\
			      & \leq &
 \P\left(\xi_{\b}^0(T^{-} e^{2\b}) \cap [-L e^{\b}, L e^\b]  \neq \xi_{\b}^{\bar A_\b}(T^{-} e^{2\b})  \cap [-L e^{\b}, L e^\b]   \right) \no \\
& & + \P( \bar C_{\eps,\b,L,T^{\pm}}  ),  \label{supB}	
\eeqn
where 
\beqnn
C_{\eps,\b,L,T^\pm} 
& = & \{\exists \ x\neq x',t\in\N : S_{\b}(x,t) \in \square_{L,T^\pm}, (x,t)\in{\cal M}^{k_\b}_{\b}, \ x,x' \in \ \xi^{\bar A_\b}_{\b}(t), 
\ |x-x'| < e^{\b} \eps, \ \} \\ 
\bar C_{\eps,\b,L,T^{\pm}} 
& = & \{\exists \ x \neq x',t\in\N : S_{\b}(x,t) \in \square_{L,T^\pm}, (x,t)\in{\cal M}_{\b}^{k_\b}, \ x,x' \in \ \bar \xi^{\bar A_\b}_{\b}(t), 
\ |x-x'| < e^{\b} \eps, \ \}, \\ 
\eeqnn
i.e., $\bar C_{\eps,\b,L,T}$ is the same event as $C_{\eps,\b,L,T}$
with a higher branching rate.  
We have
\beqnn
\P(\bar C_{\eps,\b,L,T^\pm}) & \leq & \sum_{(x,t): S_{\b}\left(x,t\right)\in\square_{L,T^\pm}} \ \  \sum_{x'\neq x: |x-x'|<e^\b\eps} 
\ \P\left( x,x' \in \ \bar \xi^{ \bar A_\b}_{\b}(t), \  (x,t)\in{\cal M}^{k_\b}_{\b}   \right) \\
&  = & k_\b  \ \sum_{(x,t): S_{\b}\left(x,t\right)\in\square_{L,T^\pm}}  \ \ \sum_{x'\neq x: |x-x'|<e^\b\eps} 
\ \P\left( x,x' \in \ \bar \xi^{\bar A_\b}_{\b}(t) \right) \\
& =&  k_\b  \ \sum_{(x,t): S_{\b}\left(x,t\right)\in\square_{L,T^\pm}} \ \ \sum_{x'\neq x: |x-x'|<e^\b\eps} 
\ \P\left( x \in \ \bar \xi^{\bar A_\b}_{\b}(t)  \right) \P\left(  x' \in \ \bar \xi^{\bar A_\b}_{\b}(t)  \right), \\
\eeqnn
where the last equality follows from the fact that at any time $t\in\N$,
$\bar \xi^{\bar A_\b}_{\b}(t)$ defines a Bernoulli field on $\Z_{even}$ (resp., $\Z_{odd}$) 
at any even (resp., odd) time. 
Next, by translation invariance,
for every $(x,t)\in\Z_{even}^2$, we have
\beq
\label{tr:inv}
\P\left( x \in \ \bar \xi^{\bar A_\b}_{\b}(t) \right) \ = \ \P\left( x_0(t) \in \ \bar \xi^{\bar A_\b}_{\b}(t)  \right)  \leq  \P\left( x_0(t) \in \ \bar \xi^{0}_{\b}(t) \right),
\eeq
where $x_0(t)=0$ if $t$ is even, and is equal to $1$ otherwise.
Hence, using that (1) for any given point $x$, there are at most $2 \eps e^\b$ 
points 
$x'$ which are at a distance $\eps e^\b$ away, (2) for any given time $t$, 
there are
at most $2 L e^\b$ point in the macroscopic box $\square_{L,T^{\pm}}$ 
and whose time coordinate is given by $t$,
and finally (3) translation invariance (\ref{tr:inv}),
we get that
\beqn
\P(\bar C_{\eps,\b,L,T^{\pm}}) & \leq &  4 \ L\ \eps \  e^{2\b} k_\b  \ \sum_{t \in\Z: T^- e^{2\b} \leq t\leq T^+ e^{2\b}}   \P\left( x_0(t) \in \ \bar \xi^{0}_{\b}(t)  \right) ^2 \no \\
& \leq &
 4 \ L\ \eps \  e^{2\b} k_\b  \ \sum_{t \in\Z:  t\leq T^+ e^{2\b}}    \ 
 \left(\P\left( x_0(t) \in \ \bar \xi^{0}_{\b}(t)  \right) \ \cdot e^{\b}\right)^2 \ \cdot \ e^{-2\b}.
\label{1.03}
\eeqn
Next, recall that 
$k_\b \cdot e^{2\b}$ is assumed to converge
to the
continuous killing parameter $k$ as $\beta\uaw\infty$, 
and by Corollary \ref{c:CLT}, 
provided that $t_\b e^{-2\b} \rightarrow t$ as $\b\uaw\infty$, 
we know that
\beqnn
\limsup_{\b\uaw\infty}   \ \P\left( x_0(t_\b) \in \ \bar \xi^{0}_{\b}(t_\b)  \right)  \cdot e^{\b}
\ = \ 
\fr{e^{-\bar b^2 t}}{\sqrt{\pi t}} \ + \ 2 \bar b \phi(\bar b\sqrt{2 t}).
\eeqnn
Hence, as $\b$ goes to $\infty$, the sum in the RHS of (\ref{1.03})
is a Riemann sum, and we get that
\beqnn
\limsup_{\b\uaw\infty} \  \P(\bar C_{\eps,\b,L,T}) & \leq & 4 \ \eps \  L  \ k \int_0^T \left( \fr{e^{-\bar b^2 t }}{\sqrt{\pi t}} \ + \ 2 \bar b \phi(\bar b\sqrt{2 t}) \right)^2 \ dt.
\eeqnn 
Letting $\eps$ go to $0$ on both sides of the inequality, (\ref{supB}) implies 
that for every $\bar b> b$,
$$
\limsup_{ \eps\daw 0} \limsup_{\b\uaw\infty}  \P(B_{\eps,\b,L,T}) \leq  \limsup_{\b\uaw\infty} 
 \P\left(\xi_{\b}^0(T^{-} e^{2\b}) \cap [-L e^{\b}, L e^\b]  \neq \xi_{\b}^{\bar A_\b}(T^{-} e^{2\b})  \cap [-L e^{\b}, L e^\b] \right)  . 
$$
Letting $\bar b$ tend to $\infty$ on 
the RHS of this inequality and using Lemma 
\ref{sparse:first} then yields our result.

\end{proof}

\bcor\label{xib}
The measure 
$$
m^{L} :=\sum_{x\in\xi^{0}(\theta)\cap[-L,L]} \delta_x
$$
consist of finitely many atoms a.s.. Furthermore,
if $L_\b = [L e^\b]$, then
$$ 
e^{-\b} m_\b^{L} : = \sum_{x\in  \xi_\b^{0}(\theta_\b)\cap [-L_\b,L_\b]} \delta_{e^{-\b}x} \rightarrow m^{L}  \ \ \mbox{in law},
$$
where the convergence is in the vague topology.
\ecor

\begin{proof}
The first  property is a direct consequence of the fact that for
any deterministic $t$, the BC point set
$\xi^0(t)$ is locally finite and that $\theta>0$ a.s..

For the second property, we consider the coupling under which the a.s. convergence statement (\ref{coucou})
holds. First, one can easily show that 
for every $x\in\xi^0(\theta)$, there must exist $x_\b\in\xi_\b^{0}(\theta_\b)$ such
that $x_\b$ converges (after rescaling) to $x$ a.s.. 
(This can be shown along the same lines as (\ref{claim-item}).)
Conversely, one can easily show that any $x_\b$ in $\xi_\b^{0}(\theta_\b)\cap [-L_\b,L_\b]$
must be close to some $x\in\xi^{0}(\theta)\cap[-L,L]$ for $\b$ large enough. 
This together
with the sparseness property proved in Lemma \ref{l:sparse}, directly implies our result.
(Indeed, having more atoms in the discrete measure, as $\beta\uaw\infty$, would imply that at least two atoms 
are at an arbitrary small distance from each other.)

\end{proof}

{\bf Proof of Proposition \ref{kill-kill}.}
First, we show that for any path $\pi\in\Net^b$
such that 
$
\sigma_\pi =0
$
and $age(\pi)>0$,
one can find a sequence $\{\pi_\b\}$ such that, after rescaling, it converges  to $\pi$ (uniformly in $\pi$)
and such that 
$\sigma_{\pi_\b}=0$. 
We can then consider a path $\tilde \pi\in\Net^b$ starting strictly below
$0$  and so that $\pi$ is the continuation of $\tilde \pi$ after time $0$. 
Since $\cUb$ converges (after rescaling) to $\Net^b$, we can then 
consider a sequence ${\tilde \pi}_\b$ converging to $\tilde \pi$ (uniformly) 
and then approximate $\pi$ with the path $\pi_\b$, the path made of the portion 
of $\tilde \pi_\b$ after time $0$ (which also belongs to $\cUb$).

Next, Lemma \ref{JAN} together with the convergence results 
of Lemma \ref{thetab}
and Corollary \ref{xib} implies that 
the sequence 
$\{(\cUb, K_\b,m_\b, \theta_\b)\}$ converges (after proper rescaling)
to $(\Net^b,K,m, \theta)$ in probability. In the following, we consider a subsequence 
$\{\bs\}$
such that 
this convergence is a.s.
and show that the second item of the proposition
is satisfied for the sequence $\{\pi_\bs\}$ -- as defined in the previous paragraph.
First, let us call 
$(x_\bs,\theta_\bs)$ the only killing point at time $\theta_\bs$
approximating the killing
point $(x,\theta)$ -- where the existence and uniqueness of $(x_\bs,\theta_\bs)$  is guaranteed by Lemma \ref{thetab}. Let us now consider  
 any path $\pi$ killed at the point $(x,\theta)$
and let us show that $\pi_\bs$ must also be killed at $(x_\bs,\theta_\bs)$, for $\tilde \beta$ large enough. 
First, we claim that the discrete approximation $\pi_\bs$
must go through the discrete killing point $(x_\bs,\theta_\bs)$. Indeed, the convergence of $m_\bs$ 
to $m$ implies that $m_\bs$ is made of finitely many atoms staying at a positive distance from each other as $\tilde \b\uaw\infty$. Since $\pi_\bs$ must go through one of those atoms, $\pi_\bs$ must go 
through the killing 
point $(x_\bs,\theta_\bs)$ for $\bs$ large enough.  In order
to prove that  $\pi_\bs$ is actually killed at $(x_\bs,\theta_\bs)$, it remains to show that $\pi_\bs$
does not encounter any killing point before time $\theta_\bs$. Since $|\pi|<L/2$ on $[0,\theta]$, 
we must have $ |S_\bs(\pi_\bs)|<2/3 L$ on $[0,\theta_\bs]$ for $\bs$ large enough (again uniformly in $\pi$).
From the definition of $\theta_\bs$,
and since  $\sigma_{\pi_\bs} = 0$,  
the path $\pi_\bs$ can not be killed before $\theta_\bs$, thus proving the claim made earlier.
The converse statement -- i.e., that the killing of a discrete path $\pi_\bs$ implying the killing 
of its continuum analog -- can be proved along the same lines. $ \ \ \ \ \square$

\subsection{Step 3 : Renewal Structure in  the Brownian Net}

For every $L_0,T_0>0$, let us define the sequence
of random variables $\{\theta^{i,L_0,T_0}\}_{i}$ as follows.
\beqn\label{thetail}
\theta^{0,L_0,T_0} & = & -T_0 \no \\ 
\forall i\geq1, \ \theta^{i,L_0,T_0} & = & \inf\{t> \theta^{i-1,L_0,T_0} \ : \ \exists x\in\xi^{\theta^{i-1,L_0,T_0}}(t)\cap[-L_0,L_0] \ \ \mbox{s.t.} \  (x,t)\in\mathcal{M}^k \}. \no
\label{sncf}
\eeqn
In this section, we will prove the following result.

\bprop\label{delta-theta}
$$
\max_{i : \theta^{i+1,L_0,T_0} \leq T_0} ( \theta^{i+1,L_0,T_0}-\theta^{i,L_0,T_0} ) \rightarrow 0  \ \mbox{in law}
$$
as $L_0$ goes to $\infty$.
\eprop

\bigskip

We proceed in several consecutive steps starting with a general result on arrays of random variables.

\blem\label{1-2-3}
Let $T>0$ and let $\{X_{k,n}\}$ be an array of random variables 
such that 
\begin{enumerate}
\ii For every $n$, $\{X_{k,n}\}_k$ is a sequence of positive
i.i.d. random variables.
\ii $\limsup_{n\uaw\infty}\E(X_{1,n})=0$, 
\ii For every $\gamma>0$, there exists $\{\alpha_n\equiv\alpha_n(\gamma)\}$
with $\alpha_n\uaw\infty$ and such that
the two following conditions hold: 
(1) $\limsup_{n\uaw\infty} \alpha_n \E(X_{1,n}) = \infty$,
(2) $\limsup_{n\uaw\infty} \alpha_n \P(X_{1,n} \geq \gamma) = 0$.  
\end{enumerate}
For $T\in(0,\infty)$, let
$$
N_n \ = \ \sup\{k \ : \ X_{1,n}+\cdots+X_{k,n} < T  \}.
$$
Then
$$
\max_{0\leq i\leq N_n} X_{i,n} \rightarrow 0 \ \ \mbox{in law as $n\rightarrow\infty$}.
$$
\elem

\begin{proof}
We start by giving some intuition behind the conditions of the  fourth item. Intuitively, the first condition on $X_{1,n}$
ensures that $\E(X_{1,n})$ is not too small and therefore
that 
$N_n$ does not grow too fast as $n\uaw\infty$.
Then, the second condition on the tail of $X_{1,n}$ ensures
that for a given value of $N_n$,
the maximum of $\{X_{1,n},\cdots,X_{N_n,n}\}$
is also not too large. More precisely,
for every $\gamma$,
\beqn
\P(\max_{i\leq N_n} X_{i,n} \geq \gamma)
& \leq & \P(N_n\geq\alpha_n) + \P(\max_{i\leq \alpha_n} X_{i,n} \geq \gamma) \no \\
&  = & \P(N_n\geq\alpha_n) +  1 - \left(\P(X_{1,n} < \gamma \right) ^{\alpha_n} \no \\
& =  &   \P(N_n\geq\alpha_n) + 1- \left(1- \frac{1}{\alpha_n} \cdot [\alpha_n\P(X_{1,n}\geq \gamma)]\right) ^{\alpha_n}  \label{small-t},
\eeqn
where we used the independence to deduce the first equality. Since
$\alpha_n \P(X_{1,n} \geq \gamma) \rightarrow 0$, it remains 
to show that the first term on the RHS of (\ref{small-t}) vanishes as $n\uaw\infty$. 
As we shall see, this will be proved by using $\alpha_n \E(X_{1,n}) \rightarrow \infty$.

By Wald's Theorem, we must have
\beqnn
\E(N_n) \ \E(X_{1,n}) = \E(X_{1,n}+\cdots+X_{N_n,n}). 
\eeqnn
Noting that $\E(X_{1,n}+\cdots+X_{N_n,n}) \leq T$, we have
\beqn
\E(N_n) \ \E(X_{1,n}) & \leq & T. \label{waldo}
\eeqn
By Chebytchev, this implies that 
\beqnn
\P(N_n\geq \alpha_n) & \leq & T/[ \E(X_{1,n}) \alpha_n]
\eeqnn 
and the conclusion follows from the fact that
$
\alpha_n \E(X_{1,n}) \rightarrow \infty
$.
\end{proof}

\blem\label{t-finite}
If $\theta^1$ is defined as in (\ref{thetaL}) with $L=1$,
then for every $M\in\N$, the $M^{th}$ moment 
of $\theta^1$ is finite.
\elem
\begin{proof}
Let us first define a path $\pi\in\Net^b$
starting at time $0$. In order to do so,
we consider the following sequence 
of  stopping times.
\beqnn
T_0 = 0, \ \ z_0 = (+1,0), \\
T_{2n+1} \ = \ \inf\{t>T_{2n} \ : \ l_{z_{2n}} =-1\}, \ \ \ z_{2n+1} \ = \  (-1,T_{2n+1}), \\ 
T_{2n+2} \ = \ \inf\{t>T_{2n+1} \ : \ r_{z_{2n+1}} =+1 \}, \ \ \ z_{2n+2} \ = \  (+1,T_{2n+2}), \\
\eeqnn
where $l_{z_{2n}}$ (resp., $r_{z_{2n+1}}$) is the leftmost (resp., rightmost) path starting
from $z_{2n}$ (resp., $z_{2n+1}$). 
Finally, we define $\pi$ 
to be the path starting starting from
$(+1,0)$ and such that
\beqnn
\pi := \left\{ \begin{array}{cc}
l_{z_{2n}} & \ \ \ \mbox{on $[T_{2n},T_{2n+1}]$}, \\
r_{z_{2n+1}} & \ \ \ \mbox{on $[T_{2n+1},T_{2n+2}]$}.
\end{array} \right.
\eeqnn
Since the Brownian net is stable under hopping, 
the path $\pi$ belongs to $\Net^b$. As a consequence,
\beqn
\P(\theta^1\geq\gamma)
& \leq & \E[ \exp(-k \int_{[0,\gamma]} 1_{\pi(t)\in[-1,1]}) dt ].				\label{noel}	
\eeqn
Thus it remains
to show that the RHS 
of the inequality is integrable.
Let us denote by $\Delta t_n$ 
the amount of time spent
by $\pi$ in the space interval $[-1,1]$, during the time interval $[T_{n},T_{n+1}]$.
$$
 \int_{[0,t]} ds 1_{\pi(s)\in[-1,1]}/t  \ \geq \ \sum_{n\leq N(t)} \Delta t_n/t,  
$$
where 
$$
N(t) \ = \ \sup\{n \ : \ T_n(t) \leq t  \}.
$$
By the strong Markov Property, the 
sequence $\{(\Delta t_n, \Delta T_n:=T_{n+1}-T_{n})\}$
consists of i.i.d.
random variables. 
Furthermore, $\Delta T_1$ is distributed as the first $0$-hitting time 
of a drifted Brownian motion with drift $-b$
starting at $1$. Such a random variable is 
known to have a finite first moment.  
As an easy consequence of the Strong Law of Large Numbers,
we get that 
\beqnn
\liminf_{t\uaw\infty} \frac{1}{t} \int_{[0,t]} 1_{\pi(s)\in[-1,1]} ds \  
&  \geq & \lim_{t\uaw\infty}  \ \sum_{n\leq N(t)} \Delta t_n/ N(t)  \cdot \ N(t)/t  \\
&   = &  m:=\E(\Delta t_1) / \E(\Delta T_1). 
\eeqnn
Not only have $\Delta T_1$ and $\Delta t_1$ a finite first moment,
but they both have exponential moments. With a little bit of extra work,
one can then easily establish the 
following large deviation estimate 
for the RHS of the previous inequality : 
\beqnn
\ \forall c <  m,  \ \ \exists a,A\in(0,\infty), \  \ \forall t>0  \   \  \   \P\left( \sum_{n\leq N(t)} \Delta t_n    \leq c t \right) \leq A \exp(-a t).
\eeqnn
This and (\ref{noel}) easily imply our lemma.
\end{proof}

{\bf Proof of Proposition \ref{delta-theta}.}
By the strong Markov Property and the translation invariance of
the Brownian net, the sequence 
$\{\theta^{i+1,L_0,T_0} - \theta^{i,L_0,T_0} \}_{i}$
is identical in distribution to
independent copies of $\theta^{L_0}$.
By Lemma \ref{1-2-3}, it is then enough to show (item 2)  
$
\lim_{L_0\uaw\infty} \E(\theta^{L_0}) = 0
$
, and (item 3) for every $\gamma>0$, there exists $\alpha_{L_0}\uaw\infty$
such that 
\be\label{co--}
\lim_{L_0\uaw\infty} \alpha_{L_0} \P(\theta^{L_0}\geq \gamma) = 0, \ \ \  \ \lim_{L_0\uaw\infty} \alpha_{L_0}  \E(\theta^{L_0})  = \infty.
\ee

\bigskip

{\it (item 3, part 1.)} Let us start by showing  that
for every $\gamma>0$, there exists $a(\gamma)\in(0,1)$
such that 
\be\label{caki}
\limsup_{L_0\uaw\infty} \
 \alpha_{L_0}  \ \P( \theta^{L_0}\geq \gamma)= 0 \ \  \mbox{with $\alpha_{L_0} :=(1/a(\gamma))^{L_0}$ }.
\ee
In order to see that, we subdivide the rectangle $[-[L_0],[L_0]]\times[0,\gamma]$
(where $[a]$ denotes the integer part of a)
into sub-rectangles of the form $[n,n+1]\times[0,\gamma]$, with $n\in\N$. In 
each of those rectangles, we consider the subset of points $z$ which can be 
connected to a point $z'$ in
the interval $[n,n+1]\times\{0\}$ by a path
in the net remaining in that rectangle.
Let us denote 
by $D_n(\gamma)$ the event that there is no killing
mark on this set of points. The
indicator random variables of the events 
$\{D_{j}(\gamma)\}$ are i.i.d. positive random variables,
with $\P( D_j(\gamma) )\in(0,1)$.
As a consequence,
$$
\P(\theta^{L_0}\geq \gamma)  = \P(\not\exists \ \ \mbox{any mark
on a path starting at time $0$ in $[-L_0,L_0]\times[0,\gamma]$}) 
\leq \P(D_1(\gamma))^{2 L_0}.
$$
Taking any $a(\gamma)\in(\P(D_1(\gamma)),1)$ yields (\ref{caki}).

\bigskip

({\it item 2.}) This easily implies the second item, i.e., that $\lim_{L_0\uaw\infty} \E(\theta^{L_0})=0$. In 
order to see that, we first note that 
$\theta^{1} \geq \theta^{L_0}$, for $L_0\geq1$,
implying that
$$
\P(\theta^{1} \geq \gamma) \geq \P(\theta^{L_0}\geq \gamma).
$$
By Lemma \ref{t-finite},
$$
\E(\theta^1) =  \int_{0}^\infty \P(\theta^{1} \geq \gamma) \ d\gamma <\infty
$$
and the Dominated Convergence Theorem
implies that
\beqnn
\lim_{L_0\uaw\infty} \E(\theta^{L_0}) = \int_0^\infty \lim_{L_0\uaw\infty}  \P(\theta^{L_0} \geq \gamma) d\gamma = 0.
\eeqnn

\bigskip

{(\it item 3, part 2.)} It remains to show (\ref{co--}) for $\alpha_{L_0}=(1/a(\gamma))^{L_0}$. 
We will actually prove a stronger result, namely that
$$
\liminf_{L_0\uaw\infty} \E(\theta^{L_0})  \ L_0^2 > 0.
$$

Let us define $\cL_{L_0}(t)$
to be the time length measure
accumulated by the set of paths starting from time $0$
in the rectangle $[-L_0,L_0]\times[0,t]$.
For a given $\gamma>0$,
we start by constructing a sequence of positive times
$\{t(L_0)\}_{L_0}$ such that 
$$
\forall L_0>0, \ \P(\cL_{L_0}(t(L_0)) \leq \gamma) \geq \half.
$$
First, it is easy to see that there exists $c>0$ (independent of $L_0$)
such that 
$$
\forall t\in[0,1], \ \ \E(\cL_{L_0}(t)) = 2 L_0 \int_0^t \left( \fr{e^{-b^2 s }}{\sqrt{\pi s}} \ + \ 2 b \phi(b\sqrt{2 s}) \right) \ ds \leq c \sqrt{t} L_0.
$$
Chebytchev's inequality then implies
that for such $t$'s 
\beqnn
\P(\cL_{L_0}(t) \geq \gamma) &  \leq & \frac{1}{\gamma} \E(\cL_{L_0}(t))
\leq   \frac{c \sqrt{t} L_0}{\gamma}.    
\eeqnn
Hence,
\beqnn
\P(\cL_{L_0}(t) \leq \gamma) \geq (1 - \frac{c \sqrt{t} L_0}{\gamma}). 
\eeqnn
Taking $t(L_0)=(\frac{\gamma}{2cL_0})^2$,
we get (as claimed earlier) that for $L_0$ large enough (so that $t(L_0)<1$),
\beqn\label{ffd}
\P(\cL_{L_0}(t(L_0)) \leq \gamma) \geq \half.
\eeqn
On the other hand, for every $t\geq0$,
\beqnn
\E(\theta^{L_0}) & \geq & \E(\theta^{L_0} \ | \ \cL_{L_0}(t) \leq \gamma ) \cdot \P( \cL_{L_0}(t) \leq \gamma) \\
& \geq & t \P(\mbox{ $\not\exists$ any mark on paths starting from $\R\times\{0\}$ in $[-L_0,L_0]\times[0,t]$} \ | \ \cL_{L_0}(t) \leq \gamma) 
\cdot  \P( \cL_{L_0}(t) \leq \gamma)\\
& \geq & t \exp(-k\gamma)   \P( \cL_{L_0}(t) \leq \gamma).
\eeqnn
Taking $t\equiv t(L_0)$ as in (\ref{ffd}), we get that 
$$
\lim_{L_0\uaw\infty} \E(\theta^{L_0}) L_0^2 \geq \frac{\gamma^2}{8c^2} \exp(-k\gamma)>0.
$$

\
\subsection{Step 4 : Approximation in a Finite Box}
\label{approx}
We now have all the ingredients to prove the invariance principle of Theorem \ref{onewebteo}.
First, by definition
of our topology (see Section \ref{the-space}), for every arbitrarily small $\gamma>0$, it is enough to construct a coupling between the discrete and continuum
level
such that 
\beqn
\P\left( \limsup_{\b\uaw\infty} \ \ \max_{\Net^{b,k}} \min_{S_\b(\cUbk)} d(\pi,\pi_\b) \leq \gamma \right) \geq 1-\gamma \label{one-half}, \\    
\P\left(\limsup_{\b\uaw\infty} \ \  \max_{S_\b(\cUbk)} \min_{\Net^{b,k}} d(\pi,\pi_\b)  \leq \gamma  \right) \geq 1-\gamma\label{no-proof}. 
\eeqn
In the following, we will construct a coupling such that (\ref{one-half}) is satisfied. We leave the reader to
convince herself that (\ref{no-proof}) is also satisfied (under the same coupling) by reasoning along the same lines.
In this section, we will focus on a weaker form of (\ref{one-half}), restricting our attention to a finite box and adding some restrictions 
on the age of a path. 
As we shall see in the next (and final)
section, this will
easily imply the desired result.

\bigskip

In order to describe more precisely the main result of this section, we start with some definitions. 
In the following, ${\cal A}^T$ will denote the set of paths 
in $\Net^{b,k}$ starting and being killed in the time interval $[-T,T]$, i.e., $[\sigma_\pi,e_\pi]\subseteq[-T,T]$. 
Next, the construction of the Brownian net with killing $\Net^{b,k}$ from  the net
$\Net^b$ -- see Section 
\ref{netconstruction} -- provides a natural coupling $(\Net^b,\Net^{b,k})$. Under this coupling,
for any path 
$\pi\in\Net^{b,k}$, one can find a pair $(\tilde \pi, e_{\pi})$,
where $\tilde \pi$ is a
path of the net $\Net^b$
starting from $\sigma_\pi$ and coinciding
with $\pi$ on $[\sigma_\pi,e_\pi]$. Let us now consider
the set of paths in $\Net^{b,k}$
with  
$
age(\tilde \pi) \geq  \eps.
$ 
In particular,
for any such path,
there must exist  another path in the Brownian net $\Net^b$ -- call it $\tilde \pi^\eps$ --
such that $\sigma_{\tilde \pi^\eps}=\sigma_\pi- \eps$
and coinciding with $\tilde \pi$ after time $\sigma_\pi$.
In the following, we will restrict our attention to
the  subset of ${\cal A}^T$
such that  $| \tilde \pi^\eps(t)| < L $ for $t\in[-T \vee (\sigma_\pi-\eps),e_\pi]$. 
We will denote this set by ${\cal A}^{T,L,\eps}$
and we will show that any
path in this set can be uniformly approximated by
a discrete path:

\bprop\label{finite-box} Under the coupling (\ref{coucou}),
$$
\lim_{\b\uaw\infty} \ \ \ \max_{ {\cal A}^{T,L,\eps}} \min_{ S_\b({\cUbk})}    \ d(\pi,\pi_\b) \ = \ 0 \ \ \mbox{in law.}
$$
\eprop

\bigskip

Let $L_0$ be such that $L_0>L$.
Recall the definition 
of the sequence $\{\theta^{i,L_0,T}\}_i$ given in (\ref{sncf}).
Since $T$ will remain fixed throughout this proof, 
we will drop the $T$ superscript to ease the notation and write
$
\theta^{i,L_0} \equiv \theta^{i,L_0,T}.
$

\begin{figure}
\centering
\includegraphics[scale=.5]{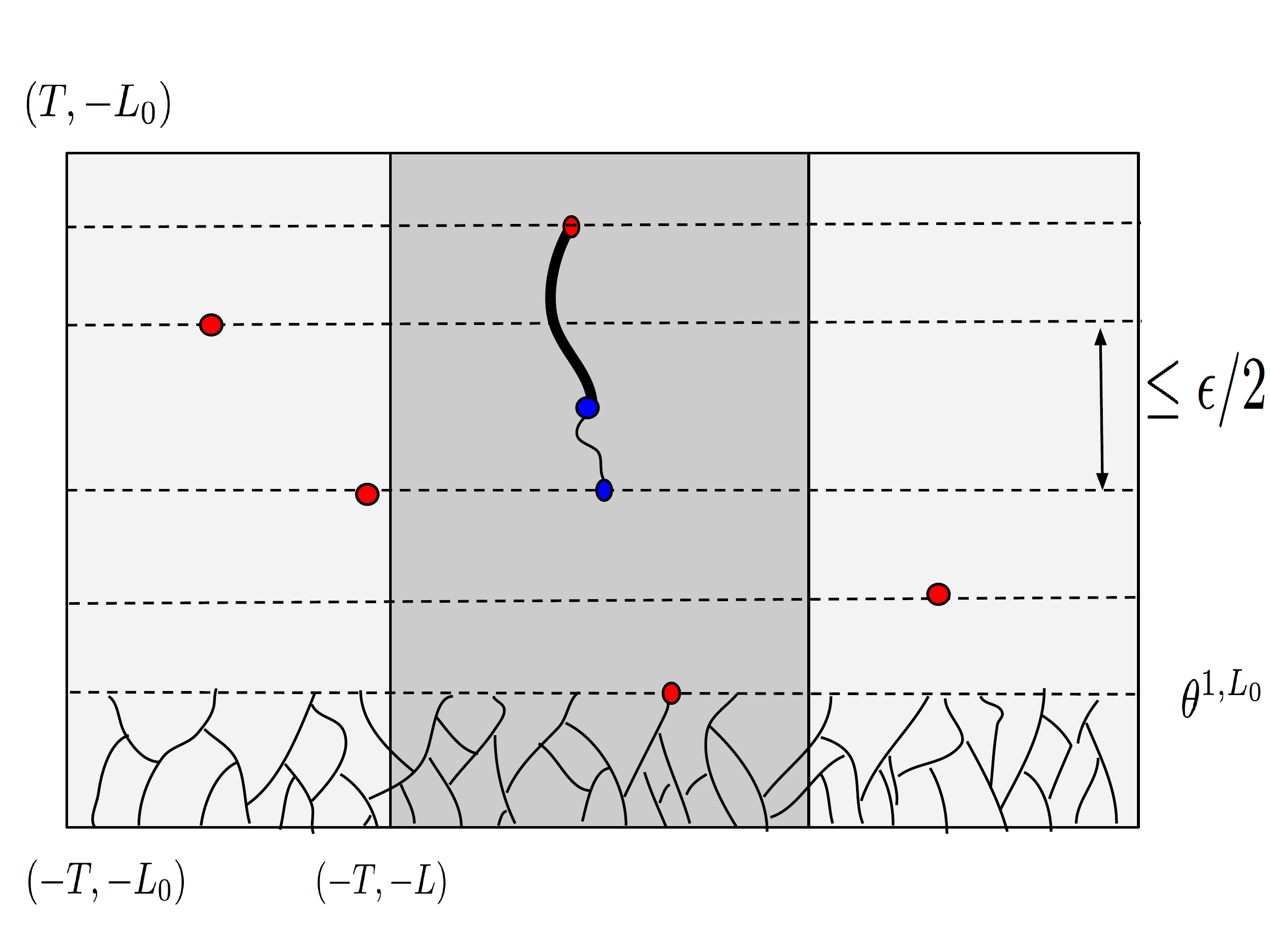}
\caption{Schematic diagram of the proof of Lemma \ref{2eps}, when 
$\{\eps \geq 2 \max_{i : \theta^{i+1,L_0}<T} (\theta^{i+1,L_0}-\theta^{i,L_0}) \}$. Red dots are killing points on $\R\times\{\theta^{i,L_0} \}$ with $x$-ccordinates in 
$\xi^{\theta^{i,L_0}}(\theta^{i-1,L_0})$.
Blue points are the starting points of the (wide) path $\pi$ (with age $\geq\eps$) and its extension $p$ (narrow path and its continuation $\pi$). }
\label{dessin}
\end{figure}

\blem\label{2eps}
On the event $\{\eps \geq 2 \max_{i : \theta^{i+1,L_0}<T} (\theta^{i+1,L_0}-\theta^{i,L_0}) \}$,
any path in ${\cal A}^{T,L,\eps}$ must be killed at some $\theta^{j,L_0}\in[\sigma_\pi,T]$.
\elem

\begin{proof}

Let us consider 
a path $\pi\in{\cal A}^{T,L,\eps}$ and
and let us denote by $i$ the index such that 
$\sigma_\pi\in[\theta^{i-1,L_0},\theta^{i,L_0})$. 
Recall that $\pi$
is characterized by a pair $(\tilde \pi ,e_\pi)$
where $\tilde \pi\in\Net^b$ and $e_\pi$
is the killing time of $\pi$.
Since
the age of $\tilde \pi$ must be greater than
$\eps$, on our event,
there must exist $p\in\Net^b$
starting at time $\theta^{i-1,L_0}$ and such 
that $p$ and $\tilde \pi$ coincide for $t\geq\sigma_{\tilde \pi}$(see Fig. \ref{dessin}). 
(We note in passing that $age(p)>0$.)
Furthermore, by definition of ${\cal A}^{T,L,\eps}$,
the path $p$ (and thus $\tilde \pi$) 
must remain in the space interval $[-L,L]$
during the time interval $[\theta^{i-1,L_0},\theta^{i,L_0}]$.
By definition
of $\theta^{i,L_0}$, it follows that 
$p$ (and thus $\tilde \pi$) must survive on $[\sigma_\pi,\theta^{i,L_0})$.
Using a simple induction on $j\geq i$,
one can then generalize the previous reasoning
to show that if $\pi$ survives up to and including  $\theta^{j,L_0}$,
it must survive up to but not including $\theta^{j+1,L_0}$. This
implies the claim made earlier, i.e.,
that if $\pi$ is killed inside $(-T,T)$, $\pi$ must be killed at some $\theta^{j,L_0}$. 

\end{proof}

Next, let us define  $\theta^{i,L_0}_\b$ analogously to $\theta^{i,L_0}$, i.e., 
\beqnn
\theta^{0,L_0}_\beta & = & - [T e^{2\b}] , \\
\forall i\geq1, \ \ \theta^{i,L_0}_\b & = & \inf\{t> \theta^{i-1,L_0}_\b \ : \ \exists x\in\xi^{\theta_\b^{i-1,L_0}}_\b(t)\cap[-L_0 e^\b,L_0 e^\b] \ \mbox{s.t.} \  (x,t)
\  \mbox{is a killing point} \}.
\eeqnn
First, by the strong Markov property
and translation invariance $\{\theta^{i+1,L_0}-\theta^{i,L_0}\}$ (resp., $\{\theta^{i+1,L_0}_\b-\theta^{i,L_0}_\b\}$) 
is a sequence of i.i.d. random variables distributed like $\theta^{L_0}$ (resp, $\theta_\b^{L_0}$) in (\ref{thetaL}). 
A straightforward extension of Proposition \ref{kill-kill} then directly yields 
the following result.
\blem
\label{kill-kill2} 
Let us consider 
the set of paths $\pi\in\Net^b$ (net with no killing) such that there exists $i\in\N$
so that 
$\sigma_\pi=\theta^{i,L_0}\in[-T,T], age(\pi)>0$ and  $\max_{t\in[-T,T]} |\pi(t)|<L_0/2$.
Under our coupling (\ref{coucou}), there 
is a subsequence $\{\bs\}$ such that the following 
properties hold a.s..
\begin{enumerate}
\item For every $i$, $e^{-2\bs} \theta_\bs^{i,L_0}$ converges to $\theta^{i,L_0}$ a.s..
\item For every $\pi$ in the set,
there exists a discrete sequence $\{\pi_\bs\}\in{\cal U}^{\tilde \b}$ 
such that (i)
$
d(\pi,S_\bs(\pi_\bs)) \rightarrow 0
$
uniformly in $\pi$, and (ii)
for $\bs$ large enough (uniformly in $\pi$),
$$
\pi   \ \mbox{is killed at some $\theta^{j,L_0}<T$} \  \Longleftrightarrow  \ \pi_\bs \  \mbox{is killed at $\theta_\bs^{j,L_0}$}. 
$$ 
\end{enumerate}
\elem

Finally, 
we will need the following easy lemma.

\blem\label{the-end-h}
Let $\{X_\b\}$ be a sequence of random variables valued in $\R$. $X_\b$ converges to $0$ in probability
if and only for every subsequence $\{\bs\}$ and for every $\gamma>0$, one can extract a further subsequence $\bss$ 
such that $\P( \limsup_{\bss\uaw\infty} |X_\bss| = 0  )\geq1-\gamma$.
\elem

\begin{proof}
It is well known that $X_\b$ converges to $0$ in probability if and only if 
every subsequence has a further subsequence converging to $0$ a.s.. The lemma 
then follows easily.
\end{proof}

{\bf Proof of Proposition \ref{finite-box}.}
In the following, we take $L_0>2 L$
and
we consider $\{\bs\}$ to be an arbitrary subsequence of $\{\b\}$.
By Lemma \ref{kill-kill2},
we can extract a further subsequence $\{\bss\}$ -- depending on $L_0$ --
such that 
the relations between 
the discrete and continuum levels listed in this lemma hold.
We will now show that
\be\label{end22}
\P(\limsup_{\bss\uaw\infty} \  \max_{{\cal A}^{T,L,\eps}} \min_{ S_\bss({\cal U}^{b_{\bss},k_{\bss}})}   d(\pi,\pi_\bss) \ = \ 0 )  \  \geq \ 1- \P(\{\eps \geq 2 \max_{i : \theta^{i+1,L_0} \leq T} (\theta^{i+1,L_0} -\theta^{i,L_0})\} ).
\ee
Since $L_0$ is arbitrarily large, our proposition follows by applying Proposition \ref{delta-theta} and
Lemma \ref{the-end-h}.

In order to show (\ref{end22}), it is sufficient to prove that for every 
realization of the event
$$
\{\eps \geq 2 \max_{i : \theta^{i+1,L_0} \leq T} (\theta^{i+1,L_0} -\theta^{i,L_0})\},
$$
we must have
$$
\limsup_{\bss\uaw\infty}  \max_{{\cal A}^{T,L,\eps}} \min_{ S_\bss({\cal U}^{b_{\bss},k_{\bss}})}   d(\pi,\pi_\bss) \ = \ 0.
$$
Let $\pi \in {\cal A}^{T,L,\eps}$ and let $i$ be such that $\sigma_\pi\in[\theta^{i,L_0},\theta^{i+1,L_0})$. 
Reasoning along the same lines as  Lemma \ref{2eps}, conditional on our event, we can find a path $\pi'\in\Net^{b,k}$ such that 
(1) $\pi'$ starts at time $\theta^{i,L_0}$ and $age(\pi')>0$, (2) $\pi$ and $\pi'$
coincide after time $\sigma_\pi$, and (3) $|\pi'|<L$ on $[\theta^{i,L_0},\theta^{i+1,L_0})$ -- see again Fig.~\ref{dessin}. Since $L_0>2L$,
the path $\pi'$ belongs to the set 
of paths (at the continuum level) 
defined in Lemma \ref{kill-kill2}. By definition of $\{\bss\}$, $\pi'$
can be approximated (after rescaling) by a sequence $\{\pi_\bss'\}$ with the properties listed in the second item of this lemma. 
In order to approximate $\pi$, we truncate 
the path $\pi_\bss'$,
only considering
the section of $\pi_\bss'$
after some time $t_\bss$
such that $e^{-2\bss} t_\bss\raw \sigma_\pi$.
This defines a new sequence of paths $\{\pi_\bss\}$
approximating the path $\pi$ and such  
that 
$$
\mbox{$\pi$ is killed at  $\theta^{j,L_0}<T$} 
\Longleftrightarrow 
\mbox{$\pi_\b$ is killed at the discrete analog $\theta^{j,L_0}_\bss$.}
$$

Since any path in
 ${\cal A}^{T,L,\eps}$ must be killed at some $\theta^{j,L_0}\in[\sigma_\pi,T]$ 
 (by Lemma \ref{2eps}), 
 and since
 $
\{e^{-2\bss} \theta^{i,L_0}_\bss\}_i \rightarrow \theta^{i,L_0}
$ a.s.,
(\ref{end22}) follows. The ends the proof of the proposition.

\subsection{Step 5 : Conclusion}
\label{convergence}

In this section,
we prove (\ref{one-half}), i.e.,
that any continuous path in the killed Brownian net
can be uniformly approximated by a discrete killed path. In the following,
we work under the coupling discussed in Lemma \ref{kill-kill2}.
We decompose the set $\Net^{b,k}$ into three sets: 
(1) the set of paths killed before $-T$ or after $+T$, 
(2) the set of paths starting and killed on $[-T,T]$-- denoted by ${\cal A}^T$ in the previous section -- and (3) 
the set of paths killed on $[-T,T]$ but starting before time $-T$. In the following, we will
show that when $T$ goes to $\infty$, each path in any of those 
subsets is well approximated by a well-chosen discrete path.

\bigskip

{\it Part 1.} We first consider the set of paths in $\Net^{b,k}$
killed before time $-T$ or after $+T$. 

First, let us take a path $\pi$ killed after time $T$. Any such path 
can be represented by a pair $(\tilde \pi,e_{\tilde \pi})$, with $\tilde \pi$
in the Brownian net (with no killing) and $e_{\pi}$
its killing time. By convergence 
of the discrete net to its continuum analog, we can find $\tilde \pi_\b$
in the discrete net approximating $\tilde \pi$ (after rescaling). The definition of our topology
(see Section \ref{the-space})
immediately implies that $\tilde \pi_\b$ is also a good approximation of the killed path $\pi$, provided that $T$ is large enough -- intuitively,
any  killing point  above $\R\times \{T\}$ is close to $\infty$ in our compactification of $\R^2$.

Analogously, since any point below the horizontal line
$\R\times\{-T\}$
is  close to $-\infty$,
a path killed before time $-T$ is well approximated by any (properly rescaled) discrete path killed before time $-T$.

\bigskip

{\it Part 2.} We now approximate the set of paths starting and killed in the interval $[-T,T]$.
As usual, $\square_{L,T}$
will refer to the box $[-L,L]\times[-T,T]$. We will also make use of the following notation.
First, we write ${\cal A}^{T,L}$
for the set of the paths 
killed during
the time interval $[-T,T]$,
and staying in the space interval $[-L,L]$
during this time period -- with no restriction on the age of the path. 
Note that ${\cal A}^{T,L}\equiv {\cal A}^{T,L,0}$, where  ${\cal A}^{T,L,\eps}$
was defined in the previous section. Next, define
\beqnn
{\cal B}^{T,L} & :  & \mbox{ paths $\pi\in{\cal A}^T$ with $|\pi(\sigma_{\pi})|<L$}, \\
\bar {\cal B}^{T,L} & :  & \mbox{ paths $\pi\in{\cal A}^T$ with $|\pi(\sigma_{\pi})|\geq L$}, \\
\bar {\cal A}^{T,L} & : &  \mbox{ paths $\pi\in{\cal A}^T$ with $|\pi|\geq L$ on $\square_{L,T}$.} 
\eeqnn
Finally, the same notation, but with an extra $\b$ subscript,
will refer to the analogous quantities 
for the {\it rescaled} set $S_\b(\cUbk)$. Next,
\beqn
\P\left( \max_{{\cal A}^T} \min_{S_\b(\cUbk) } d(\pi,\pi_\b)  \ > \ \gamma \right)
& \leq &  
\P\left( \max_{{\cal B}^{T,3L}} \min_{S_\b(\cUbk) } d(\pi,\pi_\b)   \ > \ \frac{\gamma}{2} \right),
+
\P\left( \max_{\bar {\cal B}^{T,3L}} \min_{S_\b(\cUbk) } d(\pi,\pi_\b)      \ > \ \frac{\gamma}{2} \right), \no \\
& \leq & 
\P\left( \max_{{\cal A}^{T,4L}} \min_{S_\b(\cUbk) } d(\pi,\pi_\b)   \ > \ \frac{\gamma}{2} \right)
+
\P\left( \max_{\bar {\cal A}^{T,2L}} \min_{S_\b(\cUbk) } d(\pi,\pi_\b)      \ > \ \frac{\gamma}{2} \right) 
\no \\ 
 & & +   2 \left(1 - \P(  {\cal B}^{T,3L} \subset  {\cal A}^{T,4L}, \ \ \bar{\cal B}^{T,3L} \subset  \bar {\cal A}^{T,2L} ) \right). 
\label{nounou}
\eeqn
Part 2 follows 
by taking successively $L$ and {\it then} $T$ to $\infty$
and using the three following lemmas.

\blem\label{difficult-0}
$$
\forall L,T, \ \  \ \lim_{\b\uaw\infty}  \max_{{\cal A}^{T,L}} \min_{S_\b(\cUbk)} d(\pi,\pi_\b)   = 0 \ \mbox{in law.} 
$$
\elem
\begin{proof}
In the previous section, 
we showed this result for ${\cal A}^{T,L,\eps}$ with $\eps>0$. 
We now would like to extend this result to the case $\eps=0$.

\bigskip

{\it Case 1. $\pi$ is killed after time $\sigma_{\pi}+\eps$}. Let us consider $\bar \pi$, the section 
of the path after $\sigma_\pi+\eps$. In particular, the age 
of $\bar \pi$ is equal to $\eps$ and
by Proposition \ref{finite-box}, we can approximate the path $\bar \pi$ by a discrete path $\bar \pi_\b\in S_\b(\cUbk)$.
The two main contributions to the distance between $\bar \pi_\b$ and $\pi$ 
are (1) the difference in starting time -- equal approximatively to $\eps$ and (2) the variation of $\pi$ on $[\sigma_\pi, \sigma_\pi+\eps]$.
As $\eps\daw0$, the first term vanishes.
Since the compactness of $\Net^{b,k}$ implies equicontinuity of the net paths 
in
the rectangle
$[-L,L]\times[-T,T]$ (by Arzela-Ascoli), the second term also vanishes (uniformly in $\pi$). 

\bigskip

{\it Case 2. $\pi$ is killed in $[\sigma_\pi,\sigma_{\pi}+\eps]$}. 
Let us consider the quantity
$$
\max_{ z \in\square_{L,T} } \max_{z'\in S_\b({\cal M}^{k_\b}_{\b})} \ \   \rho(z,z'), 
$$
the maximum distance between a point in the box
$\square_{L,T}$ 
from the set of rescaled discrete killing points. By a standard Borel-Cantelli argument,
this quantity a.s. goes to $0$. Thus, we can always find  a
sequence of discrete killing
points $\{z_\b\}$ converging (after rescaling)
to the starting point of $\pi$
as $\b$ goes to $\infty$. Considering the path
consisting of the single point
$z_\b$ (starting and ending at the same point),
this path is a good approximation
of $\pi$ for $\eps$ small enough.

\end{proof}

\blem\label{easy-0}
$$
\forall T>0, \  \ \lim_{L\uaw\infty}\lim_{\b\uaw\infty} 
\ \  \  \max_{\bar {\cal A}^{T,L}} \min_{S_\b(\cUbk)} d(\pi,\pi_\b) \ = \ 0,  \ \ \ \ \mbox{in law}. 
$$
\elem
\begin{proof}

Here we consider the paths $\pi\in\bar {\cal A}^{T,L}$ with $\sigma_\pi\geq -T$.
The case $\sigma_\pi<-T$ can easily be deduced from
the previous case by approximating the section of the path
after time $-T$. 

Let $\gamma>0$ be an arbitrarily small number
such that $T/\gamma\in\N$. We claim that 
for $\b$ large enough, 
for every $i,j\in\Z$ such that $i\leq j$ and $|i\gamma|,|j\gamma|\leq T$,
there exists a.s. a path $\pi_\beta$ in $\bar {\cal A}^{T,L}_\b$ with $\pi_\b>L$ 
starting from 
the interval $[i\gamma,(i+1)\gamma)$
and killed in the interval $[j\gamma,(j+1)\gamma)$.
Before proving this, we first show how this implies
our lemma.
By definition 
of the distance between paths,
for every $\pi,\pi_\b$  such that $\pi,\pi_\b>L$ for all $t\in[-T,T]$,
we must have
\beqnn
d(\pi,\pi_\b) 
& \leq & |\psi(\sigma_\pi) - \psi(\sigma_{\pi_\b}) | + |\psi(e_\pi) - \psi(e_{\pi_\b})| + (1-\tanh(L)). 
\eeqnn
where $\psi=\tanh(t)$, as in Section \ref{the-space}. 
(Note that the third term corresponds to the fact that two points with the same 
time coordinate and space coordinate $>L$ are close to each other in our
topology when $L$ is large.)
For any path $\pi\in \bar {\cal A}^{T,L}$,
let $i,j$ be such that 
$
\sigma_\pi\in [i\gamma,(i+1)\gamma)
$
and
$
e_\pi\in [j\gamma,(j+1)\gamma).
$
Finally,
let us consider any path $\pi_\beta>L$ in 
$\bar {\cal A}^{T,L}_\b$
with its starting and ending times 
belonging to the
same time intervals (which exists by the claim made earlier).
By the previous inequality, we must have
\beqnn
d(\pi,\pi_\b) 
 & \leq & 2 \gamma + (1-\tanh(L)).
\eeqnn
By symmetry, the same discrete approximation can be obtained 
for a path $\pi \in  \bar {\cal A}^{T,L}$ with $\pi<-L$.
Hence,  
we get that for $\b$ large enough
$$
 \max_{\bar {\cal A}^{T,L}} \min_{S_\b(\cUbk) } d(\pi,\pi_\b) \leq 2\gamma +  (1-\tanh(L)),
$$
implying that 
$$
\lim_{L\uaw\infty}\lim_{\b\uaw\infty} 
\ \  \  \P(\max_{\bar {\cal A}^{T,L}} \min_{S_\b(\cUbk) } d(\pi,\pi_\b) \ \leq \ 3 \gamma )=1. 
$$
Since this is true for every $\gamma>0$, the result follows.

It remains to show the claim made earlier. 
In order to prove this, for every $K\in\N\setminus\{0\}$,
let us consider the box $((2K-1)L,(2K+1)L)\times[-T,T],$
and the rescaled left-most path in $ \bar {\cal A}^{T,L}_\b$
starting at the point $(2K \cdot L,(i\gamma)_\b)$,
where $(i \gamma)_\b$ denotes the point in $\Z e^{-2\b}$ defined
as $([i\gamma e^{2\b}]+1) e^{-2\b}$.
This walk is distributed like a drifted rescaled random walk,
killed at each microscopic step with a probability equal to $k_\b$.
For $\b$ large enough (e.g., such that $e^{-2\b}<\gamma/2$)
there exists a time in the lattice $\Z e^{-2\b}$
and belonging to the interval $[j\gamma,(j+1)\gamma)$. For such $\b$,
there exists a strictly positive probability
for the event $F_\b^{K,\gamma,i,j}$  to occur, where
$
F_\b^{K,\gamma,i,j} 
$
is defined as the event
that our walk stays in the box $((2K-1)L,(2K+1)L)\times[-T,T]$ and is killed in 
the time interval $[j\gamma,(j+1)\gamma]$.
By independence of the elements of
the sequence $\{F_\b^{K,\gamma,i,j}  \}_{K\geq1}$,
there must exist a.s. a $K$ such that the event  $F_\b^{K,\gamma,i,j}$  
occurs. 
Since for $\b$ large enough,
our walk obviously starts from the interval 
$[i\gamma,(i+1)\gamma)$,
this implies the claim made earlier and ends 
the proof 
of our lemma.

\end{proof}

\blem\label{easy-1}
$$
\forall T>0, \ \lim_{L\uaw\infty}\P\left(  {\cal B}^{T,3L} \subset  {\cal A}^{T,4L}, \ \ \ \bar{\cal B}^{T,3L} \subset  \bar {\cal A}^{T,2L} \right) \ = \ 1.
$$
\elem
\begin{proof}
Let us consider $l$, the left-most path in the (non-killed) Brownian net starting from $(-T, 2.5 L)$, and $r$ the right-most
path starting from $(-T,3.5 L)$. 
We also define $\tilde l,\tilde r$ the left-most and right-most path 
starting respectively from
$(-T,-3.5L)$ and $(-T,-2.5L)$.
Finally,
let us consider the event
that 
\beqnn
\min_{[-T,T]} l  > 2 L, \max_{[-T,T]} l  < 3 L, & \mbox{and}&  \min_{[-T,T]} r >3L, \  \max_{[-T,T]} r < 4L  \ \ \mbox{and}\\ 
\max_{[-T,T]} \tilde r  < - 2 L, \min_{[-T,T]} \tilde r  > - 3 L, & \mbox{and} &  \max_{[-T,T]} \tilde l < - 3L, \  \min_{[-T,T]} \tilde l > - 4L.
\eeqnn
Using the fact that there is no path 
of the Brownian net crossing a left-most path (resp., right-most path) from the right (resp., from the 
left)
it is easy to see that under the event above,
we must have
$$
 {\cal B}^{T,3L} \subset  {\cal A}^{T,4L}, \ \ \ \bar{\cal B}^{T,3L} \subset  \bar {\cal A}^{T,2L}. 
$$
For instance, any path starting from the right of the vertical line $\{3L\}\times\R$
can not enter the box $\square_{L,2T}$
since it is blocked by the path $l$. On the other hand, any
path starting in the box $\square_{L,2T}$
is unable to exit the box $\square_{L,4T}$
since it is blocked by the paths $r$ and $\tilde l$.

Since left-most and right-most paths starting from deterministic points
are distributed like drifted Brownian motions, when $L$ goes to $\infty$ for fixed $T$,
the event above is realized with probability 
going to $1$. This ends the proof of
our lemma.

\end{proof}

\bigskip

{\it Part 3.} We end the proof of Theorem \ref{onewebteo} by considering the set 
of paths starting before time $-T$ and killed
on the time interval $[-T,T]$. Let $\pi$ be such a path 
and let us consider $\pi'$,
the section of $\pi$
after time $-T$. This path 
obviously belongs to ${\cal A}^T$  
and by the previous part, it can be approximated 
by a sequence of discrete paths $\{\pi_\b'\}$. The main 
contributions
to the distance between $\pi$ and $\pi_\b'$ are (1) the difference in their starting time
and (2) the variation of $\pi$ on $[\sigma_\pi,-T]$.
However, for $T$ large enough, any point below the line $\R\times\{-T\}$
is close to $-\infty$. Thus, as $T$ goes to $\infty$, the difference vanishes.

\section{List of Main Symbols}

$$
\begin{array}{ccc}
\Z_{even}^2 & \mbox{Even square lattice } &  \mbox{Section}  \ \ref{intro}\\
\b &\mbox{Scaling parameter} & \mbox{Section}  \ \ref{intro} \\
\cUbk & \mbox{Banching-Coalescing-Killing system} &  \mbox{Section} \ \ref{intro} \\
\cUb & \mbox{Banching-Coalescing system} &  \mbox{Section} \ \ref{intro} \\
\sigma_\pi & \mbox{Starting time of the path $\pi$} &  \mbox{Section} \ \ref{intro}\\
e_\pi & \mbox{Ending time of the path $\pi$} & \mbox{Section} \ \ref{intro}  \\
S_\b(x,t) & \mbox{Rescaling of space-time points.} & p.\pageref{scal-transf}\\
d & \mbox{Distance between paths}& p.\pageref{diistance}\\
d_{{\cal H}}&\mbox{Distance between compact sets of paths}& p.\pageref{dh}\\
{\cal W}&\mbox{Brownian web}& p. \pageref{Forward_Web} \\
(m,n) \ \mbox{points} &\mbox{Special points of the Brownian web}& p.\pageref{special_points}\\
{\Net}^b&\mbox{Brownian net}& p.\pageref{standard-net}\\
(\Wl,\Wr)&\mbox{Pair of left-most and right-most Brownian webs.}&  p.\pageref{standard-net}\\
\sim^{in}_z &\mbox{Equivalence relation for paths entering a point.}& p.\pageref{eq:in-out}\\
\sim^{out}_z &\mbox{Equivalence relation for paths exiting a point.}& p.\pageref{eq:in-out}\\
(p,p),(0,pp) \cdots &\mbox {Special points of the Brownian net}&p.\pageref{special-force}\\
{\xi^S}&\mbox{Branching Coalescing (BC) point set (with no killling) starting form $\R$ at time $0$}& p.\pageref{df-age}\\
{\Net^{b,k}}&\mbox{Brownian net with killing}&p.\pageref{onenet}\\
{{\cal M}^{k}}&\mbox{Set of killing points.}& p.\pageref{onenet} \\
{age(x,t)} & \mbox{Age of a point}& p.\pageref{age} \\
{{\cal T}^\delta} & \mbox{Time length measure on the $\delta$-shorteneing of the Brownian net.}& p.\pageref{continuous-tl}\\ 
\theta^L & \mbox{First killing time} &  p.\pageref{thetaL}\\
\theta^{i,L,T} & \mbox{Sequence of  killing times} & p.\pageref{thetail} \\

\end{array}
$$

{\bf Acknowledgments.}
The research of C.M.N. is supported in part by NSF grants OISE-0730136, DMS-1007524 and DMS-1007626.

\end{document}